\documentclass[10pt]{amsart}
\usepackage{amsmath}
\usepackage{xypic}
\usepackage{amssymb}
\usepackage{amsfonts}
\usepackage{amsthm}
\usepackage{enumerate}
\usepackage[all]{xy}
\usepackage[latin1]{inputenc}
\usepackage{graphicx}
\usepackage{latexsym}
\usepackage{accents}
\usepackage{mathdots}
\usepackage{mathrsfs}
\usepackage{color}
\input xy
\usepackage{tikz}
\usetikzlibrary{arrows}
\usetikzlibrary{decorations,arrows}
\usetikzlibrary{decorations.markings}
\usetikzlibrary{arrows.meta, bending}
\usetikzlibrary{positioning}

\makeatletter
    \newtheoremstyle{definition}
        {5pt}
        {3pt}
        {}
        {0pt}
        {\scshape}
        {.}
        {5pt}
        {\thmname{#1} \thmnumber{#2} \thmnote{[#3]}} 

\newtheoremstyle{theorems}
        {5pt}
        {3pt}
        {\itshape}
        {0pt}
        {\scshape}
        {.}
        {5pt}
        {\thmname{#1} \thmnumber{#2}\thmnote{[#3]}} 

\swapnumbers\theoremstyle{theorems}

\newtheorem{Theo}{Theorem}[section]
\newtheorem{Prop}[Theo]{Proposition}
\newtheorem{Cor}[Theo]{Corollary}
\newtheorem{Lemma}[Theo]{Lemma}
\newtheorem{Prop(BG)}[Theo]{Proposition (Bongartz-Gabriel)}
\newtheorem{Lemma(Asashiba)}[Theo]{Lemma(Asashiba)}
\newtheorem{Lemma(Gab)}[Theo]{Lemma(Gabriel)}
\newtheorem{Theo(Mil)}[Theo]{Theorem (Milicic)}

\theoremstyle{definition}
\newtheorem{Defn}[Theo]{Definition}

\newtheorem{Defn(Asashiba)}[Theo]{Definition (Asashiba)}

\newcommand{\Z}{\mathbb{Z}}
\newcommand{\T}{\mathbb{T}}

\newcommand{\cA}{\mathcal{A}}
\newcommand{\cB}{\mathcal{B}}
\newcommand{\cC}{\mathcal{C}}

\newcommand{\cF}{\mathcal{F}}

\newcommand{\cG}{\mathcal{G}}

\def\Sa{\hbox{${\mathit\Sigma}$}}

\def\La{\hbox{$\it\Lambda$}}

\newcommand{\Hom}{{\rm Hom}}
\newcommand{\GHom}{{\rm GHom}}
\newcommand{\Ext}{{\rm Ext}}
\newcommand{\End}{{\rm End}}

\newcommand{\soc}{{\rm soc \hspace{.4pt}}}
\newcommand{\rad}{{\rm rad \hspace{.4pt}}}

\def\Mod{\hbox{{\rm Mod}{\hskip 0.3pt}}}

\def\proj{\hbox{{\rm gproj}{\hskip 0.5pt}}}
\def\ModLa{\hbox{{\rm Mod}{\hskip 0.5pt}\La}}

\def\GrLa{{\rm GMod}\La}

\def\mf{\mathfrak}

\newcommand{\dt}{{\accentset{\hspace{2pt}\mbox{\large\bfseries .}}{}}}

\newcommand{\cdt}{\dt\hspace{2.5pt}}

\newcommand{\pdt}{{\hspace{.8pt}\dt\hspace{1.5pt}}}
\newcommand{\ydt}{{\hspace{0.5pt}\dt\hspace{1.5pt}}}

\newcommand{\mk}{\mathfrak}

\newcommand{\mla}{\hspace{-.5pt} \langle \hspace{-0.4pt}}
\newcommand{\mra}{ \hspace{-.8pt} \rangle}

\newcommand{\nla}{ \hspace{-1pt} \langle \hspace{-0.4pt}}
\newcommand{\nra}{ \hspace{-.5pt} \rangle}

\newcommand{\sla}{\hspace{-.6pt} \langle \hspace{-.2pt}}
\newcommand{\sra}{\hspace{.3pt}\rangle}

\newcommand{\tla}{\hspace{-.4pt} \langle \hspace{.4pt}}
\newcommand{\tra}{\hspace{-.0pt}\rangle}

\newcommand{\nid}{\hspace{-1.5pt}\mid\hspace{-1.5pt}}

\newcommand{\m}{\hspace{-.5pt}}
\newcommand{\ncirc}{\!\!\circ\m}

\begin{document}

\title[Derived categories]{\sc A representation theoretic perspective of \\ Koszul theory
}

\author[A. Bouhada]{Ales Bouhada}

\author[M. Huang]{Min Huang}

\author[Z. Lin]{Zetao Lin}

\author[S. Liu]{Shiping Liu \vspace{-10pt}
}





\address{Ales Bouhada \\ Williams School of Business, Bishop's University and Champlain College, Lennoxville, Qu\'ebec, Canada.}
\email{abouhada@ubishops.ca}

\address{Min Huang \\ School of Mathematics (Zhuhai), Sun Yat-Sen University, Zhuhai, China.}
\email{huangm97@mail.sysu.edu.cn}

\address{Zetao Lin \\ Department of Mathematical Sciences, Tsinghua University, 100084, Beijing, P. R. China.}
\email{zetao.lin@foxmail.com}

\address{Shiping Liu \\ D\'epartement de math\'ematiques, Universit\'e de Sherbrooke, Sherbrooke, Qu\'ebec, Canada.}
\email{shiping.liu@usherbrooke.ca }

\begin{abstract}

We discover a new connection between Koszul theory and representation theory.
Let $\La$ be a quadratic algebra defined by a locally finite quiver with relations. Firstly, we give a combinatorial description of the local Koszul complexes and the quadratic dual $\La^!$, which enables us to describe the linear projective resolutions and the colinear injective coresolutions of graded simple $\La$-modules in terms of $\La^!$. As applications, we obtain a new class of Koszul algebras and a stronger version of the Extension Conjecture for finite dimensional Koszul algebras with a noetherian Koszul dual. Then we construct two Koszul functors, which induce a $2$-real-parameter family of pairs of derived Koszul functors between categories derived from graded $\La$-modules and those derived from graded $\La^!$-modules. In case $\La$ is Koszul, each pair of derived Koszul functors are mutually quasi-inverse, one of the pairs is Beilinson, Ginzburg and Soergel's Koszul duality. If $\La$ and $\La^!$ are locally bounded on opposite sides, then the Koszul functors induce two equivalences of bounded derived categories: one for finitely piece-supported graded modules, and one for finite dimensional graded modules. And if $\La$ and $\La^!$ are both locally bounded, then the bounded derived cate\-gory of finite dimensional graded $\La$-modules has almost split triangles with the Auslander-Reiten translations and the Serre functors given by composites of derived Koszul functors.

\vspace{-20pt}

\end{abstract}

\maketitle

\section*{Introduction}

The history of Koszul theory traces back to Cartan and Eilenberg's computation of the cohomology groups of a Lie algebra using the Koszul resolution; see \cite[Section 8.7]{CEi}. 
This theory is connected to numerous research domains such as algebraic topology; see \cite{GKM,Pr2}, algebraic geometry; see \cite{BGS},
Hopf algebras and Lie theory; see \cite{BGS,May,MOS,RH}. 
Beilinson, Ginzburg and Soergel described the Koszul duality between a locally finite dimensional Koszul algebra $\La$ and its Koszul dual $\La^!$,
that is a pair of mutually quasi-inverse equivalences between a category derived from
graded $\La$-modules and one derived from graded $\La^!$-modules. In case $\La$ is finite dimensional and $\La^!$ is left noetherian, they obtained an equivalence of the bounded derived categories of finitely generated graded modules. Later, the Koszul duality has been generalized to positively graded Koszul categories; see \cite{MOS}. On the other hand, the representation theory of finite dimensional Koszul algebras has been studied by many representation theorists; see, for example, \cite{GMV,GMV2, GrZ, Mar, RMV1, RMV2, MvS}.


Motivated by the application of the covering technique in representation theory; see \cite{As2, BaL, BoG, Gre}, this paper aims to study Koszul algebras defined by locally finite quivers, from a novel viewpoint of connecting Koszul theory and representation theory. Our contribution is twofold. As to Koszul theory, not only our Koszul algebras have infinitely many graded simple modules, the classical Koszul duality of Beilinson, Ginzburg and Soergel is extended to a 2-real-parameter family of pairs of mutually quasi-inverse equivalences. And under a weaker hypothesis, we obtain two equivalences of bounded derived categories, one for finitely piece-supported graded modules and one for finite dimensional graded modules. In contrast to their sophisticated technique of spectral sequences, our tool is elementary: a local version of the Acyclic Assembly Lemma and an existent technique of functor extension.

As to representation theory, 
we obtain a new class of Koszul algebras and a stronger version of the Extension Conjecture for finite dimensional Koszul algebras with a noetherian Koszul dual. The Koszul functors yield an explicit graded projective resolution and a graded injective co-resolution for every finite dimensional graded modules over a Koszul algebra. This could make finite dimensional Koszul algebras a testing class for other homological conjectures such as the Finitistic Dimensional Conjecture. In the locally bounded Koszul case, we obtain an existence theorem for almost split triangles in the bounded derived category of finite dimensional graded modules, and describe the Auslander-Reiten translations and the Serre functors in terms of the derived Koszul functors. This will stimulate future study on graded Auslander-Reiten components of a hereditary or radical squared zero algebra, as did in the ungraded setting; 
see \cite{BaL2, BLP}.


In order to outline the content section by section, we let $\La$ be a graded algebra defined by a locally finite quiver with homogeneous relations, and denote by $\GrLa$ the category of unitary graded left $\La$-modules, whose subcategories of finitely piece-supported modules, of piecewise finite dimensional modules and of finite dimensional modules are written as ${\rm GMod}^{b\hspace{-3pt}}\La$, ${\rm gmod}\La$ and ${\rm gmod}^{b\hspace{-3pt}}\La$, respectively.

In Section 1, we lay the foundation for this paper. In Section 2, we introduce linear projective $n$-presentations and colinear injective $n$-copresentations; see (\ref{n-presentation}) and prove that $\La$ is quadratic if and only if every graded simple $\La$-module admits a linear projective $2$-presentation; see (\ref{qua-alg}). In particular, a Koszul algebra; see (\ref{Kalg}) is always quadratic; compare \cite[(2.3.3)]{BGS}.

In Section 3, we give a combinatorial description of the local Koszul complexes and the quadratic dual $\La^!$ in case $\La$ is quadratic; see (\ref{R-derivative}) and  (\ref{perp}). This enables us to describe linear projective resolutions and colinear injective coresolutions for graded simple modules in terms of subspaces of $\La^!$; see (\ref{k-cplx-iso}) and (\ref{Koz-cplx-dual}). And we show that $\La$ is Koszul if and only if its quadratic dual or opposite algebra is Koszul if and only if every graded simple module has a colinear injective coresolution; see (\ref{Opp-Koszul}); compare \cite[(2.2.1), (2.9.1)]{BGS}. As applications, we obtain a sufficient condition for a quadratic special multi-serial algebra to be Koszul; see (\ref{SKA}) and establish a stronger version of the Extension Conjecture; see, for definition, \cite[(2.6)]{ILP} for finite dimensional Koszul algebras with a noetherian Koszul dual; see (\ref{Ex-Conj-fd}).

In Section 4, we first develop a homotopy theory for double complexes over a concrete additive category $\cA$, including a local version of the Acyclic Assembly Lemma; see (\ref{Homology-zero}); compare \cite[(2.7.3)]{Wei}. Then, we formalize a technique of extending a functor from $\cA$ into the complex category $C(\cB)$ of a concrete additive category $\cB$ to a functor from $C(\cA)$ into $C(\cB)$; see (\ref{F-extension}). Then, we show that the extended functor descends to
categories derived from suitable subcategories of $C(\cA)$; see (\ref{Der-functor}). This is the key ingredient for constructing the derived Koszul functors.

In Section 5, we describe the generalized Koszul dualities. First in case $\La$ is quadratic, we construct two Koszul functors from $\GrLa$ to $C(\GrLa^!)$; see (\ref{F-property}). They extend to two complex Koszul functors from $C(\GrLa)$ to $C(\GrLa^!)$, which induce a 2-real-parameter family of pairs of derived Koszul functors between categories derived from 
subcategories of $C(\ModLa)$
and those derived from 
subcate\-gories of $C(\ModLa^!)$; see (\ref{F-diag}).
In case $\La$ is Koszul, 
derived Koszul functors in each pair are mutually quasi-inverse; see (\ref{Main}), and one of the pairs is the classical Koszul duality; see \cite[(2.12.1)]{BGS} and \cite[Theorem 30]{MOS}.
In case $\La$ and $\La^!$ are locally bounded on opposite sides, the Koszul functors induce two triangle-equivalences $D^{\hspace{.4pt} b\hspace{-.5pt}}({\rm GMod}^{\hspace{.5pt} b\hspace{-2.8pt}}\La) \cong D^{\hspace{.4pt}b\hspace{-.5pt}}({\rm GMod}^{\hspace{.5pt} b\hspace{-2.8pt}}\La^!)$ and $D^{\hspace{.4pt}b\hspace{-.5pt}}({\rm gmod}^{\hspace{.5pt} b\hspace{-2.8pt}}\La) \cong D^{\hspace{.4pt}b\hspace{-.5pt}}({\rm gmod}^{\hspace{.5pt} b\hspace{-2.8pt}}\La^!)$; see (\ref{Main-2}). 


In Section 6, we study almost split triangles in bounded derived categories of graded $\La$-modules in case $\La$ is Koszul. In fact, the indecomposable images of complexes in $C^{\hspace{.5pt}b\hspace{-.6pt}}({\rm gmod}^{\hspace{.5pt}b\hspace{-2.8pt}}\La^!)$ under the complex Koszul functors fit into almost split triangles in $D^{\hspace{.5pt}b\hspace{-.6pt}}({\rm gmod}\La)$; see (\ref{Main-ast-1}). And every derived-indecomposable complex in $C^{\hspace{.5pt}b\hspace{-.6pt}}({\rm gmod}^{\hspace{.5pt}b\hspace{-2.8pt}}\La)$ is the ending (respectively, starting) term of an almost split triangle in $D^{\hspace{.5pt}b\hspace{-.6pt}}({\rm gmod}\La)$ if and only if $\La^!$ is locally right (respectively, left) bounded; see (\ref{Main-ast-2}). In case $\La$ is locally bounded, $D^{\hspace{.5pt}b\hspace{-.6pt}}({\rm gmod}^{\hspace{.4pt}b\hspace{-2.8pt}}\La)$ has almost split triangles on the right (respectively, left) if and only if $\La^!$ is locally right (respectively, left) bounded. In case $\hspace{-1.5pt}\La\hspace{-1.5pt}$ and $\hspace{-1.5pt}\La^!\hspace{-1.5pt}$ are both locally bounded, the Auslander-Reiten translations and the Serre functors for $D^{\hspace{.5pt}b\hspace{-.6pt}}({\rm gmod}^{\hspace{.4pt}b\hspace{-2.8pt}}\La)$ are composites of derived Koszul functors; see (\ref{Main-ast-3}).

\vspace{-0pt}

\section{Preliminaries}

\vspace{0pt}

The objective of this section is to fix some terminology and notation, which will be used throughout the paper, and collect some preliminary results.

\vspace{2pt}

\noindent{\sc 1) Linear algebra.} Throughout this paper, $k$ denotes a commutative field, and all tensor products will be over $k$. Given a set $\mathcal S$, the $k$-vector space spanned by $\mathcal S$ will be written as $k\mathcal S$. We write ${\rm Mod}\hspace{.4pt}k$ for the category of $k$-vector spaces and ${\rm mod}\hspace{.4pt}k$ for the category of finite dimensional $k$-vector spaces. We shall make a frequent use of the exact functor $D=\Hom_k(-, k): {\rm Mod}\hspace{.4pt}k\to {\rm Mod}\hspace{.4pt}k$, which restricts to a duality $D: {\rm mod}\hspace{.4pt}k\to {\rm mod}\hspace{.4pt}k$.
%
%
%
%
%
%
%
The following statement is well-known.

\vspace{-1pt}


\begin{Lemma}\label{Tensor}

Given any $k$-vector spaces $U, V; M, N$, we have a $k$-linear map \vspace{-2pt}
$$\rho: \Hom_k(U, M) \otimes \Hom_k(V, N)\to \Hom_k(U\otimes V, M\otimes N): f\otimes g\mapsto \rho(f\otimes g),\vspace{-2pt}$$ which is natural in all variables, such that $\rho(f\otimes g)(u\otimes v)=f(u) \otimes g(v)$ for $u\in U$ and $v\in V$. And $\rho$ is an isomorphism in case $U, M\in {\rm mod} k$ 
or $V, N\in {\rm mod} k$. 

\end{Lemma}

\noindent{\sc Remark.} We shall identify $f\otimes g$ with $\rho(f\otimes g)$ in case $\rho$ is an isomorphism.

\vspace{3pt}

Since $V \m\otimes k \cong\! V \m\cong\m \Hom_k(k, V)$, we immediately obtain the following consequence.

\vspace{-1pt}

\begin{Cor}\label{Cor 1.2}

Given $U\in {\rm mod}\hspace{.4pt}k$ and $V\in {\rm Mod}\hspace{.4pt}k$, we have

\begin{enumerate}[$(1)$]

\vspace{-2pt}

\item a natural $k$-linear isomorphism $\sigma: DU \otimes V \to \Hom_k(U, V): f\otimes v \mapsto \sigma(f\otimes v),$ \vspace{.5pt} where $\sigma(f\otimes v)(u)=f(u) v,$ for $u\in U$ and $v\in V;$

\item a natural $k$-linear isomorphism $\varphi: DV \otimes DU\to D(V\otimes U): f\otimes g\mapsto \varphi(f\otimes g),$ \vspace{.5pt} where $\varphi(f\otimes g)(v\otimes u)= f(v)g(u)$, for $u\in U$ and $v\in V$.

\end{enumerate}\end{Cor}


The following statement will be needed for our later investigation.

\begin{Lemma}\label{dual-mix}

Let $f:U\to M$ and $g: N\to V$ be morphisms in ${\rm mod}\hspace{.5pt}k$. Then, we obtain a commutative diagram with
vertical isomorphisms 
\vspace{-3.5pt}
$$\xymatrixrowsep{18pt}\xymatrixcolsep{38pt}\xymatrix{
\hspace{13pt}U\otimes DV \ar[d]_{\theta_{U, V}} \ar[r]^-{f\otimes Dg} &  M\otimes DN \ar@<-2ex>[d]^{\theta_{M, N}} \hspace{18pt}\\
D(V\otimes DU)\ar[r]^-{D(g\otimes Df)} &  D(N\otimes DM).}
\vspace{-5pt}$$
\end{Lemma}

\noindent{\it Proof.} Composing 
the canonical $k$-linear isomorphism $U\otimes DV \to D^2U\otimes DV$ with the $k$-linear isomorphism $D^2U\otimes DV\to D(V\otimes DU)$ as stated in Corollary \ref{Cor 1.2}(2), we obtain a $k$-linear isomorphism $\theta_{U, V}:U\otimes DV \to D(V\otimes DU)$ such that $\theta_{U, V}(u\otimes \zeta)(v\otimes \xi)=\zeta(v) \xi(u)$, for $u\in U;$ $v\in V;$ $\zeta\in DV$ and $\xi\in DU$.
Now, it is routine to verify the commutativity of the diagram stated in the lemma commutes.
The proof of the lemma is completed.

%
%
%
%
%
%
%

\vspace{3pt}

\noindent{\sc 2) Quivers.} Throughout this paper, $Q=(Q_0, Q_1)$ denotes a locally finite quiver, where $Q_0$ is a set of vertices and $Q_1$ is a set of arrows such that at most finitely many arrows start or end at any given vertex. 
Given an arrow $\alpha: x\to y$, we write $x=s(\alpha)$ and $y=e(\alpha)$. For each $x\in Q_0$, one associates a {\it trivial path} $\varepsilon_x$ with $s(\varepsilon_x)=e(\varepsilon_x)=x$. A {\it path} of positive length $n$ is a sequence $\rho=\alpha_n \cdots \alpha_1$, with $\alpha_i\in Q_1$, such that $s(\alpha_{i+1})=e(\alpha_i)$, for $i=1, \ldots, n-1$. In this case, we call $\alpha_1$ the {\it initial arrow}{\hspace{.5pt}}; and
$\alpha_n$, the {\it terminal arrow}, of $\rho$.
Fix $x, y\in Q_0$ and an integer $n\ge 0$. We shall denote by $Q(x, y)$ the set of paths from $x$ to $y$ in $Q$ and by $Q_n$ the set of all paths of length $n$ in $Q$. Moreover, put $Q_n(x, y)=Q_n\cap Q(x, y)$ and write $Q_n(x, -)=\cup_{z\in Q_0} Q_n(x, z)$ and $Q_n(-, x)=\cup_{z\in Q_0} Q_n(z, x)$. 


The {\it opposite quiver} $Q^{\rm o}$ of $Q$ is a quiver defined in such a way that $(Q^{\rm o})_0=Q_0$ and $(Q^{\rm o})_1=\{\alpha^{\rm o}: y\to x \mid \alpha: x\to y \in Q_1\}$. A non-trivial path $\rho=\alpha_n\cdots \alpha_1$ in $Q(x, y)$, where $\alpha_i\in Q_1$, corresponds to a non-trivial path $\rho^{\rm o}=\alpha_1^{\rm o} \cdots \alpha_n^{\rm o}$ in $Q^{\rm o}(y,x)$. For convenience, we identify the trivial path at a vertex $x$ in $Q$ with the trivial path at $x$ in $Q^{\rm o}$.


\vspace{3pt}

\noindent{\sc 3) Algebras defined by quivers with relations.} In this paper, an algebra does not necessarily have an identity, and an ideal in an algebra is always two-sided.

Let $Q=(Q_0, Q_1)$ be a locally finite quiver. Write $kQ$ for the path algebra of $Q$ over $k$. An ideal in $kQ$ is called a {\it relation-ideal} if it is contained in $(kQ^+)^2$, where $kQ^+$ is the ideal in $kQ$ generated by the arrows. A non-zero element $\rho\in kQ(x, y)$ with $x, y\in Q_0$ is called {\it monomial} if $\rho=\lambda p$ with $\lambda\in k$ and $p\in Q(x,y)$, and {\it polynomial} otherwise. An element in $kQ_n$ with $n> 0$ is called {\it homogeneous}; and {\it quadratic} if $n=2$.


Let $\La=kQ/R$, where $R$ is a relation-ideal in $kQ$. Given $x, y\in Q_0$ and $n\ge 0$, we write 
$R(x, y)=R\cap kQ(x,y)$ and $R_n(x, y)=R\cap kQ_n(x,y)$, and put $R_n(x, -)=\cup_{z\in Q_0} R_n(x, z)$ and $R_n(-, y)=\cup_{z\in Q_0} R_n(z, y)$. Furthermore, set $R(x, -)=\cup_{n\ge 0} R_n(x, -)$ and $R(-, y)=\cup_{n \ge 0} R_n(-, y)$. An element $\rho=\sum_{i=1}^s \lambda_ip_i$ in $R(x, y)$ is called a {\it relation} for $\La$ if the $p_i$ are pairwise distinct paths such that $\sum_{i\in \it\Sigma}\lambda_ip_i\not\in R$ for any $\emptyset\ne \Sa \subsetneq
\{1, \ldots, s\}$; and in this case, the $\lambda_ip_i$ are called the {\it summands} of the relation $\rho$.

The algebra $\La$ is called {\it graded} or {\it quadratic} if $R$ is generated by some homogeneous relations or by some quadratic relations, respectively. Moreover, we say that $\La$ is {\it locally left bounded} if the $\La e_x$ are all finite dimensional; {\it locally right bounded} if the $e_x \La $ are all finite dimensional, and {\it locally bounded} if it is locally left and right bounded; compare \cite[(2.1)]{BoG}. Furthermore, $\La$ is called {\it special multi-serial} provided, for any arrow $\alpha$ in $Q_1$, that there exists at most one arrow $\beta$ in $Q$ such that $\beta \alpha \notin R$, and at most one arrow $\gamma$ such that $\alpha \gamma \notin R$; see \cite{VHW}.



Finally, let us fix some notation for $\La$. Write $\bar{\hspace{-1pt}\gamma}=\gamma+R\in \La$ for $\gamma\in kQ$, and $e_x = \bar{\varepsilon}_x$ for $x\in Q_0$. Then,
$\{e_x \nid x\in Q_0\}$ is a complete orthogonal set of idempotents in $\La$. The opposite algebra of $\La$ is given by $\La^{\rm o}=kQ^{\rm o}/R^{\rm o},$ where $R^{\rm o}= \{\rho^{\rm o} \mid \rho\in R\}$. We write $\bar{\gamma\hspace{1.8pt}}^{\hspace{-1pt}\rm o}=\gamma^{\rm o}+R^{\rm o}$ for $\gamma\in kQ$, but $e_x=\varepsilon_x+R^{\rm o}$ for $x\in Q_0$. 

\vspace{2pt}

\noindent{\sc 4) Almost split triangles.} In this paper, all categories are additive $k$-categories, in which morphisms are composed from the right to the left. A full subcategory of a category is called {\it strictly full} if it is closed under isomorphisms. All functors between additive $k$-categories are $k$-linear.
Let $\cA$ be a triangulated $k$-category with translation functor $[1]$. An {\it almost split triangle} in $\cA$ is an exact triangle \vspace{-6pt}
$$\xymatrix{X \ar[r]^f &Y \ar[r]^g & Z \ar[r]^-\delta &X[1]} \vspace{-5pt}$$
with $f$ minimal left almost split and $g$ minimal right almost split; see \cite{Ha1}. In this case, $X$ is called the {\it starting term}, and $Z$ the {\it ending term}, of the almost splt triangle, and we write $X=\tau Z$ and $Z=\tau^-\!X$.

An object in $\cA$ is called {\it strongly indecomposable} if it has a local endomorphism algebra.  We say that $\cA$ {\it has almost split triangles on the right} (respectively, {\it left}) if every strongly indecomposable object in $\cA$ is the ending (respectively, starting) term of an almost split triangle; and in this case, $\tau$ (respectively,
$\tau^-$) is called the {\it right} (respectively, {\it left}) {\it Auslander-Reiten translation}. And we say that $\cA$ {\it has almost split triangles} if it has almost split triangles on both sides.


Let now $\cA$ be Hom-finite and Krull-Schmidt. A functor $\mathbb{S}: \cA\to \cA$ is called a {\it left} (respectively, {\it right}) {\it Serre functor} if there exists a binatural $k$-linear isomorphism $\Hom_{\cA}(X, Y) \cong D\Hom_{\cA}(\mathbb{S}Y, X)$ (respectively, $\Hom_{\cA}(X, Y)\cong D\Hom_{\cA}(Y, \mathbb{S}X))$ for any $X, Y\in \cA$; see \cite[(I.1)]{RVDB}. As shown by Reiten and Van den Bergh,
$\cA$ has almost split triangles on the right (respectively, left) if and only if it admits a right $($respectively, left$)$ Serre functor $\mathbb{S}$; and in this case, $\tau X=\mathbb{S}(X)[-1]$ (respectively, $\tau^-X=\mathbb{S}(X)[1]$) for any indecomposable object $X\in \cA$. Moreover, $\cA$ has almost split triangles if and only if it admits a right Serre equivalence, or equivalently, a left Serre equivalence; see \cite[(I.2.3)]{RVDB}.

\vspace{2pt}

\noindent{\sc 6) Derived categories.} Let $\cA$ be a strictly full additive subcategory of an abelian $k$-category $\mk A$. We denote by $C(\cA)$ the complex category of $\cA$, whose full subcategories of bounded complexes, of bounded-below complexes and of bounded-above complexes are written as $C^b(\cA)$, $C^+(\cA)$ and $C^-(\cA)$ respectively. 
Given $* \!\in\! \{\emptyset, b, +, -\}$, we denote by $K^*(\cA)$ the $(*)$-{\it homotopy category} and by $D^*(\cA)$ the $(*)$-{\it derived category} of $\cA$.

\vspace{0pt}

A strictly full additive subcategory $\mathscr{A}$ of $C(\mk A)$ is called {\it derivable} if it is closed under shifts. In this case, $\mathscr{A}$ is closed under mapping cones. Thus, the quotient category $\mathcal{K}(\mathscr{A})$ of $\mathscr{A}$ modulo the null-homotopic morphisms is a full triangulated subcategory of $K(\mk A)$; see \cite[(II.1.7)]{Mil}, in which the quasi-isomorphisms form a localizing class compatible with the triangulation; see \cite[(III.3.1.2)]{Mil}. Therefore, the localization $\mathcal{D}(\mathscr{A})$ of $\mathcal{K}(\mathscr{A})$ at quasi-isomorphisms is a triangulated category; see \cite[(II.1.6.1)]{Mil}, which we call the {\it cate\-gory derived} from $\mathscr{A}$. 

\vspace{-.5pt}

\section{Linear projective $n$-presentations and quadratic algebras}
The main objective of this section to introduce linear projective $n$-presentations and colinear injective $n$-copresentations for graded modules and show that a graded algebra is quadratic if and only if every graded simple module admits a linear projective $2$-presentation. 

We start with some preliminaries on graded modules. Throughout this section let $\La=kQ/R$ be a graded algebra, where $Q$ is a locally finite quiver and $R$ is a homogeneous relation-ideal in $kQ$. Then $\La$ is a positively graded $k$-algebra with grading $\La=\oplus_{i\ge 0}\La_i$, where $\La_i=\{\bar{\gamma\hspace{1.6pt}} \mid \gamma\in kQ_i\}$.
Moreover, $\La^{\rm o}$ is also positively graded with grading $\La^{\rm o}=\oplus_{i\ge 0} \La_i^{\rm o}$, where $\La_i^{\rm o}=\{\bar{\gamma\hspace{1.8pt}}^{\hspace{-1pt}\rm o} \mid \gamma\in kQ_i\}$.

\vspace{1pt}

A left  $\La$-module $M$ is called {\it unitary} if $M=\oplus_{x\in Q_0} e_xM$ and {\it graded} provided $M=\oplus_{i\in \Z}M_i$,
where the $M_i$ are $k$-vector spaces such that $\La_i M_j \subseteq M_{i+j}$ for $i, j\in \Z$. Let $M$ be a unitary graded left $\La$-module. Then $M=\oplus_{(i,x)\in \Z \times Q_0} M_i(x)$, where $M_i(x)=e_xM_i$ is called the $(i, x)$-{\it piece} of $M$. The elements in $M_i$ are called {\it homogeneous} of degree $i$ and those in $M_i(x)$ are called {\it pure}.

Let $N$ be also a unitary graded left $\La$-module. A $\La$-linear morphism $f: M\to N$ is {\it graded} if $f(M_i)\subseteq N_i$ for all $i\in \Z$. In this case, we write $f_i: M_i\to N_i$ and $f_{i,x}: M_i(x)\to N_i(x)$ for the $k$-linear maps obtained by restricting $f$. In the sequel, we shall identify $f$ with a family of $k$-linear maps $f_{i,x}: M_i(x)\to g_i(x)$ with $(i, x)\in \Z\times Q_0$ such that $uf_{i,x}(m)\m=\m f_{i+j, y}(um)$, for all $u\in e_y\La_j e_x$ and $m\in M_i(x)$.

\vspace{1pt}

The unitary graded left $\La$-modules together with the graded $\La$-linear morphisms form an abelian $k$-category $\GrLa$, in which the morphism spaces will be written as $\GHom_{\mathit\Lambda}(M, N)$. A module $M\in \GrLa$ is called {\it bounded-below} if $M_i=0$ for $i\ll 0$; {\it bounded-above} if $M_i=0$ for $i\gg 0$; {\it finitely piece-supported} if $M_i(x)=0$ for all but finitely many $(i, x)\in \Z\times Q_0$; and {\it piecewise finite dimensional} if $M_i(x)$ is finite dimensional for all $(i, x)\in\mathbb{Z}\times Q_0$. We write ${\rm GMod}^+\hspace{-3.6pt}\La$, ${\rm GMod}^{-}\hspace{-3pt}\La$, ${\rm GMod}^{\hspace{.5pt}b\hspace{-2.8pt}}\La$, ${\rm gmod}\,\La$ and ${\rm gmod}^{\hspace{.5pt}b\hspace{-2.8pt}}\La$ for the subcategories of $\GrLa$ of bounded-below modules, of bounded-above modules, of finitely piece-supported modules, of piecewise finite dimensional modules and of finite dimensional modules, respectively.


Let $M\in \GrLa$. Given $s\in \Z$, the {\it grading $s$-shift} $M\sla s\sra$ of $M$ is defined by $M\sla s\sra_{i}=M_{s+i}$ for all $i\in \Z$. In particular, $M\sla s\sra_i(x)=M_{s+i}(x)$, for $(i, x)\in \Z\times Q_0$.
Given a morphism $f: M\to N$ in $\GrLa$, the {\it grading $s$-shift} $f\nla s\nra$ of $f$ is defined by $f\nla s\nra_i=f_{i+s}: M_{i+s}\to N_{i+s}$ for all $i\in \Z$. \vspace{.5pt} Moreover, given $V\in \Mod k$, we have $M\otimes V=\oplus_{i\in \mathbb{Z}} (M_i\otimes V)\in \GrLa$. Note that $(M\otimes V)_i(x)=M_i(x)\otimes V$ for all $(i,x)\in \Z\times Q_0$. It is evident that $(M\otimes V)\langle s \rangle = M\langle s \rangle \otimes V$.


\vspace{1pt}

In the study of graded $\La$-modules, an important role will be played by a contravariant functor $\mf{D} : \GrLa \to \GrLa^{\rm o}$ introduced in \cite[(2.2)]{LLi}. Given $M \!\in \m \GrLa$, we have $\mf{D}M =
\oplus_{i\in \mathbb{Z}; x\in Q_0} D(M_{-i}(x))$ such that $(u^{\rm o}  \cdot  \varphi)(v)=\varphi(u v)$, for $\varphi\in D(M_{-i}(x))$; $u\in e_x\La_je_y$ and $v\in M_{-i-j}(y)$. \vspace{.5pt} So $(\mf{D} M)_i=\oplus_{x\in Q_0}\, D(M_{-i}(x))$ and $(\mf{D}M)_i(x)=D(M_{-i}(x))$. Given a morphism $f: M\to N$ in $\GrLa$, we have a morphism $\mf{D}(f): \mf{D}N\to \mf{D}M$ such that $\mf{D}(f)_{i,x}=D(f_{-i, x})$, for all $(i, x)\in \Z\times Q_0$. For convenience,  we quote the following statement from \cite[(2.2.1), (2.2.2)]{LLi}.\hspace{-6pt}

\begin{Prop}\label{gr-duality}
Let $\La=kQ/R$ be a graded algebra with $Q$ a locally finite quiver.
\begin{enumerate}[$(1)$]
\vspace{-8pt}

\item If $M\in \GrLa$ and $s\in \Z$, then $\mf{D}(M\nla s\nra)=(\mf{D} M)\nla -s\sra$.

\vspace{.5pt}

\item If $M\in \GrLa$ and $V\in{\rm mod}k$, then $\mf{D}(M\otimes V)\cong  \mf{D}M\otimes DV$.

\vspace{.5pt}

\item The functor $\mf{D}$ restricts to a duality $\mf{D} : {\rm gmod}\La \to {\rm gmod}\La^{\hspace{-.3pt} \rm o}\hspace{-2pt}$.

\end{enumerate} \end{Prop}

\vspace{3pt}

Let $M\in \GrLa$. A $\La$-submodule $N$ of $M$ is {\it graded} if $N=\oplus_{i\in \mathbb Z}(M_i\cap N).$ In this case, if $m=\textstyle \sum_{(i,x)\in \mathbb{Z}\times Q_0}m_{i,x}\in N$ with $m_{i,x}\in M_i(x)$, then \vspace{-.5pt} $m_{i,x}\in N$ for any $(i,x)\in \Z\times Q_0.$ The quotient $M/N$ is graded as $M/N=\oplus_{i\in \Z} (M_i+N)/N$.
Recall that the {\it graded radical} ${\rm rad}M$ of $M$ is the intersection of all graded maximal submodules of $M$, the {\it graded socle} $\soc M$  of $M$ is the sum of all graded simple submodules of $M$, and the {\it graded top} ${\rm top}M$ of $M$ is the quotient $M/\rad M$. 
Note that $\rad(_{\mathit\Lambda}\La)=\oplus_{i\ge 1}\La_i=:\rad \La$ and $\rad M=(\rad \La)M$; see \cite[(2.6.2)]{LLi}, and $\soc M=\{m\in M \mid (\rad \La)m=0\}$; see \cite[(2.9.1)]{LLi}.
Further, one says that $M$ is {\it finitely generated} ({\it  in degree} $s$) if $M=\La m_1+\cdots+\La m_t$, where $m_1, \ldots, m_t$ are homogeneous (of degree $s$);
and {\it finitely cogenerated} ({\it  in degree} $s$) if $\soc M$ is finitely generated (in degree $s$) and graded essential in $M$.

\vspace{0pt}

Let $a\in Q_0$. Put $P_a=\La e_a = \oplus_{i\ge 0} \La_ie_a\in {\rm gmod}\La$, which is generated in degree $0$ by $e_a$. We denote by ${\rm GProj}\La$ and ${\rm gproj}\La$ the strictly full additive subcategories of $\GrLa$ generated by $P_a\mla -s\mra\otimes V$ with $(s, a)\in \Z\times Q_0$ and $V\in {\rm Mod}\hspace{.4pt}k$, and by $P_a\sla -s\mra$ with $(s, a)\in \Z\times Q_0$, respectively. Then, ${\rm GProj}\La$ contains only graded projective modules, and ${\rm gproj}\La$ is generated by all finitely generated graded projective modules in $\GrLa$; see \cite[(2.3.2), (2.12.2)]{LLi}.

On the other hand, we write $P_a^{\hspace{.4pt} \rm o}=\La^{\rm o} e_a\in \proj\La^{\rm o}$ and
$I_a=\mf{D}P^{\hspace{.4pt}\rm o}_a\in {\rm gmod}\La$. Then, $I_a=\oplus_{i\ge 0} (\m I_a\m)_{-i}$ with
$(I_a)_{-i}=\oplus_{x\in Q_0} (I_a)_{-i}(x)$ and $(I_a)_{-i}(x)=D(e_x \La_{i\hspace{.3pt}}^{\rm o} e_a)$, for $i\ge 0$ and $x\in Q_0$. Note that $I_a$ is co-generated in degree $0$ by $e_a^\star$, where $e_a^*\in D(e_0\La^{\rm o}e_a)$ such that $e_a^*(e_a)=1$. We denote by ${\rm GInj}\hspace{.4pt}\La$ and ${\rm ginj}\hspace{.4pt}\La$ the strictly full additive subcategories of $\GrLa$ generated by $I_a\tla s\tra\otimes V$ with $(s, a)\in  \Z \times Q_0$ and $V\in {\rm Mod}\,k$ and by $I_a\sla s\sra$ with $s\in  \Z$, respectively. Then, ${\rm GInj}\hspace{.4pt}\La$ contains only graded injective modules, and ${\rm ginj}\hspace{.4pt}\La$ is generated by all finitely cogenerated graded injective modules in $\GrLa$; see \cite[(2.4.2), (2.12.2)]{LLi}.

Let $M\in {\rm gmod}\La$. A {\it graded projective cover} over ${\rm gproj}\La$ for $M$ is an epimorphism $f:  P\to M$ with $P \in {\rm gproj}\La$ such that ${\rm Ker}(f)\subseteq \rad P$; and a {\it graded injective envelope} over ${\rm ginj}\La$ is monomorphism $g\!:\! M\m\to\m I$ with $I\!\in \m {\rm ginj}\La$ and $\soc I \subseteq {\rm Im}(g)$.
The following definition; see \cite[(2.7.1), (2.10.1)]{LLi} is important for constructing graded projective covers and graded injective envelopes.

\begin{Defn}
Let $\La=kQ/R$ be a graded algebra with $Q$ a locally finite quiver. A set $\{m_1, \ldots, m_r\}$ of pure elements in a module $M\in \GrLa$ is called a {\it top-basis} for $M$ provided that $M=\La m_1+\cdots+ \La m_r$ and $\{m_1+{\rm rad}M, \ldots, m_r+{\rm rad}M\}$ is a $k$-basis of ${\rm top} M$; and a {\it soc-basis} for $M$ provided that $\soc M$ has $\{m_1, \ldots, m_r\}$ as a $k$-basis and is graded essential in $M$.

\end{Defn}

The following statement follows from the results stated in \cite[(2.8.2), (2.11.2)]{LLi}.

\vspace{-.5pt}

\begin{Prop}\label{p-cover}

Let $\La=kQ/R$ be a graded algebra with $Q$ a locally finite quiver. Given $s\in \Z$, a module $M$ in ${\rm gmod}\La$ has

\begin{enumerate}[$(1)$]

\vspace{-2pt}

\item a graded projective cover $f \m : \m P_{a_{\hspace{-.8pt}1}} \m \sla -s\nra \m\oplus \cdots \oplus\m P_{a_{\m r\m}} \nla -s \nra \m \to \! M,$ sending $e_{a_i}\!$ to $m_i$, if and only if $\{m_1,\dots, m_r\}$ with $m_i\in M_s(a_i)$ is a top-basis for $M;$

\item a graded injective envelope $g\m:\! M \!\to\! I_{a_1} \m\mla s \mra \hspace{.3pt}\oplus \hspace{.3pt} \cdots \hspace{.3pt} \oplus \hspace{.3pt} I_{a_r} \m \mla s\mra$, sending $m_i$ to $e_{\m a_i}^\star\m\m$, if and only if $\{m_1, \dots, m_r\}$ with $m_i \in M_{-s}(a_i)$ is a soc-basis for $M.$

\end{enumerate}
\end{Prop}

\vspace{2pt}

We are ready to introduce the main notions of this section.

\vspace{-.5pt}

\begin{Defn}\label{n-presentation}

Let $\La=kQ/R$ be a graded algebra with $Q$ a locally finite quiver.
Consider a module $M\in {\rm gmod}\La$ and an integer $n\ge 0$.

\begin{enumerate}[$(1)$]

\vspace{-2pt}

\item In case $M$ is finitely generated in degree $s$, an exact sequence \vspace{-5pt}
$$\xymatrixrowsep{15pt}\xymatrixcolsep{22pt}\xymatrix{P^{-n} \ar[r]^{d^{-n}} & P^{1-n} \ar[r] & \cdots \ar[r] & P^{-1} \ar[r]^-{d^{-1}} &  P^0 \ar[r]^-{d^{\hspace{.5pt}0}} & M \ar[r] & 0}\vspace{-3pt}$$ in ${\rm gmod}\La$ is called a {\it linear projective $n$-presentation} of $M$ 
if $P^{-i} \!\m\in \hspace{-.8pt} {\rm gproj}\La$ is generated in degree $s+i$, for $i=0, \ldots, n$, and ${\rm Ker}(d^{-n})\subseteq \rad P^{-n}$.

\item In case $M$ is finitely cogenerated in degree $-s$, an exact sequence \vspace{-5pt}
$$\xymatrixrowsep{15pt}\xymatrixcolsep{22pt}\xymatrix{0\ar[r] & M \ar[r]^{d^0} & I^0 \ar[r]^{d^1} &I^1 \ar[r] & \cdots \ar[r] & I^{n-1} \ar[r]^{d^n} & I^n}\vspace{-3pt}$$ in ${\rm gmod}\La$ is called a {\it colinear injective $n$-copresentation} of $M$ 
if $I^{i} \!\m\in \hspace{-.8pt} {\rm ginj}\La$ is cogenerated in degree $-s-i$, for $i=0, \ldots, n$, and $\soc I^n\subseteq {\rm Im}(d^n)$.

\end{enumerate}\end{Defn}

\vspace{2pt}

The following statement is important for our later investigation.

\vspace{-1pt}

\begin{Lemma}\label{n-presentation-dual}

Let $\La=kQ/R$ be a graded algebra with $Q$ a locally finite quiver.
Given a module $M\in {\rm gmod}\La$, a sequence \vspace{-5pt} $$\xymatrixrowsep{15pt}\xymatrixcolsep{22pt}\xymatrix{P^{-n} \ar[r]^-{d^{-n}} & P^{-n+1} \ar[r] & \cdots \ar[r] & P^{-1} \ar[r]^-{d^{-1}} &  P^0 \ar[r]^-{d^{\hspace{.5pt}0}} & M \ar[r] & 0}
\vspace{-3pt}$$ in ${\rm gmod}\La$ is a linear projective $n$-presentation of $M$ if and only if \vspace{-3pt}
$$\xymatrixrowsep{15pt}\xymatrixcolsep{22pt}\xymatrix{0\ar[r] & \mk DM \ar[r]^-{{\mathfrak D} d^{\hspace{.5pt}0}} & \mk DP^0 \ar[r]^-{{\mathfrak D} d^{-1}} & \mk DP^{-1} \ar[r] & \cdots \ar[r] & \mk DP^{1-n} \ar[r]^-{{\mathfrak D} d^{-n}} & \mk DP^{-n}}
\vspace{-2pt}$$ is a colinear injective $n$-copresentation of $\mk DM$.

\end{Lemma}

\noindent{\it Proof.} By Proposition \ref{gr-duality}(2), we have a duality $\mk D: {\rm gmod}\La \to {\rm gmod}\La^{\rm o}$. Thus, one of the two sequences stated in the lemma is exact if and only if the other one is exact.
Moreover, $P^{-i}$ is generated in dgeree $s+i$ if and only if $P^{-i}\cong \oplus_{i=1}^r P_{a_i}\langle -s-i \rangle $ with $a_i\in Q_0$ if and only if $\mk D P^{-i}\cong \oplus_{i=1}^r I_{a_i}\langle s+i \rangle $ if and only if $\mk DP^{-i}$ is cogenerated in degree $-s-i$. Finally, ${\rm Ker}(d^{-n})\subseteq \rad P^{-n}$ if and only if $d^{-n}$ is right minimal if and only if $\mk D d^n$ is left minimal if and only if $\soc (\mk{D}P^{-n}) \subseteq {\rm Im}(\mk{D}d^{-n})$. The proof of the lemma is completed.

\vspace{2pt}

Given $u \in e_a\La_{s-t}e_b$, the right multiplication by $u$ yields a graded $\La$-linear morphism $P[u]: P_a\sla -s \sra\to P_b \hspace{.3pt}\tla -t \tra: v\mapsto vu.$
Note that this notation does not distinguish $P[u]$ from its grading shifts. It is known that every graded simple module in $\GrLa$ is isomorphic to $S_a\langle s\rangle$ for some $(s, a)\in \Z \times Q_0$; see \cite[(2.5.1)]{LLi}. The following statement is the promised characterization of quadratic algebras.

\begin{Theo}\label{qua-alg}

Let $\La=kQ/R$ be a graded algebra with $Q$ a locally finite quiver. Then $\La$ is a quadratic algebra if and only if every graded simple $\La$-module admits a linear projective $2$-presentation.

\end{Theo}

\noindent{\it Proof.} It suffices to show, for any $a\in Q_0$, that $R(a, -)$ is generated by quadratic relations if and only if $S_a$ admits a linear projective $2$-presentation. Write $Q_1(a, -)=\{\alpha_i: a \to b_i \mid i=1, \ldots, r\}$. Then, $\{\bar\alpha_1, \cdots, \bar\alpha_r\}$ is  a top-basis for ${\rm rad}P_a$. Since $\bar\alpha_i \!\in\! ({\rm rad}P_a)_1(b_i)$, by Proposition \ref{p-cover}(1), we obtain a projective presentation \vspace{-7pt}
$$\xymatrixcolsep{21pt}\xymatrixrowsep{18pt}\xymatrix{P_a^{-1}=P_{b_1}\m\mla-1\sra \oplus \cdots \oplus P_{\hspace{.4pt}b_r}\m\mla-1\sra \ar[r]^-{d^{-1}_a} & P_a \ar[r]^{p_a} & S_a \ar[r] & \hspace{-1pt} 0} \vspace{-3pt}$$ of $S_a$ over ${\rm gproj}\La$, where $d^{-1}_a=\left(P[\bar{\alpha}_1], \cdots, P[\bar{\alpha}_r]\right)$ with ${\rm Ker}(d_a^{-1})\subseteq \rad P_a^{-1}$.

Suppose that $R(a, -)$ is generated by quadratic relations. Since $Q$ is locally finite, $R(a,-)$ has a minimal generating set $\{\rho_1, \ldots, \rho_s\}$, where $\rho_j\in kQ_2(a, c_j)$ with $c_j\in Q_0$. Write $\rho_j=\sum_{i=1}^r\gamma_{ij}\alpha_i$, where $\gamma_{ij}\in kQ_1(b_i, c_j)$, for $j=1, \ldots, s$. Considering $\bar\gamma_{ij}\in P_{b_i}\mla-1\sra_2(c_j)$, we have morphisms $P[\bar{\gamma\hspace{1pt}}_{\hspace{-1pt}ij}]\m:\m P_{c_j}\m\nla-2\tra \!\to\! P_{b_i}\m\sla-1\tra$ in ${\rm gproj}\La$, for $1\le i\le r$ and $1\le j\le s$. Putting $d^{-2}_a=(P[\bar{\gamma\hspace{1pt}}_{\hspace{-1pt}ij}])_{r\times s}$, we claim that $S_a$ admits a linear projective $2$-presentation \vspace{-5pt}
$$\xymatrixcolsep{21pt}\xymatrixrowsep{18pt}\xymatrix{P_{c_1} \m\mla - 2\tra\oplus \cdots \oplus P_{c_s}\m\mla - 2 \tra \ar[r]^-{d^{-2}_a} & P_{b_1} \m\mla-1\sra \oplus \cdots \oplus P_{\hspace{.4pt}b_r}\m\mla-1\sra \ar[r]^-{d^{-1}_a} & P_a \ar[r]^{p_a} & S_a \ar[r] & \hspace{-1pt} 0.}
\vspace{-2pt}$$ 

Set \vspace{.5pt} $u_j=d_a^{-2}(e_{c_j})=(\bar{\gamma}_{1j}, \ldots, \bar{\gamma}_{rj})\in {\rm Ker}(d^{-1}_a)_2$, for $j=1, \ldots, s$. Since $\rad P_a^{-1}$ is generated in degree $2$, the $u_1, \ldots, u_s$ are top-elements of ${\rm Ker}(d^{-1}_a).$ Consider $v=(\bar{\delta}_1, \ldots, \bar{\delta}_r)\in {\rm Ker}(d^{-1}_a).$  \vspace{.5pt} To show that $v\in \sum_{j=1}^s \La u_j$, since ${\rm Ker}(d^{-1}_a)$ is graded, we may assume that $\delta_i\in kQ_p(b_i, c)$, where $c\in Q_0$ and $p\ge 1$, for $i=1, \ldots, r$. Then, $\sum_{i=1}^r \delta_i \alpha_i \in R_{p+1}(a, c)$. Since $R(a, -)$ is generated by the $\rho_j$, we may write $$\textstyle{\sum}_{i=1}^r \delta_i \alpha_i={\sum}_{j=1}^s \omega_j \rho_j + {\sum}_{i=1}^r \eta_i \alpha_i ={\sum}_{i=1}^r ({\sum}_{j=1}^s \omega_j\gamma_{ij} + \eta_i)\alpha_i,$$ where $\omega_j\in kQ_{p-1}(c_j, c)$ and $\eta_i\in R_p(b_i, c)$. Therefore, $\delta_i = {\sum}_{j=1}^s \omega_j\gamma_{ij} + \eta_i,$ for $i=1, \ldots, r$. This yields $v={\sum}_{j=1}^s \bar\omega_j u_j$. Further, assume that $\sum_{i=1}^s \lambda_j u_j=0$ with $\lambda_j\in k$. Then $\sum_{j=1}^s\lambda_j \bar \gamma_{ij}=0$. Since $R\subseteq (kQ^+)^2$, we have
$\sum_{j=1}^s\lambda_j \gamma_{ij}=0$, for $i=1, \ldots, r,$ and hence, $\sum_{j=1}^s\lambda_j\rho_j=0$.  By the minimality, $\lambda_1=\cdots=\lambda_s=0.$ This shows that $\{u_1, \ldots, u_s\}$ is a top-basis for ${\rm Ker}(d^{-1}_a)$. By Proposition \ref{p-cover}, $d^{-2}_a$ co-restricts to a graded projective cover of ${\rm Ker}(d_a^{-1})$. This establishes our claim.

Conversely, suppose that $S_a$ admits a linear projective $2$-presentation \vspace{-4pt}
$$\xymatrixcolsep{25pt}\xymatrixrowsep{18pt}\xymatrix{P^{-2} \ar[r]^-{d^{-2}} & P^{-1} \ar[r]^-{d^{-1}} & P^0 \ar[r]^{d^0} & S_a \ar[r] & 0.}\vspace{-4pt}$$
By the uniqueness of graded projective cover, we may assume that $d^0=p_a$ and $d^{-1}\!=\!d_a^{-1}$. Being generated in degree two, $P^{-2}\!=\!P_{c_1}\m\mla-2\nra \oplus \cdots \oplus P_{c_s}\m\mla-2\nra$ where $c_1, \ldots, c_s\in Q_0$.
Then $u_j\!:\,=d^{-2}(e_{c_j}) \m\in\m {\rm Ker}(d_a^{-1})_2(c_j)\subseteq e_{c_j} \m \La_1e_{b_1} \oplus \cdots \oplus e_{c_j} \m \La_1e_{b_r}$, which can be written as $u_j=(\bar\gamma_{1j}, \ldots, \bar \gamma_{rj})$, where $\gamma_{ij}\in kQ_1(b_i, c_j)$. By Proposition \ref{p-cover}(1), $\{u_1, \ldots, u_s\}$ is a top-basis for ${\rm Ker}(d_a^{-1}).$ By the definition of $d_a^{0}$, we see that $\eta_j:=\sum_{i=1}^r \gamma_{ij}\alpha_i \in R_2(a, c_j)$, for $j=1, \ldots, s.$

\vspace{1pt}

Suppose that $\rho\in R_n(a, c)$, for some $n\ge 2$ and $c\in Q_0$. Then, $\rho=\sum_{i=1}^r\, \gamma_i \alpha_i$ for some $\gamma_i\in kQ_{n-1}(b_i,c)$. Observing that $(\bar{\gamma}_1, \ldots, \bar{\gamma}_r)\in {\rm Ker}(d^{-1}_a)$, we may write $(\bar{\gamma}_1, \ldots, \bar{\gamma}_r)=\sum_{j=1}^s \bar\delta_{\hspace{-1pt}j} \hspace{1pt} u_j$, \vspace{1pt} for some $\delta_j\in kQ_{n-2}(c_j, c).$ So, $\gamma_i=\sigma_i+\sum_{j=1}^s \delta_j \gamma_{ij}$, where $\sigma_i\in R_{n-1}(b_i, c),$ for $i=1, \ldots, r$. If $n=2$, since $R_1(-, c)=0$, we have
$\gamma_i=\sum_{j=1}^s \delta_j \gamma_{ij}$, for $i=1, \ldots, r$, and consequently, $\rho=\sum_{j=1}^s \delta_j \eta_j.$ By induction on $n$, we see that $\rho\in \sum_{j=1}^s (kQ)\eta_j.$ This shows that $R(a, -)$ is generated by the quadratic relations $\eta_1, \ldots, \eta_s$. The proof of the proposition is completed.

\vspace{2pt}

More generally, we have the notions of linear projective resolution and co\-linear injective coresolution.

\begin{Defn}\label{linear complex}

Let $\La=kQ/R$ be a graded algebra with $Q$ a locally finite quiver. Given a module $M\in {\rm gmod}\La$, a graded projective resolution \vspace{-5pt}
$$\xymatrixrowsep{15pt}\xymatrixcolsep{20pt}\xymatrix{\cdots
\ar[r] &P^{-n} \ar[r]^{d^{-n}} & P^{1-n} \ar[r] & \cdots \ar[r] & P^{-1} \ar[r]^{d^{-1}} &  P^0 \ar[r] & 0 \ar[r] & \cdots }\vspace{-2pt}$$ of $M$ over ${\rm gproj}\La$ is called {\it linear} if $P^{-n}$ is generated in degree $s+n$ with $s$ a constant, for all $n\ge 0$. And a graded injective coresolution \vspace{-5pt}
$$\xymatrixrowsep{15pt}\xymatrixcolsep{20pt}\xymatrix{\cdots \ar[r] & 0\ar[r] & I^0 \ar[r]^{d^1} &I^1 \ar[r] & \cdots \ar[r] & I^{n-1} \ar[r]^{d^n} &  I^n \ar[r] & \cdots }\vspace{-3pt}$$ of $M$ over ${\rm ginj} \La$ is called {\it colinear} if $I^n$ is cogenerated in degree $t-n$ with $t$ a constant, for all $n\ge 0$.

\end{Defn}



The following definition of a Koszul algebra is essentially the same as the classical one; see \cite[(1.2.1)]{BGS} and \cite[(5.4)]{MOS}.

\begin{Defn}\label{Kalg}

Let $\La=kQ/R$ be a graded algebra with $Q$ a locally finite quiver. We call $\La$ a {\it Koszul algebra} if every graded simple $\La$-module admits a linear projective resolution over ${\rm gproj}\hspace{.5pt}\La$.

\end{Defn}

\noindent{\sc Remark.} (1) By Theorem \ref{qua-alg}, a Koszul algebra is quadratic; compare \cite[(2.3.3)]{BGS}.

\noindent (2) It is clear that $\La$ is Koszul if and only if $S_a$ admits a linear projective resolution over ${\rm gproj}\hspace{.5pt}\La$, for every $a\in Q_0$.

\vspace{2pt}

\noindent{\sc Example.} Given a locally finite quiver $Q$, the path algebra $kQ$ is Koszul. Indeed, every
$S_a$ with $a\in Q_0$ admits a linear projective resolution \vspace{-4pt}
$$\xymatrixcolsep{26pt}\xymatrixrowsep{18pt}\xymatrix{
\cdots \ar[r] & 0 \ar[r] & P_{\hspace{.4pt} b_1}\m \mla-1\sra \oplus \cdots \oplus P_{\hspace{.4pt}b_r} \hspace{-.8pt} \mla-1\sra \ar[rr]^-{\left(P[\bar{\alpha}_1], \cdots, P[\bar{\alpha}_r]\right)} && P_a \ar[r] & 0,}\vspace{-4pt}$$
where $\alpha_i: a \to b_i$, $i=1, \ldots, r$ are the arrows in $Q_1(a, -).$

\section{Local Koszul complexes and Koszul duals}

Throughout this section, let $\La=kQ/R$ be a quadratic algebra, where $Q$ is a locally finite quiver. First, we give a combinatorial description of the local Koszul complexes and the quadratic dual $\La^{\hspace{.3pt}!}$. Then, we describe linear projective resolutions and colinear injective coresolutions of graded simple modules in terms of the subspaces of $\La^{\hspace{.3pt}!}$. This leads to a number of equivalent conditions for $\La$ to be Koszul. As applications, we obtain a new class of Koszul algebras and a stronger version of the Extension Conjecture for certain finite dimensional Koszul algebras.

Let us start with some notation and terminology. Given an arrow $\alpha: y\to x$ in $Q$, we have a $k$-linear {\it derivation} $\partial_\alpha: kQ\to kQ$, sending a path $\rho$ to $\delta$ if $\rho=\alpha \delta$, and to $0$ if $\alpha$ is not a terminal arrow of $\rho$. In particular, $\partial_\alpha$ vanishes on $kQ_0$ and sends $kQ_n$ to $kQ_{n-1}$ for any $n>0$. Fix $a\in Q_0$ and $n>0$.
Given $\alpha\in Q_1(y, x)$, we have a graded $\La$-linear morphism $P[\bar \alpha]: P_x\tla\!-n\tra\to P_y \tla\m 1-n\tra$, the right multiplication by $\bar \alpha$, and a $k$-linear map $\partial_\alpha: kQ_n(a, x)\to kQ_{n-1}(a,y)$. Since $Q$ is locally finite, for any $x, y\in Q_0$, we have a morphism \vspace{-2pt}
$$\partial_a^{-n}(y,x) \!=\! \textstyle\sum_{\alpha\in Q_1(y,x)}\! P[\bar{\alpha}] \m\otimes \m \partial_\alpha\m:\m P_x\sla -n\nra \m \otimes\m kQ_n(a, x) \!\to\! P_y\sla 1\!-\!n\nra \m \otimes\m kQ_{n-1}(a, y)\vspace{-2pt}$$ in ${\rm gproj}\La$. The following statement is useful for later calculation.

\begin{Lemma}\label{diff}

Let $\La=kQ/R$ be a quadratic algebra with $Q$ a locally finite quiver. If 
$u \!\in \!P_x\sla -n \nra,$ $\delta\in kQ_{n-1}(a, y)$ and $\zeta\in kQ_1(y, x)$, then
$\partial_a^{-n}(y,x)(u \m\otimes\m \zeta \delta ) = u \hspace{.5pt} \bar{\zeta} \m \otimes \m \delta.$

\end{Lemma}

\noindent{\it Proof.} Let $\delta\in kQ_{n-1}(a, y)$ and $\zeta\in kQ_1(y, x)$. Then, $\partial_\alpha(\zeta \delta)=\partial_\alpha(\zeta)\,\delta$ for any $\alpha \! \in \! Q_1$.
Write $\zeta=\sum_{\beta \in Q_1(y,x)}\lambda_\beta \beta,$ where $\lambda_\beta\in k.$
For any $u\in P_x\mla -n \nra$, we have
$\partial_a^{-n}(y,x)(u\otimes \zeta \delta)  = \sum_{\alpha,\beta\in Q_1(y,x)} \lambda_\alpha u \hspace{.4pt} \bar \alpha \otimes \partial_{\alpha }(\beta \delta)= \sum_{\alpha\in Q_1(y,x)} \lambda_\alpha u \hspace{.4pt}  \bar \alpha \otimes \delta = u  \hspace{.5pt} \bar\zeta \otimes \delta.\vspace{.5pt}$ The proof of the lemma is completed.


\vspace{3pt}

Fix $a, x\in Q_0$. For $n=0, 1,$ we put $R^{(n)}(a, x)=kQ_n(a,x);$ and for any integer $n\ge 2,$ we define $R^{(n)}(a,x)=\cap_{\hspace{.6pt}0\le j\le n-2}\, kQ_{n-2-j}(-, x) \m\cdot\m R_2 \m\cdot\m kQ_j(a, -).$ In particular, $R^{(2)}(a,x)=R_2(a, x)$. Moreover, set $R^{(n)}(a, -)=\oplus_{x\in Q_0}\, R^{(n)}(a, x)$. 

\begin{Lemma}\label{R-derivative}

Let $R$ be a quadratic ideal of $kQ$, where $Q$ is a locally finite quiver. Consider $a, x\m\in\m Q_0$ with $Q_1(-, x)\!=\!\{\alpha_i\m:\m y_i\m\to\m x \nid i=1, \ldots, r\}$ and an integer $n\ge 1$.

\begin{enumerate}[$(1)$]

\vspace{-1pt}

\item If $\gamma\in R^{(n)}(a, x)$ and $\alpha \in Q_1(y, x)$, then $\partial_\alpha(\gamma)\in R^{(n-1)}(a, y);$ and consequently, $\gamma=\sum_{i=1}^r \alpha_i \gamma_i$, for some $\gamma_i\in R^{(n-1)}(a, y_i)$.

\vspace{1pt}

\item If $\rho=\sum_{i=1}^r \zeta_i \rho_i$ with $\rho_i \in R^{(n-1)}(a, y_i)$ and $\zeta_i\in kQ_1(y_i, x)$, then $\rho \in R^{(n)}(a, x)$ if and only if
$\rho \in R_2(-, x)\cdot kQ_{n-2}(a, -)$.

\end{enumerate}

\end{Lemma}

\noindent{\it Proof.} Let $\gamma\in R^{(n)}(a, x)$ and $\alpha \in Q_1(y, x)$. Since $\partial_\alpha(\gamma)\in kQ_{n-1}(a, y)$, we may assume that $n\ge 3$. For any $0\le j\le n-3,$ \vspace{.5pt} we may write $\gamma \m=\m {\sum}_{i=1}^r \beta_i \zeta_i \rho_i \delta_i,$ where $\beta_i\in Q_1(-, x)$; $\zeta_i\m\in\m kQ_{n-3-j}(-, y_i);$ $\rho_i \m\in\m R_2$; $\delta_i \m\in\m kQ_j(a,-).$ Assume that $\beta_i=\alpha$ if and only if $1\le i\le s$. Then, $\partial_\alpha(\gamma)= \sum_{i=1}^s \zeta_i \rho_i \delta_i\in kQ_{n-3-j}(-, y) \m\cdot\m R_2 \cdot kQ_j(a, -).$
So, $\partial_\alpha(\gamma)\in R^{(n-1)}(a, y).$ The first part of Statement (1) is established, and the second part follows immediately from it. Moreover, Statement (2) follows directly from the definition of $R^{(n)}(a,x)$.
The proof of the lemma is completed.

\vspace{2pt}

Fix $a\in Q_0$. For $n\ge 0$, put $\mathcal{K}_a^{-n}=\oplus_{x\in Q_0} P_x\sla -n\nra\otimes R^{(n)}(a, x)\in {\rm gproj}\La$. Given $n\ge 1$ and $x, y\in Q_0$, by Lemma \ref{R-derivative}(1), we get a graded $\La$-linear morphism \vspace{-3pt}
$$\partial_a^{-n}(y,x)= {\textstyle\sum}_{\alpha\in Q_1(y,x)} P[\bar{\alpha}] \m\otimes \m \partial_\alpha: P_x \sla -n \sra \m\otimes\m R^{(n)}(a, x) \!\to\!
P_y\sla 1\!-\!n\sra \m\otimes\m R^{(n-1)}(a, y).\vspace{-3pt}$$

Write \vspace{1pt} $\mathcal{K}_a^{1-n}=\oplus_{y\in Q_0} P_y\sla 1-n\nra\otimes R^{(n-1)}(a, y)$ and consider the graded $\La$-linear morphism $\partial_a^{-n}=(\partial_a^{-n}(y,x))_{(y,x)\in Q_0\times Q_0}: \mathcal{K}_a^{-n}\to \mathcal{K}_a^{1-n}$. \vspace{-5.5pt} We obtain a sequence
$$\xymatrixcolsep{20pt}\xymatrix{\mathcal{K}_a^\pdt: \hspace{-5pt} & \cdots \ar[r] & \mathcal{K}_a^{-n} \ar[r]^-{\partial_a^{-n}} & \mathcal{K}_a^{1-n} \ar[r] &
\cdots \ar[r] & \mathcal{K}_a^{-1} \ar[r]^{\partial_a^{-1}} & \mathcal{K}_a^0 \ar[r] & 0 \ar[r] & \cdots \vspace{-8pt}}$$
in ${\rm gproj}\La$, which is a complex as shown below, called the {\em local Koszul complex at} $a$; compare \cite[(2.6)]{BGS}.

\begin{Lemma}\label{K-cplx}

Let $\La=kQ/R$ be a quadratic algebra with $Q$ a locally finite quiver. If $a \!\in\! Q_0$, then $\mathcal{K}_a^\cdt\hspace{-3pt}$ is a complex such that ${\rm Ker}(\partial_a^{-n}\m) \!\subseteq\! {\rm rad}\mathcal{K}_a^{-n}\!,$ for all $n>0$.

\end{Lemma}

\noindent{\it Proof.} Fix $n\ge 1$. By definition, $\mathcal{K}_a^{-n}$ is generated in degree $n$. Let $w\!\in\!(\mathcal{K}_a^{-n})_n$ be such that $\partial_a^{-n}(w)=0$. Then, $w\!=\!\sum_{x\in Q_0} e_x \m\otimes\m \gamma_x$ where $\gamma_x\in R^{(n)}(a, x)$. By Lemma \ref{R-derivative}(1), $\gamma_x=\sum_{\beta\in Q_1(-, x)} \beta \delta\m_\beta$ with $\delta\m_\beta\in R^{(n-1)}(a, -)$. And by Lemma \ref{diff},
$\textstyle \partial_a^{-n}(w)=\sum_{x\in Q_0; \beta \in Q_1(-, x)} \bar \beta \otimes \delta\m_\beta=0.$ Since the $\bar \beta$ are
$k$-linearly independent, $\delta\m_\beta=0$ for all $\beta\in Q_1(-, x)$. Thus, $w=0$. 
This shows that ${\rm Ker}(\partial_a^{-n})\subseteq {\rm rad} \mathcal{K}_a^{-n}$.

Next, let $u \in P_{x}\sla -n\mra$ and $\xi \in R^{(n)}(a, x)$. We may assume that $\xi=\rho \hspace{1pt} \delta$, where $\delta \in k Q_{n-2}(a, z)$ and $\rho \in R_{2}(z, x)$. Write $\rho=\sum_{i=1}^{s} \gamma_i\zeta_i$ with $\gamma_i \in kQ_1(y_i, x)$ and $\zeta_i\in kQ_1(z, y_i)$.
By Lemma \ref{diff}, $(\partial_{a}^{1-n} \m\circ \partial_{a}^{-n})(u\otimes \xi)
= u \hspace{.4pt} (\textstyle\sum_{i=1}^{s} \bar \gamma_i\hspace{.4pt} \bar\zeta_i \otimes \delta) = u \hspace{.4pt} \bar \rho \otimes \delta=0.$
%
%
The proof of the lemma is completed.


\vspace{2pt}

Observe that $\mathcal{K}_a^{\hspace{.5pt} 0}=P_a\otimes k\varepsilon_a$, for any $a\in Q_0$. Thus, we have a graded projective cover $\partial_a^{\hspace{.5pt} 0}: \mathcal{K}_a^{\hspace{.5pt} 0} \to S_a$, sending $e_a\otimes \varepsilon_a$ to $e_a+{\rm rad}P_a$.

\begin{Prop}\label{Koz-pres}

Let $\La=kQ/R$ be a quadratic algebra with $Q$ a locally finite quiver. If $a\in Q_0$ and $n>0$, then $S_a$ has a linear projective $n$-presentation over ${\rm gproj}\La$ if and only if the following sequence is exact$\hspace{.8pt}:$ \vspace{-5pt}
$$\xymatrixrowsep{20pt}\xymatrixcolsep{22pt}\xymatrix{\mathcal{K}_a^{-n} \ar[r]^-{\partial_a^{-n}} & \mathcal{K}_a^{1-n} \ar[r] & \cdots \ar[r] & \mathcal{K}_a^{-1} \ar[r]^{\partial_a^{-1}} & \mathcal{K}_a^{\hspace{.5pt} 0} \ar[r]^{\partial_a^{\hspace{.5pt} 0}} & S_a \ar[r] & 0.} \vspace{-3pt}$$

\end{Prop}

\noindent{\it Proof.} Since $\mathcal{K}_a^{-n}$ is generated in degree $n$, it suffices to show the necessity. Let \vspace{-4pt}
$$\xymatrix{P^{-n} \ar[r]^-{d^{-n}} & P^{1-n} \ar[r] & \cdots \ar[r] & P^{-1} \ar[r]^-{d^{-1}} & P^0 \ar[r]^-{d^{\hspace{.5pt} 0}} & S_a \ar[r] & 0}\vspace{-4pt}$$ be a linear projective $n$-presentation of $S_a$ over ${\rm gproj}\hspace{.5pt}\La$. Then, there exists a graded $\La$-linear isomorphism $f^0: P^0\to \mathcal{K}_a^{\hspace{.5pt} 0}$ such that $d^{\hspace{.5pt} 0}=f^0 \m\circ \partial_a^{\hspace{.5pt} 0}$. Assume that we have a commutative diagram with vertical isomorphisms \vspace{-3pt}
$$\xymatrixrowsep{22pt}\xymatrixcolsep{22pt}\xymatrix{ P^{-p} \ar[r]^-{d^{-p}} & P^{1-p} \ar[r]^-{d^{1-p}} \ar[d]_{f^{1-p}} & \cdots \ar[r] &  P^0 \ar[d]^{f^{0}}\ar[r] \ar[r]^-{d^0} & S_a \ar[r] \ar@{=}[d] &  0 \\ \mathcal{K}_a^{-p} \ar[r]^-{\partial_a^{-p}} & \mathcal{K}_a^{1-p} \ar[r]^-{\partial_a^{1-p}} & \cdots \ar[r] & \mathcal{K}_a^0 \ar[r]^{\partial_a^0} & S_a \ar[r] & 0} \vspace{-4pt}$$
for some $1\le p\le n$. Then, $f^{1-p}\circ d^{-p}$ co-restricts to a graded projective cover of ${\rm Ker}(\partial_a^{1-p})$. We shall obtain an isomorphism $f^{-p}\m:\m P^{-p}\!\to\! \mathcal{K}_a^{-p}$ in ${\rm gproj}\La$ such that $f^{1-p} \circ d^{-p}= \partial^{-p}_a \circ f^{-p}.$ This amounts to show that $\partial_a^{-p}$ co-restricts to a graded projective cover of ${\rm Ker}(\partial_a^{1-p})$, that is,
${\rm Ker}(\partial_a^{1-p})\subseteq {\rm Im}(\partial_a^{-p})$ by Lemma \ref{K-cplx}.

Since $P^{-p}$ is generated in degree $p$, by Proposition \ref{p-cover}(1), ${\rm Ker}(\partial_a^{1-p})$ has a top-basis $T^{p-1}$ contained in $(\mathcal{K}_a^{1-p})_{p}$. Choose a $k$-basis $\{\rho_1, \ldots, \rho_t\}$ of $R^{(p-1)}(a, -)$, where $\rho_j\in R^{(p-1)}(a, y_j)$. Then, $\mathcal{K}_a^{1-p}=\oplus_{j=1}^t P_{y_j} \m \sla 1-p \sra\otimes k\rho_j$.
Consider $u \!\in\! T^{p-1}$. \vspace{-1pt} Then, $u\in (\mathcal{K}_a^{1-p})_{p}(z)\!=\!\oplus_{j=1}^t e_z\La_1 e_{y_j} \otimes k\rho_j$, for some $z\!\in \! Q_0$. Write $u=\sum_{j=1}^t u_j$, where $u_j=\bar{\gamma\hspace{1.2pt}}_{\hspace{-1pt}j} \otimes \rho_j$ with $\gamma_j \in kQ_1(y_j, z),$ for $j=1,\ldots, t.$ Now, choose a $k$-basis $\{\xi_1, \ldots, \xi_s\}$ of $R^{(n-1)}(a, -)$, where $\xi_i\in R^{(n-1)}(a, x_i)$. By Lemma \ref{R-derivative}(2), $\rho_j={\sum}_{i=1}^s \zeta_{ij} \, \xi_i,$ where $\zeta_{ij} \in kQ_1(x_i, y_j).$ Thus, $u=\sum_{i=1}^s \sum_{j=1}^t\bar\gamma_j\otimes \zeta_{ij} \xi_i$. This yields $\partial_a^{-n}(u)=\textstyle\sum_{i=1}^s (\textstyle\sum_{j=1}^t \bar\gamma_j\bar \zeta_{ij}) \otimes \xi_i=0.$  As a consequence,  ${\sum}_{j=1}^t \bar{\gamma}_j \bar{\zeta}_{ij} =0$, and hence, $\eta_i={\sum}_{j=1}^t\gamma_j \zeta_{ij} \in R_2(x_i, z)$, for $i=1, \ldots, s$. Put $\omega = \sum_{j=1}^t \gamma_j \rho_j$, where $\gamma_j \in kQ_1(y_j, z)$ and $\rho_j\in R^{(p-1)}(a, y_j).$ Since $\omega=\sum_{i=1}^s\eta_i \hspace{.4pt} \xi_i$, where $\eta_i\in R_2(x_i, z)$ and $\xi_i\in kQ_{n-1}(a, x_i)$, by Lemma \ref{R-derivative}(2), $\omega\in R^{(p)}(a, z).$ So, $e_z\otimes \omega \in \mathcal{K}_a^{-p}$ and $\partial^{-p}_a(e_z\otimes \omega)= \partial^{-p}_a ({\textstyle\sum}_{j=1}^te_z\otimes \gamma_j\rho_j)= {\textstyle\sum}_{j=1}^t\bar \gamma_j \otimes \rho_j={\textstyle\sum}_{j=1}^t u_j=u.$ This shows that ${\rm Ker}(\partial_a^{1-p})\subseteq {\rm Im}(\partial^{-p}_a)$.
The proof of the proposition is completed.

\vspace{2pt}

As an immediate consequence, we obtain the following statement, which genera\-lizes the result stated in \cite[(2.6.1)]{BGS}.

\begin{Theo}\label{Koz-proj-rls}

Let $\La=kQ/R$ be a quadratic algebra with $Q$ a locally finite quiver.
Then $\La$ is Koszul if and only if $\mathcal{K}_a^\cdt\!$ is a graded projective resolution of $S_a$ for $a\in Q_0$.

\end{Theo}

\vspace{2pt}

Next, we shall define the quadratic dual of $\La$ by $Q^{\rm o}$; compare 
\cite{BGS, Mar, MOS}. We need some notation. Let $n\ge 0$. Given $\xi\in Q_n$, let $\xi^*\in D(kQ_n)$ such, for any $\eta\in Q_n$,  that $\xi^*(\eta)=1$ if $\eta=\xi$; and $\xi^*(\eta)=0$ otherwise. Given $\gamma=\sum\lambda_i \xi_i\in kQ_n$ with $\lambda_i\in k$ and $\xi_i\in Q_n$, we write $\gamma^*=\sum \lambda_i \hspace{.5pt} \xi_i^*$. This yields a $k$-isomorphism $\psi_n: kQ_n^{\rm o} \to D(kQ_n): \gamma^{\rm o}\to \gamma^*$.

\vspace{2pt}

If $\xi\in kQ_n(x,y)$, then the restriction of $\xi^*$ to $kQ_n(x,y)$ is also written as $\xi^*$. Since $Q$ is locally finite, $\{\xi^* \mid \xi\in Q_n(x, y)\}$ is the dual basis of $Q_n(x, y)$ in $D(kQ_n(x, y))$. For a subspace $U$ of $kQ_n(x, y)$, we denote by $U^\perp$ the subspace of $D(kQ_n(x, y))$ of $k$-linear functions vanishing on $U$. The following statement is evident.

\begin{Lemma}\label{dual-comp}

Let $Q$ be a locally finite quiver with $x, y, z\in Q_0$ and $s, t\ge 0$.

\begin{enumerate}[$(1)$]

\vspace{-2pt}

\item If $\xi \!\in\! kQ_s(x, y)$ and $\zeta\!\in\! Q_1(y, z)$, then $(\zeta\xi)^*\hspace{-1.5pt}(\eta)=\xi^*\hspace{-1.5pt}(\partial_\zeta(\eta)),$ for all $\eta \hspace{-1.5pt}\in kQ_{s+1}.$

\vspace{.8pt}

\item If $\xi \in kQ_s(x,  y)$ and $\zeta \in kQ_t(y, z),$ then $(\zeta\xi)^*(\gamma\delta)=\zeta^*(\gamma) \hspace{.5pt} \xi^*(\delta)$ for all $\delta \in kQ_s$ and $\gamma  \in kQ_t$.

\item If $U$ and $V$ are $k$-vector subspaces of $kQ_s(x, y)$, then $(U+V)^\perp =U^\perp \cap V^\perp$ and $(U\cap V)^\perp=U^\perp+V^\perp.$

\vspace{1pt}

\end{enumerate} \end{Lemma}

%

Let $R$ be a quadratic ideal in $kQ$. The {\it quadratic dual} of $R$ is the ideal $R^{\hspace{.6pt}!}$ in $kQ^{\hspace{.6pt}\rm o}$ generated by the $R^{\hspace{.6pt}!}_2(y, x)$ with $x, y\in Q_0$, where
$R^{\hspace{.6pt}!}_2(y, x)$ stands for the $k$-subspace of $kQ^{\hspace{.6pt}\rm o}_2(y, x)$ of elements $\rho^{\rm o}$ with $\rho\in kQ_2(x, y)$ such that $\rho^*\in R_2(x, y)^\perp$.

\begin{Defn}\label{qalg}

Let $\La=kQ/R$ be a quadratic algebra, where $Q$ is a locally finite quiver. Then  $\La^{\hspace{.5pt}!}=kQ^{\rm o} / R^{\hspace{.6pt}!}$ is called the {\it quadratic dual} of $\La$.

%
%
%
%
\end{Defn}

The quadratic dual of a quadratic ideal in $kQ$ is described explicitly as follows.

\begin{Lemma}\label{perp}

Let $R$ be a quadratic ideal in $kQ$, where $Q$ is a locally finite quiver. If $\sigma\in kQ_n(x,y)$ with $n\ge 0$, then $\sigma^{\rm o}\in R^!_n\hspace{-1pt}(y,x)$ if and only if $\sigma^*\in R^{(n)}(x,y)^\perp.$

\end{Lemma}

\vspace{-1pt}

\noindent{\it Proof.} Fix $\sigma\in kQ_n(x,y)$ with $x, y\in Q_0$ and $n\ge 0$. We only need to consider the case where $n\ge 3$.
By definition, we have $R^{(n)}(x, y) = \cap_{j=0}^{n-2}R^{(n,j)}(x,y)$ where
\vspace{-3pt} $$\textstyle R^{(n\m,j)\hspace{-1pt}}(x, \hspace{-.8pt} y)\hspace{-2.5pt}=\hspace{-3pt} {\sum}_{a, \hspace{-.5pt} b\in Q_0}\hspace{-.8pt} kQ_{n-\m 2-\m j}(b,\hspace{-.8pt} y) \cdot \hspace{-.3pt} R_2(a, \hspace{-.8pt} b) \cdot \hspace{-.3pt} kQ_j(x, \hspace{-.8pt} a)\m;
\vspace{-3pt} $$ \vspace{.5pt} and $R^{!}_n\m (y, x) \hspace{-2pt}=\hspace{-3pt} {\sum}_{j=0}^{n-2} R^{!}_{n,j}(y, x)$ where \vspace{-3pt}  $$R^{\hspace{.5pt}!}_{n,j}(y, x)={\textstyle\sum}_{a,b\in Q_0} kQ^{\rm o}_j(a, x) \m\cdot\m R^{\hspace{.5pt}!}_2(b,a) \m\cdot\m kQ^{\rm o}_{n-2-j}(y,b).$$

\vspace{-3pt}

For the necessity, we may assume that $\sigma^{\rm o}=\gamma^{\rm o} \eta^{\rm o} \delta^{\rm o}$, where $\delta \in kQ_{n-2-j}(b, y),$ $\gamma \in kQ_j (x, a)$, and $\eta \in kQ_2(a, b)$ such that $\eta^{\rm o}\in R^{\hspace{.5pt}!}_2(b, a)$, for some $a, b\in Q_0$ and $0\le j\le n-2$. Given $w \in R^{(n)}(x, y)$, we may write $w = \sum_{i=1}^s \delta_i \hspace{.5pt} \eta_i \hspace{.5pt} \gamma_i,$ for some $\gamma_i\in kQ_j(x,a_i)$, $\eta_i\in R_2(a_i, b_i)$ and $\delta_i\in k_{n-j-2}Q(b_i, y)$, where $a_i, b_i\in Q_0$. Since $\eta^*\in R_2(a,b)^\perp$, we see that $\eta^*(\eta_i)=0$, for $i=1, \ldots, s$.
By Lemma \ref{dual-comp}(1), $\sigma^*(w)=(\delta\eta\gamma)^*(w) = \sum_{i=1}^s\delta^*\hspace{-1pt}(\delta_i) \, \eta^*\hspace{-1pt}(\eta_i)  \,\gamma^*\hspace{-1pt}(\gamma_i) =0.$ That is, $\sigma^*\in R^{(n)}(x,y)^\perp.$

\vspace{1pt}

Conversely, since $R^{(n)}\m(x,y)^\perp\!\m=\!\m{\sum}_{j=0}^{n-2}\, R^{(n,j)}(x,y)^\perp$; see (\ref{dual-comp}), we may assume that
$0\ne \sigma^* \!\in\! R^{(n,p)}\m(x,y)^\perp$ for some $0\le p\le n-2$, and show that $\sigma^{\rm o}\!\in\! R^{\hspace{.5pt}!}_{n,p}(y,x)$.
Write $\sigma=\sum_{i=1}^t \sigma_i$, where $\sigma_i \in kQ_{n-2-p}(b_i, y)\cdot kQ_2(a_i, b_i)\cdot kQ_p(x, a_i)$. By Lemma \ref{dual-comp}(1), $\sigma_i^*\in (kQ_{n-p-2}(b_j, y)\cdot kQ_2(a_j, b_j)\cdot kQ_p(x, a_j))^\perp$ for $j\ne i$, and hence, $\sigma_i^*\in R^{(n,p)}(x,y)^\perp,$ for $i=1, \ldots, t$. So, we may assume that $\sigma=\delta_1\zeta_1 \gamma_1$, for some $\delta_1\in kQ_{n-p-2}(b_1, y)$, $\zeta_1\in kQ_2(a_1, b_1)$ and $\gamma_1\in kQ_p(x, a_1)$. Since $\sigma^*\ne 0$, by Lemma \ref{dual-comp}(1), $\delta_1^*$ and $\gamma_1^*$ are non-zero. Hence,
$\delta_1^*(\nu_1)=\gamma_1^*(\mu_1)=1$, for some $\nu_1\in kQ_{n-2-p}(b_1, y)$ and $\mu_1\in kQ_p(x, a_1)$.

Choose a basis $\{\rho_{1}, \ldots, \rho_{r}, \rho_{r+1}, \ldots, \rho_{s}\}$ of $kQ_2(a_1, b_1)$, where $\{\rho_{1}, \ldots, \rho_{r}\}$ is a basis of $R_2(a_1, b_1)$. There exists a basis $\{\eta_{1}, \ldots, \eta_{r}, \eta_{r+1}, \ldots, \eta_{s}\}$ of $kQ_2(a_1, b_1)$ such that $\{\eta_{1}^*, \ldots, \eta_{r}^*, \eta_{r+1}^*, \ldots, \eta_s^*\}$ is the dual basis \vspace{.5pt} of $\{\rho_{1}, \ldots, \rho_{r}, \rho_{r+1}, \ldots, \rho_{s}\}$, and hence, $\{\eta_{r+1}^{\rm o}, \ldots, \eta_{s}^{\rm o}\}$ is a basis of $R^{\hspace{.5pt}!}_2(b_1, a_1)$. Writing $\zeta_1=\sum_{j=1}^s \lambda_j\eta_j$ where $\lambda_j \in k$, we get $\sigma^*=\sum_{j=1}^s\lambda_j(\delta_1\eta_j\gamma_1)^*\in
(Q_{n-2-p}(b_1,y) \cdot  R_2(a_1,b_1) \cdot kQ_p(x,a_1))^\perp$. In view of Lemma \ref{dual-comp}(1), we see that \vspace{-3pt}
$$\textstyle 0 = \sigma^*(\nu_1\rho_i\mu_1) = \sum_{j=1}^s \lambda_j (\delta_1\eta_j \gamma_1)^*(\nu_1\rho_{i}\mu_1) = \sum_{j=1}^s \lambda_j  \delta_1^*(\nu_1) \eta_j^*(\rho_i) \gamma_0^*(\mu_1)= \lambda_i, \vspace{-3pt}$$ for $i=1, \ldots, r$.
Therefore, $\sigma^*=\sum_{j=r+1}^{s} \lambda_j (\delta_1\eta_j\gamma_1)^*$, and consequently,
we have
$\sigma^{\rm o} =\sum_{j=r+1}^{s} \lambda_j\gamma_1^{\rm o} \eta_j^{\rm o} \delta_1^{\rm o}\in R^{\,!}_{n,p}(y,x).$ The proof of the lemma is completed.

\vspace{5pt}

The following statement in particular justifies the terminology of quadratic dual; compare \cite[(2.8.1)]{BGS}.

\begin{Prop}\label{q-dual}

Let $\La=kQ/R$ be a quadratic algebra with $Q$ a locally finite quiver. Then $\La^!$ and $\La^{\rm o}$ are quadratic with $(\La^!)^!= \La$ and $(\La^{\rm o})^!=(\La^!)^{\rm o}$.

\end{Prop}

\vspace{-2pt}

\noindent{\it Proof.} Clearly, $\La^{\rm o}$ and $\La^!$ are quadratic with $(\La^!)^!=k((Q^{\rm o})^{\rm o})/(R^{\hspace{.5pt}!})^!=kQ/(R^{\hspace{.5pt}!})^!$ and $(\La^{\rm o})^!=k((Q^{\rm o})^{\rm o})/(R^{\hspace{.5pt}!})^{\rm o}=kQ/(R^{\hspace{.5pt}!})^{\rm o}$. 
Let $\gamma, \rho\in kQ_2(x, y)$ with $x, y\in Q_0$. It is easy see that $(\gamma^{\rm o})^*(\rho^{\rm o})=\rho^*(\gamma)=\gamma^*(\rho)$.

Now, $\gamma\in (R^{\hspace{.5pt}!})^!_2(x, y)$ if and only if $(\gamma^{\rm o})^*(\rho^{\rm o})=0$ for all $\rho^{\rm o}\in R^{\hspace{.5pt}!}_2(y, x)$, or equivalently, $\rho^*(\gamma)=0,$ for all $\rho^*\in R_2(x, y)^\perp$, that is, $\gamma\in R_2(x, y)$. This implies that $(R^{\hspace{.5pt}!})^!=R$, and hence, $(\La^!)^!=kQ/R=\La$.

Next, $\gamma\in (R^{\rm o})^!_2(x, y)$ if and only if, $(\gamma^{\rm o})^*(\rho^{\rm o})=0$ for all $\rho^{\rm o} \in R_2^{\rm o}(y, x)$, if and only if $\gamma^*(\rho)=0$, for all $\rho\in R_2(x, y)$, or equivalently, $\gamma^{\rm o}\in R^{\hspace{.5pt}!}_2(y, x)$, that is, $\gamma\in (R^{\hspace{.5pt}!})^{\rm o}_2(x, y)$. This implies that $(R^{\rm o})^!=(R^{\hspace{.5pt}!})^{\rm o}.$ Now, $(\La^{\rm o})^!=kQ/(R^{\hspace{.5pt}!})^{\rm o}=(\La^!)^{\rm o}.$
The proof of the proposition is completed.

\vspace{3pt}

We shall give an alternative description of the local Koszul complexes in terms of $\La^{\hspace{.6pt}!}$. We need to fix some notation.
Write $\bar{\gamma\hspace{.8pt}}^!=\gamma^{\rm o}+ R^{\hspace{.5pt}!}$ for $\gamma\in kQ^+$, but $e_x=\varepsilon_x+R^{\hspace{.6pt}!}$ for $x\in Q_0$. Then $\La^!=\oplus_{n\ge 0}\La^!_n$, where $\La^!_n=\{\bar{\gamma\hspace{.6pt}}^! \mid \gamma \in kQ_n\}$ for $n\ge 0$. \vspace{1pt} Given any $x\in Q_0$, we shall write $P^{\hspace{.6pt}!}_x=\La^{\hspace{.2pt}!}e_x$, and $S_x^!=P^{\hspace{.6pt}!}_x/\rad P^{\hspace{.6pt}!}_x$, and $I_x^!=\mk D ((\La^!)^{\rm o}e_x)$.

\vspace{2pt}

Fix $a \in Q_0$. Given $n\in \Z$, we set $\mathcal{P}_a^{-n} \! = \! \oplus_{x\in Q_0} P_x\mla -n\mra \otimes D(e_a\La^!_ne_x)\in {\rm gproj}\La$. Given $\alpha \!\m\in\!\m Q_1(y, x)$, we have a morphism
$P[\bar \alpha]:  P_x \sla -n\nra \to  P_y \sla 1 -n\nra$ in ${\rm gproj}\La$, that is
the right multiplication by $\bar \alpha$; and a $k$-linear map $P[\bar\alpha^!]: e_a\La^{\hspace{.5pt}!}_{n-1} e_y \to e_a\La^!_n e_x$, that is
the right multiplication by $\bar \alpha^!$.
Thus, we obtain a graded $\La$-linear morphism \vspace{-8pt}
$$P[\bar\alpha]\otimes DP[\bar\alpha^!]: P_x \sla -n\nra \otimes D(e_a\La^!_{n}e_x)\to P_y \sla 1-n\nra \otimes D(e_a\La^!_{n-1}e_y).\vspace{-0pt}$$

Write $\mathcal{P}_a^{1-n}=\oplus_{y\in Q_0} P_y\tla \m 1-\m n\nra \otimes D(e_a\La^!_{n-1}e_y)$ \vspace{1pt} and consider the graded morphism $\ell^{-n}_a=(\sum_{\alpha\in Q_1(y,x)} \hspace{-2pt} P[\bar\alpha]\otimes DP[\bar\alpha^!])_{(x,y)\in Q_0\times Q_0}: \mathcal{P}_a^{-n}\to \mathcal{P}_a^{1-n}$. We obtain a sequence
\vspace{-9pt} $$\xymatrixcolsep{22pt}\xymatrix{\mathcal{P}_a^\pdt: \quad \cdots \ar[r] & \mathcal{P}_a^{-n} \ar[r]^-{\ell_a^{-n}} & \mathcal{P}_a^{1-n} \ar[r] &
\cdots \ar[r] & \mathcal{P}_a^{-1} \ar[r]^{\ell_a^{-1}} & \mathcal{P}_a^{\hspace{.4pt}0} \ar[r] & 0 \ar[r] & \cdots} \vspace{-1pt}$$ in ${\rm gproj}\La$, which is a complex as shown below.

\begin{Lemma}\label{k-cplx-iso}

Let $\La=kQ/R$ be a quadratic algebra with $Q$ a locally finite quiver.
Then $\mathcal{K}_a^\pdt\m\cong \m \mathcal{P}_a^\pdt,$ for all $a\in Q_0$.

\end{Lemma}

\noindent{\it Proof.} Fix $a, x\in Q_0$ and $n\ge 0$. Then, $D(kQ_n(a, x))=\{\gamma^* \mid \gamma\in kQ_n(a, x)\}$ and $e_a\La^!_ne_x=\{\bar{\gamma\hspace{.6pt}}^!=\gamma^{\rm o}+R^{\hspace{.5pt}!} \mid \gamma\in kQ_n(a, x)\}$. By Lemma \ref{perp}, $\gamma^*\in  R^{(n)}(a,x)^\perp$ if and only if $\gamma^{\rm o}\in R^{\hspace{.5pt}!}_n(x, a)$. Therefore, we obtain a $k$-bilinear form  \vspace{-1pt}
$$\langle -, - \rangle: R^{(n)}(a,x)\times e_a\La^!_ne_x \to k: (\delta, \bar{\gamma\hspace{.8pt}}^!) \mapsto \gamma^*(\delta),
 \vspace{-1pt}$$
which is non-degenerate on the right. If $\delta\in R^{(n)}(a,x)$ is non-zero, then $\gamma^*(\delta)\ne 0$, that is, $\langle\delta, \gamma^!\rangle\ne 0$, for some $\gamma\in kQ_n(a, x)$. Hence, $\langle -, - \rangle$ is non-degenerate. This yields a $k$-isomorphism
$\phi^{n}_x: R^{(n)}(a,x)\to D(e_a\La^!_ne_x): \delta \to \langle\delta, -\rangle.$ Given $\alpha\in Q_1(y, x)$ and $n> 0$,
we claim that the following diagram commutes: \vspace{-3pt}
$$\xymatrixrowsep{16pt}\xymatrixcolsep{30pt}\xymatrix{
\hspace{5pt} R^{(n)}(a,x) \ar[r]^-{\partial_\alpha} \ar[d]_-{\phi^{n}_x} & R^{(n-1)}(a, y)\ar[d]^{\phi^{n-1}_y} \hspace{8pt} \\
D(e_a\La^!_{n}e_x) \ar[r]^-{DP[\bar\alpha^!]} &  D(e_a\La^{!}_{n-1}e_y).
}\vspace{-5pt} $$
Indeed, for any $\rho\in R^{(n)}(a, x)$ and $\zeta\in kQ_{n-1}(a, y)$, by Lemma \ref{dual-comp}(2), we have \vspace{-2pt}
$$DP[\bar\alpha^!](\phi^{n}_x(\rho))(\bar\zeta\hspace{.6pt}^!) = \phi^{n}_x(\rho)(\bar\zeta\hspace{.6pt}^!\bar\alpha^!) = (\alpha\zeta)^*(\rho)= \zeta^*(\partial_\alpha(\rho)) =\phi^{n-1}_y(\partial_\alpha(\rho))(\bar\zeta\hspace{.6pt}^!). \vspace{-2pt}$$

Thus, the above commutative diagram commutes. It is now easy to see that the graded $\La$-linear isomorphisms
$\oplus_{x\in Q_0} ({\rm id}\otimes \phi^n_x)$ with $n\in \Z$
gives rise to a complex isomorphism $\mathcal{K}_a^\pdt\m \cong \m \mathcal{P}_a^\pdt$. The proof of the lemma is completed.

\vspace{3pt}

Next, we shall consider colinear injective coresolutions of graded simple mo\-dules. Given $u \in e_a\La_{t-s} e_b$, the right multiplication by $u^{\rm o}$ yields a morphism $P[u^{\rm o}]:\m P_b^{\hspace{.5pt} \rm o} \tla -t \tra \!\to\! P_a^{\hspace{.5pt} \rm o}\tla -s \tra$ in ${\rm gproj}\La^{\rm o}$. Then, $I[u]=\mf{D}(P[u^{\rm o}]): I_a \sla s\sra  \to I_b\tla t \tra$ is a morphism in $ {\rm ginj}\La$. Note that this notation does not distinguish $I[u]$ from its grading shifts. For simplicity, we put $\hspace{1.5pt}\hat{\hspace{-3pt}\La}=(\La^!)^{\rm o}=kQ/(R^{\hspace{.6pt}!})^{\rm o}$. Write $\hspace{2pt}\hat{\hspace{-2.5pt}\gamma}=\gamma+ (R^{\hspace{.6pt}!})^{\rm o}$ for $\gamma\in kQ^+$ and $e_x=\varepsilon_x+(R^{\hspace{.5pt}!})^{\rm o}$ for $x\in Q_0$. In this way, $\hspace{1.5pt}\hat{\hspace{-3pt}\La}=\oplus_{n\ge 0}\hspace{3pt}\hat{\hspace{-3pt}\La}_n$, where $\hspace{1.5pt}\hat{\hspace{-3pt}\La}_n=\{\hspace{2pt}\hat{\hspace{-2.5pt}\gamma} \m \mid \!\gamma\in kQ_n\}$ for all $n\ge 0$.

Fix $a\in Q_0$. Given $n\in \Z$, we set $\mathcal{I}_a^n=\oplus_{x\in Q_0} I_x\mla n\mra \otimes e_x\La^!_ne_a\in {\rm ginj}\La$. For $n\ge 0$ and $\alpha\in Q_1(x, y)$, we have a morphism $I[\bar \alpha]: I_y \hspace{.5pt} \langle n-1\nra\to I_x \langle n \rangle$ in ${\rm ginj}\La$. And the left multiplication by $\bar \alpha^!$ yields a $k$-linear map $P_a^{\hspace{.4pt}!}(\bar\alpha^!): e_y\La^!_{n-1}e_a\to e_x\La^!_{n}e_a$. This yields a graded $\La$-linear morphism
$$I[\bar{\alpha}] \otimes P_a^{\hspace{.4pt}!}(\bar \alpha^!): I_y\mla n-1\mra \otimes e_y\La^!_{n-1}e_a\to I_x\mla n \mra \otimes e_x\La^!_ne_a. \vspace{-1pt} $$

Write \vspace{1.5pt} $\mathcal{I}_a^{n-1}=\oplus_{y\in Q_0} I_y \sla n-1\mra \otimes e_y\La^!_{n-1}e_a$
and consider the graded morphism $d^n_a\!=\!\m(\textstyle \sum_{\alpha\in Q_1(x,y)}\m I[\bar{\alpha}]\otimes P_a^{\hspace{.4pt}!}\m(\alpha^!)\m)_{(x,y)\in Q_0\times Q_0}\m: \mathcal{I}_a^{n-1} \to \mathcal{I}_a^n.$ We obtain a sequence \vspace{-6pt}
$$\xymatrixcolsep{20pt}\xymatrix{\mathcal{I}^\pdt_a: & \cdots\ar[r] & 0 \ar[r] & \mathcal{I}_a^{\hspace{.3pt}0} \ar[r]^{d^1_a} & \mathcal{I}_a^1 \ar[r] & \cdots \ar[r] & \mathcal{I}_a^{n-1} \ar[r]^-{d^n_a} & \mathcal{I}_a^n\ar[r] & \cdots}$$ in ${\rm ginj}\La$, which is a complex as shown below.

\begin{Prop}\label{Koz-cplx-dual}

Let $\La=kQ/R$ be a quadratic algebra with $Q$ a locally finite quiver. If $a\in Q_0$, then $\mathcal{I}_a^\ydt$ is a complex, and it is a graded injective coresolution of $S_a$ if and only if $S_a$ admits a colinear injective coresolution over ${\rm ginj}\La$.

\end{Prop}

\noindent{\it Proof.} By Proposition \ref{q-dual}, $(\La^{\rm o})^{\hspace{.6pt}!}= (\La^{\hspace{.6pt}!})^{\rm o} = \hspace{2.3pt}\hat{\hspace{-3pt}\La}=\{\hat \gamma \mid \gamma \in kQ\}$, where $\hspace{2pt}\hat{\hspace{-2.5pt}\gamma}=\gamma+ (R^{\hspace{.6pt}!})^{\rm o}$ for $\gamma\in kQ$. Fix $a\in Q_0.$ As stated in Lemma \ref{k-cplx-iso}, we have a complex
\vspace{-4pt}
$$\xymatrixcolsep{23pt}\xymatrix{\mathcal{P}_{\! a^{\rm o}}^\pdt: \hspace{5pt} \cdots \ar[r] & \mathcal{P}_{\! a^{\rm o}}^{-n} \ar[r]^-{\ell^{-n}} & \mathcal{P}_{\! a^{\rm o}}^{1-n}\ar[r] & \cdots \ar[r] & \mathcal{P}_{\! a^{\rm o}}^{-1}\ar[r]^{\ell^{-1}} & \mathcal{P}_{\! a^{\rm o}}^0 \ar[r] & 0 \ar[r] & \cdots, } \vspace{-3pt}$$
where \vspace{1pt} $\mathcal{P}_{\! a^{\rm o}}^{-n} \! = \! \oplus_{x\in Q_0} P_x^{\hspace{.4pt}\rm o}\mla - \m n \nra \m \otimes \! D(e_a \hat{\hspace{-.5pt}\La}_{n}e_x\m)$; \hspace{.5pt} $\mathcal{P}_{\! a^{\rm o}}^{1-n} \! = \! \oplus_{y\in Q_0} P_y^{\hspace{.4pt}\rm o}\mla 1 - n \mra \m \otimes \! D(e_a\; \hat{\hspace{-3pt}\La}_{n-1}e_y\m)$ and $\ell^{-n} = (\sum_{\alpha\in Q_1(x,y)} \! P[\bar\alpha^{\rm o}] \otimes DP[\hat\alpha])_{(y,x)\in Q_0\times Q_0}$ for all $n\ge 1$.

First, we show that $\mf{D}(\mathcal{P}_{\! a^{\rm o}}^\pdt) \!\cong \! \mathcal{I}_a^\pdt$.
Given any $n\ge 0$, since $Q$ is locally finite, $e_a\,\hat{\hspace{-3pt}\La}_{n}$ is finite dimensional. By Proposition \ref{gr-duality}, $\mf{D}(\mathcal{P}_{\! a^{\rm o}}^{-n})= \oplus_{x\in Q_0} I_x\nla n\nra \otimes D^2(e_a\hspace{2.3pt}\hat{\hspace{-3pt}\La}_{n}e_x)$ and
$\mf{D}(\ell^{-n})=(I[\bar \alpha]\otimes D^2\m P[\hat\alpha])_{(y,x)\in Q_0\times Q_0}$. 
Moreover, since $\hspace{1.5pt}\hat{\hspace{-3pt}\La}=(\La^!)^{\rm o},$ we have a $k$-linear isomorphism $\sigma^n_x: e_a\hspace{2pt}\hat{\hspace{-3pt}\La}_ne_x\to e_x\La^!_ne_a$, sending $\hat{\gamma\hspace{1.7pt}}\mapsto \bar{\gamma\hspace{.7pt}}^!$. Composing this with the canonical $k$-isomorphism
$\varphi^n_x: D^2(e_a\hspace{2pt}\hat{\hspace{-3pt}\La}_ne_x) \to e_a \hspace{2pt}\hat{\hspace{-3pt}\La}_ne_x$, we obtain a $k$-linear isomorphism $\theta^n_{x}=\sigma^n_x \!\circ\! \varphi^n_x: D^2(e_a\hspace{2pt}\hat{\hspace{-3pt}\La}_ne_x) \to e_x\La^!_ne_a$. It is easy to verify that
 \vspace{-3pt} $$\xymatrixrowsep{18pt}\xymatrixcolsep{32pt}\xymatrix{
D^2(e_a \hspace{2pt} \hat{\hspace{-3pt}\La}_{n-1} e_y) \ar[r]^-{\varphi^{n-1}_y} \ar[d]_{D^2\! P[\hat\alpha]} & e_a \hspace{2pt}\hat{\hspace{-3pt}\La}_{n-1} e_y \ar[r]^{\sigma_y^{n-1}} \ar[d]^{P[\hat\alpha]} & e_y \La^!_{n-1} e_a\ar[d]^-{P_a^{\hspace{.4pt}!}(\bar\alpha^!)} \\
D^2(e_a \hspace{2.5pt}\hat{\hspace{-3pt}\La}_{n} e_x)\ar[r]^{\varphi^n_x}
& e_a \hspace{2.5pt}\hat{\hspace{-3pt}\La}_{n} e_x \ar[r]^{\sigma_x^n} & e_x \La^!_{n} e_a
}\vspace{-1pt}$$ commutes for every $\alpha\in Q_1(x,y)$.
This yields a complex isomorphism $\mf{D}(\mathcal{P}_{\! a^{\rm o}}^\pdt) \!\cong \! \mathcal{I}_a^\pdt$, \vspace{.5pt} given by $\oplus_{x\in Q_0} ({\rm id}\otimes \theta^n_x): \oplus_{x\in Q_0} I_x \tla n \nra \!\otimes\! D^2(e_a\hspace{2.5pt}\hat{\hspace{-3pt}\La}_{n} e_x) \! \to\m  \oplus_{x\in Q_0} I_x \tla n\nra \m \otimes\m e_x \La^!_{n} e_a$ with $n\in \Z$.
So, $\mathcal{I}_a^\ydt$ is a complex. Since $\mathcal{I}_a^n$ is co-generated in degree $-n$, we see that $\mathcal{I}_a^\ydt$ is a colinear injective coresolution if it is a graded injective co-resolution of $S_a$.

Suppose that $S_a$ has a colinear injective coresolution $\mathcal{I}^\pdt$ over ${\rm ginj}\hspace{.5pt}\La$. In view of Lemma \ref{n-presentation-dual}, we see that $\mf{D}(\mathcal{I}^\pdt)$ is a linear projective resolution of $S_a^{\rm o}$ over ${\rm gproj }\La^{\rm o}.$ By Theorem \ref{Koz-proj-rls} and Lemma \ref{k-cplx-iso}, $\mf{D}(\mathcal{I}^\pdt)\cong \mathcal{P}_{\! a^{\rm o}}^\pdt$. Hence, $\mathcal{I}^\pdt\cong \mf{D}^2(\mathcal{I}^\pdt)\cong \mf{D}(\mathcal{P}_{\! a^{\rm o}}^\pdt) \cong \mathcal{I}_a^\pdt$. The proof of the proposition is completed.

\vspace{2pt}

Given $M, N\in \GrLa$, we write ${\rm GExt}_{\mathit\Lambda}^n(M, N)$ for the $n$-th graded extension group in $\GrLa$;
see \cite[Section III.5]{Mac}.
The following statement includes the classical results stated in \cite[(2.2.1), (2.10.2)]{BGS}. \vspace{-1pt}

\begin{Theo}\label{Opp-Koszul}

Let $\La=kQ/R$ be a quadratic algebra with $Q$ a locally finite quiver. The following statements are equivalent.

\begin{enumerate}[$(1)$]

\vspace{-1pt}

\item The algebra $\La$ is Koszul.

\vspace{0pt}

\item The opposite algebra $\La^{\rm o}$ is Koszul.

\vspace{0pt}

\item  The quadratic dual $\La^!$ is Koszul.

\vspace{0pt}

\item Every graded simple $\La$-module has a colinear injective coresolution over ${\rm ginj}\La$.

\end{enumerate}

\end{Theo}

\noindent{\it Proof.} By Lemma \ref{n-presentation-dual}, Statements (2) and (4) are equivalent. By Proposition \ref{q-dual}, if Statement (1) implies (4),
then Statements (1) and (2) are equivalent; and if Statement (1) implies (3), then Statements (1) and (3) are equivalent. Thus, it suffices to show that Statement (1) implies Statements (3) and (4).

Assume that $\La$ is Koszul. Fix $a\in Q_0$. Since $\La^{\rm o}$ is quadratic; see (\ref{q-dual}), by Theorem \ref{qua-alg}, $S_a^{\rm o}$ has a linear projective $2$-presentation over ${\rm gproj}\La^{\rm o}$. And by Lemma \ref{n-presentation-dual}, we may assume that $S_a$ has a colinear injective $(n-1)$-copresentation \vspace{-5pt}
$$\xymatrixcolsep{20pt}\xymatrix{0\ar[r] &S_a \ar[r]^{d^{\hspace{.3pt}0}} & I^0 \ar[r]^{d^1} & I^1 \ar[r] & \cdots \ar[r] & I^{n-2} \ar[r]^{d^{n-1}} & I^{n-1}}
\vspace{-3pt}$$ over ${\rm ginj}\La$, for some $n\ge 3$. Denote by $c^n: I^{n-1}\to C^{n}$ the cokernel of $d^{n-1}$. Given any $b\in Q_0$ and $p\in \Z$,
since $\soc I^{n-1}\subseteq {\rm Ker}(c^n)$, it is well-known that ${\rm GExt}_{{\mathit\Lambda}}^{n}(S_b \tla p\tra, S_a)\cong {\rm GHom}_{{\mathit\Lambda}}(S_b \tla p\tra, C^{n})$; see \cite[(III.6.4)]{Mac}. Since $\La$ is Koszul, $S_b$ has a linear projective resolution \vspace{-2pt}
$$\xymatrixcolsep{20pt}\xymatrix{\cdots \ar[r] & P^{-i} \ar[r] & P^{1-i} \ar[r] & \cdots \ar[r] & P^{-1} \ar[r] & P_b \ar[r] & S_b\ar[r] &0.}
\vspace{-2pt}$$  over ${\rm proj}\La$. So, ${\rm GHom}_{{\mathit\Lambda}}(S_b\tla p\tra, C^{n}\m\m) \m\cong\m {\rm GExt}_{{\mathit\Lambda}}^{n}(S_b, S_a\sla -p \hspace{.4pt}\sra) \!\cong\m {\rm GHom}_{{\mathit\Lambda} \m}(P^{-n}\!\m, S_a\sla \m-p \hspace{.4pt}\sra\m).$
Since $P^{-n}$ is finitely generated generated in degree $n$ and $S_a\sla -p \sra$ is generated in degree $p$, we deduce that ${\rm GHom}_{{\mathit\Lambda}}(S_b\tla p\tra, C^{n}\m\m)=0$ for all but finitely many $b\in Q_0$ and for all $p\ne n$. Hence, ${\rm soc} \,C^n$ is finitely generated in degree $-n$. Since $I^{n-1}$ is bounded above, so is $C^n$. Hence, $\soc C^n$ is graded essential in $C^n$; see \cite[(2.9.2)]{LLi}. That is, $C^n$ is finitely cogenerated in degree $-n$. By Proposition \ref{p-cover}(2), there exists a graded injective envelope $q^n: C^n\to I^n$, where $I^n\in {\rm ginj}\La$ is cogenerated in degree $-n$. Thus, $S_a$ has a colinear injective $n$-copresentation over ${\rm ginj}\La$. By induction, Statement (4) holds.
It remains to show that Statement (3) holds.

Since $\La$ is quadratic, $\La^!$ is quadratic with $(\La^!)^!=\La$; see (\ref{q-dual}). By Lemma \ref{k-cplx-iso}, the local Koszul complex at $a$ of $\La^!$ is isomorphic to the complex \vspace{-4pt}
$$\xymatrixcolsep{20pt}\xymatrix{\mathcal{P}_{a^!}^\pdt: & \cdots \ar[r] & \mathcal{P}_{a^!}^{-n}\ar[r]^{\ell^{-n}} & \mathcal{P}_{a^!}^{1-n}\ar[r] & \cdots \ar[r] & \mathcal{P}_{a^!}^{-1}\ar[r]^{\ell^{-1}} & \mathcal{P}_{a^!}^0
\ar[r] & 0 \ar[r] & \cdots,} \vspace{-2pt}$$
where \vspace{.5pt} $\mathcal{P}_{a^!}^{-n}\!=\!\oplus_{x\in Q_0} P^{\hspace{.5pt}!}_x \tla \m - \m n\nra\otimes D(e_a\La_ne_x)$ and $\mathcal{P}_{a^!}^{1-n}\!=\!\oplus_{y\in Q_0} P^{\hspace{.5pt}!}_y\tla 1-n\nra\otimes D(e_a\La_{n-1}e_y),$ and \vspace{.5pt}
$\ell^{-n} 
= ({\textstyle\sum}_{\alpha\in Q_1(x,y)}P[\bar\alpha^!]\otimes DP[\bar\alpha])_{(y,x)\in Q_0\times Q_0}$.
Fix $n>0$. We claim that  \vspace{-5pt}
$$\xymatrixrowsep{2pt}\xymatrixcolsep{22pt}\xymatrix{
\hspace{-10pt} (\dag) \hspace{20pt} \oplus_{z\in Q_0}e_b\La^!_{s-1} e_z\otimes D(e_a\La_{n+1}e_z)
\ar[r]^-{\ell^{-n-1}_{n+s,b}} & \oplus_{x\in Q_0}e_b\La^!_s e_x\otimes D(e_a\La_n e_x) \\
\hspace{165pt}\ar[r]^-{\ell^{-n}_{n+s,b}} &
\oplus_{y\in Q_0} e_b \La^!_{s+1} e_y\otimes D(e_a\La_{n-1}e_y)\hspace{-15pt}} \vspace{-4pt}$$
is exact for any $(s, \m b)\!\in\m \Z\times Q_0$, where $\ell^{-n}_{n+s,b}\! = \! ({\textstyle\sum}_{\alpha\in Q_1(x,y)}P[\bar\alpha^!]\otimes DP[\bar\alpha])_{(y,x)\in Q_0\times Q_0}\m.$

If $s<0$, then $e_b\La^!_s e_x=0$, and $(\dag)$ is evidently exact. In case $s=0$, it becomes \vspace{-6pt}
$$\xymatrix{0 \ar[r] & e_b\La^!_0 e_b \otimes D(e_a \La_n e_b) \ar[r]^-{\ell_{n,b}^{-n}} & \oplus_{y\in Q_0} e_b \La^!_1 e_y\otimes D(e_a\La_{n-1}e_y),}\vspace{-3pt}$$
where $\ell_{n,b}^{-n}=(\ell_{n,b}^{-n}(y,b))_{y\in Q_0}$ with $\ell_{n,b}^{-n}(y,b)=\sum_{\alpha\in Q_1(b,y)} P[\bar \alpha^!] \otimes
DP[\bar \alpha]$. Consider $0\ne f\in D(e_a \La_n e_b)$. Since $n>0$, there exist
$\beta\in Q_1(b,y)$ and $u\in e_a\La_{n-1} e_y$ with $y\in Q_0$ such that $f(u \bar\beta) \ne 0$, that is, $(DP[\bar\beta])(f)(u)\ne 0$. So, $(DP[\bar\beta])(f)\ne 0$. Now, $\ell_{n,b}^{-n}(y,b)(e_b\otimes f)={\textstyle\sum}_{\alpha\in Q_1(b,y)}\bar\alpha^!\otimes (DP[\bar\alpha])(f),$
which is non-zero. 
Thus, $\ell_{n,b}^{-n}$ is a monomorphism. That is, the sequence $(\dag)$ is exact in this case.

Let $s>0$. By Theorem \ref{Koz-proj-rls} and Lemma \ref{k-cplx-iso}, $S_b$ has a linear projective resolution $\mathcal{P}_b^\pdt$, which is exact in degree $-s$. Writing $\mathcal{P}_b^{-s}=\oplus_{x\in Q_0} P_x \sla -s \sra \otimes D(e_b\La^!_se_x)$, and restricting it to the $(n+s, a)$-piece, we obtain an exact sequence \vspace{-6pt}
$$\xymatrixrowsep{2pt}\xymatrixcolsep{28pt}\xymatrix{
(\ddag) \hspace{20pt} \oplus_{y\in Q_0}e_a \La_{n-1} e_y \otimes D(e_b\La^!_{s+1} e_y)
\ar[r]^-{\ell_{b, n+s, a}^{-s-1}} & \oplus_{x\in Q_0} e_a \La_n e_x \otimes D(e_b\La^!_s e_x)\hspace{20pt} \\
\hspace{175pt} \ar[r]^-{\ell_{b,n+s,a}^{-s}} & \oplus_{z\in Q_0} e_a\La_{n+1} e_z \otimes D(e_b\La^!_{s-1}e_z),} \vspace{-6pt}$$
where $\ell_{b,n+s, a}^{-s-1} = (\sum_{\alpha\in Q_1(x,y)}\hspace{-2pt}P[\bar\alpha] \otimes DP[\bar\alpha^!])_{(x,y)\in Q_0\times Q_0}.$
Applying $D$ to the sequence $(\ddag)$, by Lemma \ref{dual-mix}, we obtain an exact sequence isomorphic to $(\dag)$. This proves our claim. Thus, $\mathcal{P}_{a^!}^\pdt$ is a graded projective resolution of $S^!_a$. By Theorem \ref{Koz-proj-rls}, $\La^!$ is Koszul. The proof of the theorem is completed.

\vspace{2pt}

\noindent{\sc Remark.} In case $\La$ is Koszul, one calls $\La^!$ the {\it Koszul dual} of $\La$.


\vspace{2pt}

\noindent{\sc Example.} Consider $\La=kQ/(kQ^+)^2$, where $Q$ is a locally finite quiver. Then $\La^!=kQ^{\rm op}$, which is Koszul. By Theorem \ref{Opp-Koszul}, $\La$ is Koszul.

\vspace{2pt}

We shall conclude this section with two applications. Some quadratic algebras are known to be Koszul, such as symmetric algebras and exterior algebras; see \cite[page 476]{BGS}, and monomial quadratic algebras; see \cite[(2.19)]{Mar}. Using our description of the local Koszul complexes, we will provide a new class of Koszul algebras. For this purpose, we introduce a condition $(*)$ for a quadratic special multi-serial algebra $\La=kQ/R$ as follows:

\vspace{.5pt}

Let $\sum_{i=1}^s\lambda_i \beta_i\alpha_i$ be a polynomial relation in $R_2(x,z)$ with $\alpha_i, \beta_i \in Q_1$ and $\lambda_i\in k$ such that $\zeta\beta_i\notin R_2(-, b)$ for some $\zeta \in Q_1(z,b)$ and $1\le i\le s$. If $\gamma \in Q_1(a, x)$, then $\alpha_j\gamma$ is a summand of a relation in $R_2(a, -)$, for all $1\le j\le s$ with $j\ne i$.

\vspace{2pt}

\noindent{\sc Remark.} A quadratic special multi-serial algebra $\La$ satisfies the condition $(*)$ if and only if $\La^{\rm o}$ satisfies its dual, which is left for the reader to formulate explicitly.

\begin{Lemma}\label{msa-k-relation}

Let $\La=k Q/R$ be a quadratic special multi-serial algebra with $Q$ a locally finite quiver, satisfying the condition $(*)$. Consider a polynomial relation
$\sum_{i=1}^s\lambda_i \beta_i\alpha_i\in R_2(x,z)$ with $\lambda_i\in k,$ $\alpha_i\in Q_1(x,y_i)$ and $\beta_i\in Q_1(y_i, z)$, such that $\zeta\beta_1\notin R_2$ for some $\zeta \!\in\! Q_1(z,-)$. If $\xi \!\in \! R^{(n-1)}\m(a, x)$ with $n\m \ge \! 1$ then, for each $2\le i\le s$, there exists some $\eta_i\in kQ_n(a, y_i)$ such that $\alpha_i \xi + \eta_i\in R^{(n)}(a, y_i)$ and $\beta_i \eta_i\in R_2(-, z) \cdot kQ_{n-1}(a, -)$.

\end{Lemma}

\noindent{\it Proof.} Let $\xi \!\in \! R^{(n-1)}\m(a, x)$. If $n=1$, then we take $\eta_i=0$, for $2\le i\le s$. Assume that $n\ge 2$. Choose a $k$-basis $\{\xi_1, \ldots, \xi_t\}$ of $R^{(n-2)}(a,-)$, where $\xi_j\in R^{(n-2)}(a, b_j)$. By Lemma \ref{R-derivative}, we may write $\xi=\sum_{j=1}^t \sigma_j \, \xi_j$ for some $\sigma_j \in kQ_1(b_j, x)$.

Fix $1< i\le s$. If $\alpha_i\sigma_j\in R_2(b_i, y_l)$ for all $1\le j \le t$, then $\alpha_i \xi\in R^{(n)}(a, y_i)$; see (\ref{R-derivative}), and set $\eta_i=0$. Otherwise, let $J_i$ be the set of $j\in \{1, \ldots, t\}$ such that $\alpha_i \hspace{.4pt}\sigma_j\notin R_2(b_j, y_i)$. Fix $j\in J_i$. Since $\La$ is special multi-serial, $\sigma_j=\lambda_j\theta_j+ \delta_j$, where $\lambda_j\in k$ and $\theta_j, \delta_j\in Q_1(b_j, x)$ such that $\lambda_j\alpha_i \theta_j\notin R_2(b_j, y_i)$ and $\alpha_i\delta_j\in R_2(b_j, y_i)$. By the condition $(*)$, there exists a polynomial relation $\omega_j=\lambda_j\alpha_i\theta_j+\sum_{l=1}^{r_j}\lambda_{jl}\alpha_{il}\theta_{jl}$ in $R_2(b_j, y_i)$, where $\lambda_{jl}\in k;$ $\theta_{jl}\in Q_1(b_j, c_{jl})$ and $\alpha_{il}\in Q_1(c_{jl}, y_i)$. Again since $\La$ is special multi-serial, $\alpha_{il}\ne \alpha_i$ for $1\le l\le r_j$.
Since $\beta_i\alpha_i\notin R_2(x,z)$, we have $\beta_i \alpha_{il}\in R_2(c_{jl}, z)$ for $1\le l\le r_j$. By the induction hypothesis, we have $\eta_{jl} \in kQ_{n-1}(a, c_{jl})$ such that $\xi_{jl}=\lambda_{jl}\theta_{jl}\xi_j+\eta_{jl} \!\in\! R^{(n-1)}(a,c_{jl})$ and $\alpha_{il} \eta_{jl}$ lies in $R_2(-, y_i) \cdot kQ_{n-2}(a, -)$, for $1 \!\le \! l \!\le\! r\!_j$.
Set $\eta_i\!=\! \sum_{j\in J_i;1\le l\le r_j} \!\m \alpha_{il}\xi_{jl} \!\in \! kQ_n(a, y_i)$. Then, $\beta_i\eta_i =\sum_{j\in J;1\le l\le r_j} \alpha_i\beta_{jl}\xi_{jl}
\! \in \! R_2(-,\m z) \cdot kQ_{n-1}(a, \m -)$. Consider $\chi_i=\alpha_i\xi+\eta_i.$ Then $\chi_i=\alpha_i\xi+\sum_{j\in J_i;1\le l\le r_j} \! \alpha_{jl}\xi_{jl}$ with $\xi, \xi_{jl}\in R^{(n-1)}(a, -)$. On the other hand,  we can verify that \vspace{-3pt}
$$\textstyle\chi_i \!=\!\m \sum_{j\in I_i}\! (\omega_j + \alpha_i\delta_j) \xi_j + \sum_{j\notin I_i} \hspace{-3pt} \alpha_i\sigma_j\xi_j +\sum_{j\in J_i; 1\le l\le r_j} \!\m \alpha_{il}\eta_{jl}\in R_2(-, y_i)\cdot kQ_{n-2}(a,-).\vspace{-3pt}$$
Thus,  $\chi_i\in R^{(n)}(a, y_i)$ by Lemma \ref{R-derivative}. The proof of the lemma is completed.

\vspace{3pt}

The following is the promised new class of Koszul algebras.

\begin{Theo}\label{SKA}

Let $\La=k Q/R$ be a quadratic special multi-serial algebra with $Q$ a locally finite quiver. If the condition $(\ast)$ or its dual is satisfied, then $\La$ is Koszul.

\end{Theo}

\noindent{\it Proof.} 
By Theorem \ref{Opp-Koszul}, we only need to consider the case where the condition $(*)$ is satisfied. By Theorem \ref{Koz-proj-rls}, it amounts to show, for any $a\in Q_0$, that $\mathcal{K}_a^\pdt$ is exact in degree $-n$ for all $n\ge 1$. By Theorem \ref{qua-alg} and Proposition \ref{Koz-pres}, we may assume $n\ge 2$.
By Lemma \ref{K-cplx}, it suffices to prove that $\operatorname{Ker}\left(\partial_{a}^{-n}\right)\subseteq \operatorname{Im}\left(\partial_{a}^{-n-1}\right)$. For this purpose, we recall that $\mathcal{K}^{-n}_a=\oplus_{y \in Q_0} P_y\mla -n\mra \otimes R^{(n)}(a, y)$.

\vspace{1pt}

Consider $0\ne u\in \operatorname{Ker}(\partial_{a}^{-n})\subseteq \rad \mathcal{K}^{-n}_a$. Since $\partial_{a}^{-n}$ is graded, we may assume that
$u\in \oplus_{y\in Q_0} P_y \sla \m -n\nra\m_m(b)\otimes R^{(n)}(a, y)=\oplus_{y\in Q_0} e_b\La_{m-n}e_y\otimes R^{(n)}\m(a, \m y),$ for some $b\!\in\! Q_0$ and $m> n$.
Let $s(\hspace{.5pt} \ge 1)$ be minimal such that we can write $u=\sum_{l=1}^s \bar{\theta\hspace{1.6pt}}_{\hspace{-2pt}l} \otimes \rho_l$, for some $\theta_l\in Q_{m-n}(y_l, b)$ and $\rho_l\in R^{(n)}(a, y_l)$. Then, $\bar\theta_1, \ldots, \bar\theta_s$ are $k$-linearly independent in $\hspace{.4pt}e_b\La_{m-n}$. Choose a $k$-basis $\{\xi_1, \ldots, \xi_t\}$ of $R^{(n-1)}(a, -)$, where $\xi_j\in R^{(n-1)}(a, x_j)$. And since $\La$ is special multi-serial, $\hspace{.4pt}e_b\La_{m-n-1}$ has a $k$-basis $\{\bar\eta_1, \ldots, \bar\eta_r\}$, where $\eta_i$ is a path in $Q_{m-n-1}(z_i, b)$ with $z_i\in Q_0$.
Then, $\rho_l={\sum}_{j=1}^t \zeta_{lj} \, \xi_j$ with $\zeta_{lj} \in kQ_1(x_j, y_l);$ see (\ref{R-derivative}) and $\bar \theta_l=\sum_{i=1}^r \bar \eta_i \hspace{.4pt} \bar \delta_{il}$ with $\delta_{il}\in kQ_1(y_l, z_i)$, for $l=1, \ldots, s$. We shall divide our argument into several statements.

\vspace{1pt}

(1) {\it For any $1\le j\le t$, we have $\sum_{l=1}^s \bar\theta_l \bar\zeta_{lj}=\sum_{i=1}^r \! \sum_{l=1}^s \bar\eta_i\bar\delta_{il} \bar\zeta_{lj}=0$.}

\vspace{1pt}

Since $u=\sum_{j=1}^t \! \sum_{l=1}^s \bar\theta_l\otimes \zeta_{lj} \xi_j,$ we have $\partial_a^{-n}(u)=\textstyle\sum_{j=1}^t (\sum_{l=1}^s \bar\theta_l\bar \zeta_{lj}) \otimes \xi_j=0$; see (\m\ref{diff}\m). Since the $\xi_j$ are linearly independent, ${\sum}_{l=1}^s \bar\theta_l \bar{\zeta}_{lj} \! = \! \sum_{i=1}^r\! \sum_{l=1}^s \bar \eta_i \bar \delta_{il} \bar \zeta_{lj} \m=\m 0$, for $1\le j\le t$. This establishes Statement (1).

\vspace{1pt}

(2) {\it If $m=n+1$, then $u\in {\rm Im}(\partial_a^{-n-1})$.}

Let $m=n+1$. Then $\hspace{.4pt}e_b\La_{m-n-1}=e_b\La_0=k e_b$. In particular, $r=1$ and $\eta_1=\varepsilon_b$. By Statement (1),
$\sum_{l=1}^s \delta_{1l} \zeta_{lj}=\sum_{i=1}^r \eta_i (\sum_{l=1}^s \delta_{il} \zeta_{lj}) \m\in\! R_2(x_j, b)$, for $j=1, \ldots, t.$
Set $\chi_1=\sum_{l=1}^s \delta_{1l} \rho_l=\sum_{l=1}^s(\sum_{j=1}^t \delta_{1l}\zeta_{lj})\xi_j\in R_2(-, z_1)\cdot kQ_{n-1}(x_j, -)$. By Lemma \ref{R-derivative}(2), $\chi_1\in R^{(n+1)}(a, z_1)$ such that $\partial_a^{-n-1}(\bar\eta_1 \otimes \chi_i)=\sum_{l=1}^s\bar \theta_l\otimes \rho_l=u$. This establishes Statement (2).

Now, assume that $m\ge n+2$. Since $\delta_{il}\in kQ_1(y_l, z_i)$ and $\eta_i\in Q_{m-n-1}(z_i, b)$ is non-trivial, $\bar \theta_l=\sum_{i=1}^r \bar \eta_i \hspace{.4pt} \bar \delta_{il}$. Since $\La$ is special multi-serial, we may assume that $\delta_{il}$ is monomial, for $i=1, \ldots, r; l=1, \ldots, s$. We need to consider another derivation $\partial^\alpha: kQ\to kQ$ for any $\alpha\in Q_1$, which sends a path $\rho$ to $\eta$ if $\rho=\eta\alpha$; and $0$ if $\alpha$ is not an initial arrow of $\rho$.

\vspace{.5pt}

(3) {\it If $\delta_{il}\zeta_{lj}\notin R_2(x_j, z_i)$, then $\zeta_{lj}$ has a summand $\lambda_{lj}\alpha_{lj}$, where $\lambda_{lj}\in k$ and $\alpha_{lj}\in Q_1(x_j, y_l)$, such that $\lambda_{lj}\delta_{il}\alpha_{lj}$ is a summand of a polynomial relation in $R_2$.}

Suppose that $\delta_{ip}\zeta_{pq}\notin R_2(x_q, z_i)$ for some $1 \m\le \m i \m\le\m r; 1 \m\le \m p \m\le \m s; 1 \m\le \m q \m\le \m t$. Then,
$\delta_{ip}$ is a non-zero monomial in $kQ_1(y_{p}, z_i)$ and $\zeta_{pq}$ has a non-zero summand $\lambda_{pq}\alpha_{pq}$, where $\alpha_{pq}\in k$ and $\alpha_{pq}\in Q_1(x_q, y_p)$, such that $\lambda_{pq}\delta_{ip} \alpha_{pq}\notin R_2(x_{q}, \! z_{i})$. By Statement (1), we may write
$\textstyle\sum_{l=1}^s \theta_l \zeta_{lq} = \sum_{j=1}^d \nu_j\hspace{.4pt} \omega_j \kappa_j,$ where
$\kappa_j \in Q_{n\m_j}(x_{q}, -)$ with $n\m_j$ some non-negative integer, $\omega_j\in R_2$ and $\nu_j\in kQ_{m-n-n_j-1}(-, b)$.

Assume, for each $1 \le j \le d$, that either $n_j>0$ or $\partial^{\alpha_{pq}}(\omega_j)\!=\!0$. Applying $\partial^{\alpha_{pq}}$ to the above equation, we obtain $\sum_{l=1}^s \lambda_l \theta_l \!\in\! R_{m-n}(-, b)$, where $\lambda_l \!=\! \partial^{\alpha_{pq}}\m(\m\zeta_{lq}\m) \m\in\! k$. Since $\lambda_{p}=\lambda_{pq}\ne 0$, contrary to  $\bar\theta_1, \ldots, \bar \theta_s$ being $k$-linearly independent. Thus, we may assume that $n_1=0$ and $\alpha_{pq}$ is the initial arrow of a monomial summand of $\omega_1\in R_2(x_q, -)$. Since $\La$ is special multi-serial with $\lambda_{pq}\delta_{ip}(\alpha_{pq}\notin R_2(x_q, z_i)$, we see that $\lambda_{pq}\delta_{ip}(\alpha_{pq})$ is a summand of $\omega_1$, which is a polynomial relation in $R_2(x_q, z_i).$ This establishes Statement (3).


(4) {\it For each $1\le i\le r$, there exists some element $\chi_i \in R^{(n+1)}(a, z_i)$ such that $\partial_a^{-n-1}(\bar\eta_i\otimes \chi_i) = \sum_{l=1}^s \bar\eta_i \bar\delta_{il}\otimes \rho_l$.}

\vspace{1pt}

Fix $1\le i\le r$. If $\bar\eta_i\bar \delta_{i1}=0$ for all $l=1, \ldots, s$, then we take $\chi_i=0$. Otherwise, denote by $L$ the set of $l\in \{1, \ldots, s\}$ such that $\bar\eta_i \bar \delta_{il}\ne 0$; and for $l\in L$, denote by $J_l$ the set of $j\in \{1, \ldots, t\}$ such that $\delta_{il}\zeta_{lj}\notin R_2(x_j, z_i)$.

Fix $(j, l)\in L\times J_{l}.$
Since $\La$ is special multi-serial, we write $\zeta_{lj}=\alpha_{lj}+\sigma_{lj}$, where $\sigma_{lj}$ is such that $\delta_{il}\sigma_{lj}\in R_2(x_j, z_i)$ and $\alpha_{lj}$ is a monomial such that $\delta_{il}\alpha_{lj}\not\in R_2(x_j, z_i)$. By Statement (3), we have a polynomial relation $\omega_{lj}=\delta_{il}\alpha_{lj}+\sum_{p=1}^{r_{lj}}\gamma_{il}^p\beta_{lj}^p$ in $R_2(x_j, z_i)$, where $\beta_{lj}^p\in Q_1(x_j, c_{lj}^p)$ with $c_{lj}^p\in Q_0$ and $\gamma_{il}^p\in kQ_1(c_{lj}^p, z_i)$ is monomial. Since $\eta_i$ is a non-trivial path with $\bar \eta_i \bar \delta_{il}\ne 0$, we have $\bar \eta_i \bar \gamma_{il}^p=0$ for all $1\le p\le r_{lj}$. By Lemma \ref{msa-k-relation}, there exists $\xi_{lj}^p\in kQ_n(a, c_{lj}^p)$ such that $\rho_{lj}^p=\beta_{lj}^p\xi_j+\xi_{lj}^p\in R^{(n)}(a, c_{lj}^p)$ and $\gamma_{il}^p\xi_{lj}^p\in R_2(-, z_i) \cdot kQ_{n-1}(a,-)$, for each $1\le p\le r_{lj}$.

Put $\chi_i = \sum_{l\in L} \! \delta_{il}\rho_l+\sum_{l\in L; j\in J_l; 1\le p \le r_{lj}}\! \gamma_{il}^p \rho_{lj}^p$, where $\rho_l, \rho_{lj}^p\in R^{(n)}(a,-)$. Since $\rho_l={\sum}_{j=1}^t \zeta_{lj} \, \xi_j$, a routine verification shows that \vspace{-4pt}
$$\textstyle \chi_i \! = \! \sum_{l\in L; \hspace{.4pt} j\in J_l} \! (\omega_{lj} +  \delta_{il}\sigma_{lj}) \xi_j
+  \sum_{l\in L\m ; \hspace{.4pt} j\notin J_l} \! \delta_{il} \zeta_{lj} \xi_j
+ \sum_{\hspace{.4pt} l\in L; \hspace{.5pt} j\in J_l; \hspace{.5pt} 1\le p \le t_{lj}} \!\! \gamma_{il}^p\hspace{.5pt} \xi_{lj}^p. \vspace{-3pt} $$

Since $\delta_{il}\zeta_{lj}\in R_2(x_j, z_i)$ for $(l,j) \in L\times J_l$, we get $\chi_i\in R_2(-, z_i)\cdot kQ_{n-1}(a,-)$. By Lemma \ref{R-derivative}(2), $\chi_i\in R^{(n+1)}(a, z_i)$, and hence, $\bar\eta_i\otimes \varphi_i\in \mathcal{K}_a^{-n-1}$. Further, since $\bar\eta_i \bar \gamma_{il}^p=0$ for $(l,j)\in L\times J_l$ and $1\le p \le r_{lj}$,
and $\bar\eta_i \bar \delta_{il}=0$ for $l\notin L_i$, we deduce that
$\textstyle \partial_a^{-n-1}(\bar\eta_i\otimes \chi_i) = \sum_{l\in L} \! \bar\eta_i \bar \delta_{il}\otimes \rho_l = \sum_{l=1}^s \bar\eta_i \bar \delta_{il}\otimes \rho_l.$ This establishes Statement (4).
%

\vspace{.5pt}

Finally, $w=\sum_{i=1}^r\bar\eta_i\otimes \chi_i\in K^{-n-1}_a$ is such that
$\textstyle\partial_a^{-n-1}(w) 
=\sum_{l=1}^s \bar\theta_l\otimes \rho_l=u.$ The proof of the theorem is completed.

\vspace{3pt}

\noindent{\sc Example.} Consider the quadratic special biserial algebra $\La=kQ/R$, where \vspace{-0pt}
$$\xymatrixcolsep{28pt}\xymatrixrowsep{9pt}\xymatrix@!=4pt{
&&&3\ar[dr]^{\gamma}&&\\
Q: &1\ar[r]^\alpha &2\ar[dr]_-{\zeta}\ar[ur]^-{\beta}&&5\ar[r]^\delta & 6 \\
&&&4\ar[ur]_{\eta}&&&
} \vspace{-1pt}$$
and $R=\langle \zeta \alpha, \delta\gamma, \gamma\beta+\eta\zeta\rangle$. Clearly, $S_1$ has a linear projective $2$-presentation \vspace{-7pt}
$$\xymatrix{P_4\mla -2 \rangle \ar[r]^-{\partial_1^{-2}}\otimes k\langle \zeta\alpha \hspace{.4pt} \rangle & P_2 \mla -1  \rangle \otimes k\langle \alpha \hspace{.4pt} \rangle \ar[r]^-{\partial_1^{-1}}&P_1\otimes k \langle \varepsilon_1 \rangle \ar[r]^-{\partial_1^0}&S_1\ar[r] &0.} \vspace{-5pt}$$

Note that ${\rm Ker}(\partial_1^{-2})=k\langle \bar\delta\bar\eta\otimes \zeta\alpha\rangle$, which is generated in degree $4$. Thus, $S_1$ has no linear projective resolution. So, $\La$ is not Koszul. In fact, $\La$ does not satisfy condition $(*)$ or its dual.

\vspace{3pt}

Finally, applying the description of the linear projective resolutions and the colinear injective coresolutions for graded simple modules, we obtain a stronger version of the Extension Conjecture for certain finite dimensional Koszul algebras. 

\begin{Theo}\label{Ex-Conj-fd}

Let $\La=kQ/R$ be a finite dimensional Koszul algebra such that $\La^!$ is left or right noetherian. If ${\rm Ext}_{\mathit\Lambda}^1(S_a, S_a)\ne 0$ with $a\in Q_0$, then ${\rm Ext}_{\mathit\Lambda}^n(S_a, S_a)\ne 0$ for every integer $n\ge 1.$

\end{Theo}

\noindent{\it Proof.} Let ${\rm Ext}_{\mathit\Lambda}^1(S_a, S_a)\ne 0$ for some $a\in Q_0$. Then $Q$ has a loop $\sigma$ at $a$. Suppose first that $\La^!$ is left noetherian. Since $\La$ is finite dimensional, $\La_t e_a=0$ for some $t>0$. By Proposition \ref{Koz-cplx-dual} and Theorem \ref{Opp-Koszul}, \vspace{.5pt} $S_a^{\hspace{.5pt}!}$ has a graded injective coresolution $\mathcal{I}_a^\pdt$ with $\mathcal{I}^n_a=\oplus_{x\in Q_0} I^!_x\tla n\nra \otimes e_x \La_n e_a$, for $n\in \Z$. \vspace{.5pt} In particular, $\mathcal{I}^n_a=0$ for $n>t$.

\vspace{1pt}

Consider $M^{(i)}=\La^! (\bar\sigma^!)^i\in {\rm gmod}\La^!$ for $i\ge 1$. Since $\La^!$ is left noetherian, by Proposition \ref{p-cover}, $M^{(i)}$ has a minimal graded projective resolution \vspace{-4pt}
$$\xymatrixcolsep{20pt}\xymatrix{P^{\pdt, i}: \quad \cdots\ar[r] & P^{-n,i} \ar[r] & \cdots \ar[r] & P^{-1,i}\ar[r] & P^{\hspace{.5pt}0,i}\ar[r] & 0}
\vspace{-4pt}$$ over
${\rm gproj}\La^!$. For any $s\in \Z$, we have ${\rm GHom}_{\mathit\Lambda^!}(P^{-n,i}\hspace{-2pt}, S_a^{\hspace{.5pt}!}\nla s \sra)\cong {\rm GExt}^n_{\m\mathit\Lambda^!}(M^{(i)}, S_a^{\hspace{.5pt}!}\sla s \sra)$, which is a subquotient of ${\rm GHom}_{\mathit\Lambda^!}(M^{(i)}, \mathcal{I}^n_a\sla s \sra).$ Thus ${\rm GHom}_{\mathit\Lambda^!}(P^{-n,i}, S_a^{\hspace{.5pt}!}\sla s \sra)=0$. That is, $P_a^!\sla s \sra$ is not a direct summand of $P^{-n,i}$ for all $n>t$ and all $s\in \Z$. Forgetting the gradation, we see that $P^{\pdt,i}$ is a projective resolution of $M^{(i)}$ over ${\rm proj}\La^!$ such that $P^!_a$ is not direct summand of $P^{-n,i}$ for all $n>t$. In other words, $P^{\pdt,\hspace{.5pt} i}$ is an $e_a$-bounded projective resolution of $\La^! (\bar\sigma^{\hspace{.5pt}!})^i$ over ${\rm proj}\La^!$; see, for definition,  \cite{ILP}, for every $i\ge 1$.
If $\bar\sigma^{\hspace{.5pt}!}$ is nilpotent, as argued in the proof of \cite[(1.6)]{ILP}, we conclude that $\bar\sigma^{\hspace{.5pt}!}\in [\La^!\!, \La^!]+
\sum_{a\ne x\in Q_0} \!\La^! e_x\La^!$, \vspace{.5pt} where $[\La^!, \La^!]$ is the commutator group of $\La^!$. This implies that $\bar\sigma^{\hspace{.5pt}!}\!\!\in {\rm rad}^2\!\m\La^!$, contrary to $R^!$ being generated in degree $2$. Thus, $0\ne (\bar\sigma^{\hspace{.5pt}!})^n\in e_a\La^!_ne_a$ for every $n\ge 0$. By Lemma \ref{k-cplx-iso}, $S_a$ admits a linear projective resolution $\mathcal{P}_a^\pdt$, which is defined by $\mathcal{P}_a^{-n}=\oplus_{x\in Q_0} P_x\tla -n\nra \otimes D(e_a \La^!_n e_x)$ for $n\ge 1$. \vspace{1pt} Since $e_a \La^!_n e_a\ne 0$, we see that $P_a\tla -n\nra$ is a direct summand of $\mathcal{P}_a^n$, for $n\ge 1$. \vspace{1pt} That is, ${\rm GExt}_{\m\mathit\Lambda}^n(S_a, S_a\sla -n\sra)\ne 0$, and hence, ${\rm Ext}_{\m\mathit\Lambda}^n(S_a, S_a)\ne 0$ for $n\ge 1$.

Suppose now that $\La^!$ is right noetherian. Then, $(\La^!)^{\rm o}$ is left noetherian. By Theorem \ref{Opp-Koszul} and Proposition \ref{q-dual}$,\La^{\rm o}$ is a finite dimensional Koszul algebra with $(\La^{\rm o})^!=(\La^!)^{\rm o}.$ Since $\sigma^{\rm o}$ is a loop in $Q^{\rm o}$ at $a$, \vspace{1pt} ${\rm Ext}_{\mathit\Lambda^{\rm o}}^n(S_a^{\rm o}, S_a^{\rm o})\ne 0$, for every $n\ge 1.$ Since $\La$ is finite dimensional, we have a duality
$\mk D: {\rm mod}^{\hspace{.5pt}b\hspace{-2.8pt}}\La^{\rm o}\to {\rm mod}^{\hspace{.5pt}b\hspace{-2.8pt}}\La.$ Therefore, ${\rm Ext}_{\mathit\Lambda}^n(S_a, S_a)\ne 0$, for every $n\ge 1.$ The proof of the theorem is completed.

\vspace{4pt}

\noindent{\sc Remark.} Note that special multiserial algebras defined by finite quivers are left and right noetherian; see \cite{VHW, LLi}. Thus, the Extension Conjecture holds for finite dimensional Koszul algebras with a special multiserial Koszul dual.

\vspace{4pt}

\noindent{\sc Example.} Consider the special multi-serial algebra $\La=kQ/R$, where \vspace{-5pt}
$$\begin{tikzpicture}[-{Stealth[inset=0pt,length=3pt,angle'=35,round,bend]},scale=0.33]
\node (b) at (3.3,0){$1$};
\draw (3.65,0.3) arc (160:-160:1);
\node[left] at (6.9,0){$\alpha$};
\node (a) at (0,0) {$2$};
\draw (0.6,0.2) -- (2.9,0.2);
\node[above] at (1.67, 0.05){$\beta$};
\draw (2.92,-0.2) -- (0.52,-0.2);
\node[above] at (1.6,-1.7){$\gamma$};
\draw (-0.45,0.25)--(-2.8,0.2);
\node[above] at (-1.6,0.1) {$\delta$};
\node (c) at (-3.2,0) {$3$};
\draw (-5.84,-0.2)-- (-3.6,-0.2);
\node[above] at (-4.59,-1.8){$\zeta$};
\draw  (-5.84,0.2)--(-3.6,0.2);
\node[above] at (-4.5,0.15){$\eta$};
\node (d) at (-6.15,0) {$4$};
\draw (-6.55,0.3) arc (25:340:1);
\node[above] at (-8.9,-0.7){$\sigma$};
\node[above] at (-13,-1){$Q:$};
\end{tikzpicture}\vspace{-7pt}$$
and $R=\langle \alpha^2+\beta\gamma, \alpha\beta, \gamma\beta, \gamma \alpha, \eta\sigma, \sigma^2\rangle$. Satisfying the condition $(*)$, $\La$ is Koszul; see (\ref{SKA}). Note that $\La^!=kQ^{\rm o}/R^!$, where $R^!=\langle (\alpha^{\rm o})^2 - \gamma^{\rm o}\beta^{\rm o}, \gamma^{\rm o}\delta^{\rm o}, \sigma^{\rm o}\zeta^{\rm o}\rangle,$ which is noetherian. By Theorem \ref{Ex-Conj-fd}, $\Ext^i_{\!\m\mathit\Lambda}(S_1, S_1)\ne 0$ and $\Ext^i_{\!\m\mathit\Lambda}(S_4, S_4)\ne 0$ for $i\ge 1$.



\vspace{-1pt}

\section{Double complexes and extension of functors}

The main objective of this section is to formalize a technique of extending a functor from an additive category to a complex category to the complex category, which has been already used in various settings; see, for example, \cite{BaL2, BGS, 
Ric}. This technique is essential for us to construct derived Koszul functors. An additive $k$-category is called {\it concrete} if the objects are equipped with a $k$-vector space structure and morphisms are $k$-linear maps.

Throuout this section, $\cA, \cB, \cC$ stand for concrete additive $k$-categories. Let $(X^\ydt, d\hspace{-1.5pt}_X^{\,\pdt})$ be a complex over $\cA$. Given $n\in \Z$, we shall write $X^\cdt[n]$ for the {\it $n$-shift} of $X^\cdt$ and ${\rm H}^n(X^\cdt)$ for its $n$-th homology group. The {\it twist complex} $\mathfrak{t}(X^\cdt)$ of $X^\cdt$ is defined by $\mathfrak{t}(X^\cdt)^n=X^n$ and $d_{\mathfrak{t}(X^\cdt)}^n=-d^n_{\m X}$; see \cite{BaL2}.
This induces an automorphism $\mathfrak{t}$ of $C(\cA)$, called the {\it twist functor}. And for a morphism $f^\ydt: X^\ydt\to Y^\ydt$ in $C(\cA)$, we shall write $C_{\hspace{-1pt}f^\pdt}$ for its {\it mapping cone}; see \cite[(III.1.5)]{Mil}.

First, we develop a homotopy theory of double complexes. Let $(M^{\ydt\hspace{.5pt} \ydt}, v_{\hspace{-1pt}_M}^{\ydt\hspace{.5pt} \ydt}, h_{\hspace{-1pt}_M}^{\ydt\hspace{.5pt} \ydt})$ be a double complex  \vspace{1pt} over $\cA$, where $v_{\hspace{-1pt}_M}^{\ydt\hspace{.5pt} \ydt}$ is the vertical differential and $h_{\hspace{-1pt}_M}^{\ydt\hspace{.5pt} \ydt}$ is the horizontal one. Given $i, j\in \Z$, we call $(M^{i, \ydt}, v_{\hspace{-1pt}_M}^{i,\,\ydt\,})$ and $(M^{\cdt, \hspace{.5pt} j}, h_{\hspace{-1pt}_M}^{\cdt, \hspace{.5pt} j})$ the {\it $i$-th column} and the {\it $j$-th row} of $M^{\ydt\hspace{.5pt} \ydt}$, respectively. A double complex morphism $f^{\cdt \cdt}\!:\m M^{\ydt\hspace{.5pt}\ydt}\to N^{\ydt\hspace{.5pt}\ydt}$ consists of morphisms $f^{i,j}: M^{i,j}\to N^{i,j}$ in $\cA$ with $i, j\in \Z$ making \vspace{-1pt}  
$$
\xymatrixrowsep{15pt}\xymatrixcolsep{12pt}\xymatrix@=8pt{
&N^{i,j+1}\\
M^{i,j+1} \ar[ur]^{f^{i,j+1}}&  \\
& N^{i,j} \ar[rr]^--{h_{\hspace{-1pt}_N}^{i,j}}\ar[uu]_{v_{\hspace{-1pt}_N}^{i,j}}&& N^{i+1,j} \\
M^{i,j} \ar[uu]^{v_{\hspace{-1pt}_M}^{i,j}} \ar[ur]^-{f^{i,j}} \ar[rr]^--{\hspace{5pt}h_{\hspace{-1pt}_M}^{i,j}} && M^{i+1,j}\ar[ur]_-{f^{i+1,j}}
}\vspace{-1pt}$$ commute for all $i, j\in \Z$, that is, $f^{i,\ydt}: M^{i,\ydt}\to N^{i,\ydt}$ and $f^{\ydt,j}: M^{\cdt,j}\to N^{\cdt,j}$ are complex morphisms, for all $i,j\in \Z$. The double complexes over $\cA$ together these morphisms form an additive $k$-category written as $DC(\cA)$.

Now, assume that $\cA$ has countable direct sums. Given $M^{\cdt \cdt}\in DC(\cA)$, its {\it total complex} $\mathbb{T}(M^{\ydt\hspace{.6pt} \ydt})\in C(\cA)$ is defined by $\T(M^{\ydt\hspace{.6pt} \ydt})^n=\oplus_{i\in \Z} \, M^{i, n-i}$ and
\vspace{-3pt} $$d_{\hspace{.6pt} \T(M^{\ydt\hspace{.6pt} \ydt})}^{\,n}=(d_{\hspace{.6pt}\T(M^{\ydt\hspace{.6pt} \ydt})}^{\,n}(j,i))_{(j,i)\in \Z\times \Z}:
\oplus_{i\in \mathbb{Z}} \, M^{i, n-i} \to \oplus_{j\in \mathbb{Z}} \, M^{j, n+1-j}, \vspace{-3pt}$$
where $d_{\hspace{.6pt}\T(M^{\ydt\hspace{.6pt} \ydt})}^{\,n}(j,i)\!:\m  M^{i, n-i} \to M^{j, n+1-j}$ is defined such that $d_{\hspace{.6pt}\T(M^{\ydt\hspace{.6pt}\ydt})}^{\,n}(i,i) = v_{\hspace{-1pt}_M}^{i, n-i};$ $d_{\T(M^{\ydt\hspace{.6pt}\ydt})}^{\,n}(i+1,i)=h_{\hspace{-1pt}_M}^{i, n-i}$ and $d_{\hspace{.6pt}\T(M^{\ydt\hspace{.6pt}\ydt})}^{\,n}(j,i) = 0$ if $j\notin i$ or $i+1$. And given a morphism $f^{\pdt\hspace{.4pt}\pdt}\m:\! M^{\ydt\hspace{.6pt}\ydt} \!\to\! N^{\ydt\hspace{.6pt}\ydt}$ in $DC(\cA)$, we define its {\it total morphism} $\T(f^{\ydt \hspace{.4pt} \ydt})\m:\! \T(M^{\ydt\hspace{.5pt} \ydt}) \!\to \!\T(N^{\ydt\hspace{.5pt} \ydt})$ by \vspace{-3pt} $$\T(f^{\ydt\hspace{.4pt} \ydt})^n=\left(\T\left(f^{\ydt \hspace{.5pt} \ydt}\right)^n\hspace{-3pt}(j, i)\right)_{(j,i)\in \mathbb{Z}\times \mathbb{Z}}: \oplus_{i\in \Z} \, M^{i, n-i} \to \oplus_{j\in \Z} \,  N^{j, n-j},
\vspace{-3pt}$$
where $\T(f^{\ydt\hspace{.4pt}\ydt})^n(j, i): M^{i, n-i} \to N^{j, n-j}$ is such that $\T(f^{\ydt\hspace{.4pt}\ydt})^n(i, i)=f^{i, n-i}$ and $\T(f^{\ydt\hspace{.4pt}\ydt})^n(j, i)=0$ for all $j\ne i$.
%
%
This yields clearly a functor $\mathbb{T}: DC(\cA) \to C(\cA)$.

\vspace{2pt}

We shall study when the total complex of a double complex is acyclic. Consider a double complex $M^{\pdt \hspace{.4pt} \pdt}\in DC(\cA)$. Given $n\in \Z$, the $n$-{\it diagonal} of $M^{\pdt \hspace{.4pt} \pdt}$ consists of the objects $M^{i, n-i}$ with $i\in \Z$. We shall say that $M^{\pdt \hspace{.4pt} \pdt}$ is {\it $n$-diagonally bounded} (respectively, {\it bounded-above}, {\it bounded-below}) if $M^{i, n-i}=0$ for all but finitely many (respectively, positive, negative) integers $i$. Moreover, $M^{\pdt \hspace{.4pt} \pdt}$ is called {\it diagonally bounded} (respectively, {\it bounded-above},  {\it bounded-below}) if it is $n$-diagonally bounded (respectively, bounded-above, bounded-below) for every $n\in \Z$. 
We obtain a local version of the well-known Acyclic Assembly Lemma; see \cite[(2.7.3)]{Wei} as follows.

\begin{Prop}\label{Homology-zero}

Let $\cA$ be a concrete additive $k$-category with countable direct sums. If $M^{\ydt \hspace{.6pt} \ydt}\in DC(\cA)$ and $n\in \Z$, then ${\rm H}^n(\T(M^{\ydt \hspace{.4pt} \ydt}))=0$ in case

\vspace{-2pt}

\begin{enumerate}[$(1)$]

\item $M^{\ydt \hspace{.4pt} \ydt}$ is $n$-diagonally bounded-below with ${\rm H}^{n-j}(M^{\ydt\hspace{.7pt}, \hspace{.5pt} j})=0$ for all $j\in \Z;$ or

\vspace{1pt}

\item $M^{\ydt \hspace{.4pt} \ydt}$ is $n$-diagonally bounded-above with ${\rm H}^{n-i}(M^{i, \ydt})=0$ for all $i\in \Z.$

\end{enumerate} \end{Prop}

\noindent{\it Proof.} Let $(M^{\cdt \hspace{.4pt}\cdt}\hspace{-3pt}, v^{\cdt \cdt}\hspace{-3pt}, h^{\cdt \cdt}) \!\in\!\m DC(\cA)$ such that Statement (1) holds, say $\! M^{i, n-i}\!\!=\!0,$  for all $i<t$, where $t$ is some negative integer. Consider $c=(c_{i, n-i})_{i\in \mathbb{Z}}\in {\rm Ker}(d^n_{\hspace{.6pt}\mathbb{T}(M^{\cdt \cdt})})$ with $c_{i, n-i} \!\in\! M^{i, n-i}.$ Then, $v^{i, n-i}(c_{i, n-i})+h^{i-1, n-i+1}(c_{i-1, n-i+1})=0$, for $i\in \Z$. We may assume that $c_{i, n-i}=0$ for $i> 0$. Then, $h^{0,n}(c_{0,n})=-v^{1, n-1}(c_{1, n-1})=0$. Since ${\rm H}^0(M^{\ydt\hspace{.5pt}, n})=0$, we have $c_{0,n}=h^{-1,n}(x_{-1,n})$, for some $x_{-1,n}\in M^{-1, n}$. So \vspace{-2pt}
$$h^{-1, n+1}(c_{-1, n+1} - v^{-1, n}(x_{-1,n})) = h^{-1, n+1}(c_{-1, n+1}) + v^{0, n}(c_{0,n})=0.\vspace{-2pt}$$
Since ${\rm H}^{-1}(M^{\ydt\hspace{.6pt},n+1})=0$, we see that $c_{-1, n+1}-v^{-1, n}(x_{-1,n})=h^{-2,n+1}(x_{-2,n+1})$ for some
$x_{-2,n+1}\in M^{-2, n+1}$. Continuing this process, we get $x_{i, n-1-i}\in M^{i, n-1-i}$ such that
$c_{i, n-i}=v^{i, n-1-i}(x_{i,n-1-i})+h^{i-1,n-i}(x_{i-1,n-i})$, for $i=-1, -2, \ldots, t$. Since $M^{t-1, n-t+1}=0$, we see that $v^{t-1, n-t}(x_{t-1, n-t})=0=c_{t-1, n-1+1}$. Setting $x=(x_{i, n-1-i})_{i\in \mathbb{Z}}$ with $x_{i, n-1-i}=0$ for $i\geq 0$ or $i<t-1$, we obtain $c=d^{n-1}_{\hspace{.6pt}\mathbb{T}(M^{\cdt \cdt})}(x)$. The proof of the proposition is completed.

\vspace{3pt}

Given a double complex $(M^{\ydt\hspace{.5pt} \ydt}, v_{\hspace{-1pt}_M}^{\ydt\hspace{.5pt} \ydt}, h_{\hspace{-1pt}_M}^{\ydt\hspace{.5pt} \ydt})$, \vspace{.5pt} we define its {\it horizontal shift} $M^{\ydt\hspace{.5pt} \ydt}[1]$ to be $(X^{\ydt\hspace{.5pt}\ydt}, v_{\hspace{-1pt}_X}^{\ydt\hspace{.5pt}\ydt}, h_{\hspace{-1pt}_X}^{\ydt\hspace{.5pt}\ydt})$ with $X^{i,j}=M^{i+1,j}$, $v_{\hspace{-1pt}_X}^{i, j}=-v_{\hspace{-1pt}_M}^{i+1,j}$ and $h_{\hspace{-1pt}_X}^{i,j}=-h_{\hspace{-1pt}_M}^{i+1, j}$, \vspace{.5pt}  for all $i, j\in \Z$.
Moreover, \vspace{.5pt}  a morphism $f^{\cdt \cdt}: M^{\cdt \cdt}\to N^{\cdt \cdt}$ is called {\it horizontally null-homotopic} if
there exist morphisms $u^{i, j}: M^{i,j}\to N^{i-1, j}$ such that $u^{i+1, j} \circ h_{\hspace{-1pt}_M}^{i,\hspace{.5pt} j}+ h_{\hspace{-1pt}_N}^{i-1,\hspace{.5pt} j}\circ u^{i,j}=f^{i,j}$ and
$v_{\hspace{-1pt}_N}^{i-1, j} \circ  u^{i,j}+u^{i,j+1} \circ  v_{\hspace{-1pt}_M}^{i, j} =0$, for all  $i, j\in \Z$.

\vspace{-.5pt}

\begin{Lemma}\label{Tot}

Let $\cA$ be a concrete additive $k$-category with countable direct sums.

\begin{enumerate}[$(1)$]

\vspace{-1pt}

\item If $M^{\ydt\hspace{.5pt} \ydt}\in DC(\cA)$, then $\T(M^{\ydt\hspace{.5pt} \ydt}[1])=\T(M^{\ydt\hspace{.5pt} \ydt})[1].$

\vspace{.5pt}

\item If $f^{\ydt \pdt}: M^{\ydt\,\ydt}\to N^{\ydt\,\ydt}$ is horizontally null-homotopic, then $\T(f^{\ydt \pdt})$ is null-homotopic.

\end{enumerate}\end{Lemma}

\noindent{\it Proof.} Statement (1) can be shown by a routine verification. Let $f^{\pdt\hspace{.4pt} \pdt}: M^{\ydt\hspace{.6pt}\ydt}\to N^{\ydt\hspace{.6pt}\ydt}$ be horizontally null-homotopic. Then, there exist $u^{i, j}: M^{i,j}\to N^{i-1, j}$ such that $f^{i,j}=u^{i+1, j} \circ h_{\hspace{-1pt}_M}^{i,\hspace{.5pt} j}+ h_{\hspace{-1pt}_N}^{i-1,\hspace{.5pt} j} \circ u^{i,j}$  \vspace{.5pt}  and $v_{\hspace{-1pt}_N}^{i-1, j} u^{i,j}+ u^{i,j+1} v_{\hspace{-1pt}_M}^{i, j} =0$, for all $i,j\in \Z$.

Given $n\in \Z$, set $h^n=(h^n(j,i))_{(j,i)\in \Z\times \Z}: \oplus_{i\in \Z} M^{i, n-i}\to \oplus_{j\in \Z} N^{j, n-j},$ where
$h^n(n-i,i)=u^{i,n-i}$ and $h^n(j,i)=0$ for $j\ne n-i$.
It is easy to verify that $\T(f^{\ydt \hspace{.5pt} \ydt})^n=h^{n+1}\circ d^n_{\hspace{.6pt}\T(M^{\ydt \hspace{.5pt} \ydt})}+ d^{n-1}_{\hspace{.6pt}\T(N^{\ydt \hspace{.5pt} \ydt\hspace{.5pt}})} \circ h^n.$ The proof of the lemma is completed.

\vspace{3pt}

Given a morphism $f^{\cdt\hspace{.3pt}\cdt}: M^{\cdt\hspace{.3pt}\cdt}\to N^{\cdt\hspace{.3pt}\cdt}$ in $DC(\cA)$. We define its {\it horizontal cone} $H_{\hspace{-1.3pt}f^{\cdt\hspace{.3pt}\cdt}}$ to be the double complex $(H^{\cdt\hspace{.3pt}\cdt},v^{\cdt\hspace{.3pt}\cdt},h^{\cdt\hspace{.3pt}\cdt})$ with
$H^{i,j}=M^{i+1,j}\oplus N^{i,j}$ and \vspace{-1pt}
$$
v^{i,j}=\left(\begin{array}{cc}
-v_{\hspace{-1pt}_M}^{i+1,j} &  0  \vspace{1.5pt} \\
         0   & v_{\hspace{-1pt}_N}^{i,j}
\end{array}\right), \;
h^{i,j}=\left(
\begin{array}{cc}
-h_{\hspace{-1pt}_M}^{i+1,j} &  0 \vspace{1.5pt} \\
f^{\,i+1,j} &  h_{\hspace{-1pt}_N}^{i,j}
\end{array}\right).
\vspace{-1pt}
$$

\vspace{-0pt}

Observe that the $j$-th row of $H_{\hspace{-1.3pt}f^{\ydt \pdt}}$ is the mapping cone of $f^{\ydt\hspace{.3pt},\hspace{.5pt} j}: M^{\ydt, \hspace{.5pt}j}\to N^{\ydt, \hspace{.5pt}j}$. In a similar fashion, we may define the {\it vertical cone} $V_{\hspace{-.5pt}f^{\ydt \pdt}}$ of $f^{\ydt \pdt}$ in such a way that its $i$-th column is the mapping cone of $f^{i,\ydt}: M^{i,\ydt}\to N^{i, \ydt}$. By a routine verification, we can verify the following statement.

%
%
%
%
%
%
%
%


%


\begin{Lemma}\label{Cone}

Let $\cA$ be a concrete additive $k$-category having countable direct sums. If $f^{\cdt\hspace{.3pt}\cdt}: M^{\cdt\hspace{.3pt}\cdt}\to N^{\cdt\hspace{.3pt}\cdt}$ is a morphism in $DC(\cA)$, then $\T(H_{\hspace{-1.3pt}f^{\cdt\hspace{.3pt}\cdt}})= C_{\hspace{.4pt}\mathbb{T}(f^{\cdt\hspace{.3pt}\cdt})} =\T(V_{\! f^{\cdt\hspace{.3pt}\cdt}}).$

\end{Lemma}

\vspace{2pt}

The following statement tells us when the total morphism of a double complex morphism is a quasi-isomorphism.

\begin{Lemma}\label{V-cone}

Let $\cA$ be a concrete additive category having countable direct sums. Consider a morphism $f^{\ydt \pdt}: M^{\ydt\hspace{.6pt}\ydt}\to N^{\ydt\,\ydt}$ in $DC(\cA)$ such that $f^{i,\cdt}\hspace{-2pt}:\hspace{-1.5pt} M^{i,\cdt} \hspace{-2.2pt} \to \hspace{-2.5pt} N^{i,\cdt}$ is a quasi-isomorphism, for every $i\in \Z$. If $M^{\ydt\hspace{.6pt}\ydt}$ and $N^{\ydt\,\ydt}$ are diagonally bounded-above, then $\T(f^{\ydt \pdt})$ is a quasi-isomorphism.

%
%
%
%

\end{Lemma}

\noindent{\it Proof.} Assume that $M^{\ydt\hspace{.6pt}\ydt}$ and $N^{\ydt\,\ydt}$ are diagonally bounded-above. Then, the vertical cone $V_{f^{\ydt \pdt}}$ of $f^{\ydt \pdt}$ is also diagonally bounded-above. Given $i\in \Z$, since $f^{i,\cdt}\hspace{-2pt}:\hspace{-1.5pt} M^{i,\cdt} \hspace{-2.2pt} \to \hspace{-2.5pt} N^{i,\cdt}$ is a quasi-isomorphism, its cone is acyclic, that is, the $i$-th column of $V_{\! f^{\cdt \cdt}}$ is acyclic.
By Proposition \ref{Homology-zero}, $\T(V_{f^{\cdt \cdt}})$ is acyclic. By Lemma \ref{Cone}, $C_{\hspace{.5pt}\T(f^{\cdt \cdt})}$ is acyclic, and hence, $\T(f^{\cdt \cdt})$ is a quasi-isomorphism. The proof of the lemma is completed.

\vspace{2pt}

Now, suppose that $\cB$ has countable direct sums. Consider a functor \vspace{-1pt}
$$F: \cA\to C(\cB): M\to F(M)^\ydt\hspace{.4pt}; f\mapsto F(f)^\ydt.
\vspace{-1pt}$$
In order to extend $F$ to $C(\cA)$, we first construct a functor
$F^{DC}: C(\cA)\to DC(\cB)$. Given a complex $M^\ydt\in C(\cA)$, applying $F$ component-wise yields a double complex \vspace{-13pt}
$$\hspace{-30pt}\xymatrixrowsep{20pt}\xymatrixcolsep{18pt}\xymatrix{
                             &                            & \vdots        & & \vdots        \\
F(M^\cdt)^\cdt: \hspace{5pt} & \cdots\ar[r] & F(M^i)^{j+1} \ar[u] \ar[rr]^{F(d_{\hspace{-.5pt}M}^i)^{j+1}} && F(M^{i+1})^{j+1} \ar[u]\ar[r] & \cdots\\
& \cdots\ar[r] & F(M^i)^j \ar[rr]^{F(d_{\hspace{-.5pt}M}^{\hspace{.5pt}i})^j} \ar[u]^{(-1)^i d^j_{F(M^i)}} && F(M^{i+1})^j \ar[u]_{(-1)^{i+1} d_{F(M^{i+1})}^j} \ar[r] & \cdots  \\
&& \ar[u] && \ar[u] }\vspace{-16pt}$$
$$\hspace{-30pt}\xymatrixrowsep{20pt}\xymatrixcolsep{18pt}\xymatrix{
&\hspace{11pt}& \vdots   & \hspace{23pt} && \vdots } \vspace{1.5pt} $$
whose $i$-th column is $\mathfrak{t}^i(F(M^i)^\cdt)$, that is the $i$-th twist of $F(M^i)^\cdt$. And given a morphism $f^\ydt: M^\ydt\to N^\ydt$ in $C(\cA)$, we obtain a commutative diagram \vspace{-1pt}
$$\xymatrixcolsep{14pt}\xymatrixrowsep{8pt}\xymatrix{
& F(N^i)^{j+1}\\
F(M^i)^{j+1} \ar[ur]^{F(f^i)^{j+1}}&  \\
& F(N^i)^j \ar[rr]^{F(d_N^i)^j}\ar[uu]_{(-1)^id^j_{F(N^i)}}&& F(N^{i+1})^j,\\
F(M^i)^j \ar[uu]^{(-1)^id^j_{F(M^i)}}\ar[ur]^-{F(f^i)^j} \ar[rr]^-{F(d_M^i)^j} && F(M^{i+1})^j\ar[ur]_-{F(f^{i+1)^j}}
} \vspace{-1pt} $$ for $i, j\in \Z$. So, $F(f^\ydt)^\ydt=(F(f^i)^j)_{i,j\in \mathbb Z}: F(M^\ydt)^\ydt \to F(N^\ydt)^\ydt$ is a morphism in $DC(\cB)$.

\vspace{1pt}

\begin{Prop}\label{cal-F}

Let $\cA, \cB$ be concrete additive $k$-categories such that $\cB$ has coun\-table direct sums. Then every functor $F: \cA\to C(\cB)$ induces a functor \vspace{-1pt}$$F^{DC}: C(\cA) \to DC(\cB): M^\ydt\mapsto F(M^\ydt)^\ydt\hspace{.5pt}; \hspace{.5pt} f^\ydt\mapsto F(f^\ydt)^\ydt.
\vspace{-2pt}
$$

\begin{enumerate}[$(1)$]

\item If $M^{\cdt}$ is a complex in $C(\cA)$, then $F^{DC}(M^{\cdt}) = F^{DC}(M^{\cdt})[1].$

\item If $f^{\ydt}$ is a morphism in $C(\cA)$, then $F^{DC}(C_{\hspace{-1pt}f^{\ydt}})=H_{F^{DC}(f^{\pdt})};$ and in case $f^{\ydt}$ is null-homotopic, $F^{DC}(f^{\ydt})$ is horizontally null-homotopic.

\end{enumerate}

\end{Prop}

\noindent{\it Proof.} Statement (1) and the first part of Statement (2) can be shown by a routine verification. Let $f^\ydt: M^\ydt\to N^\cdt$ be a null-homotopic morphism in $C(\cA)$
with morphisms $u^i: M^i\to N^{i-1}$ in $\cA$ such that $f^i=u^{i+1}\circ d_{\hspace{-1pt}_M}^i+d_{\hspace{-1pt}_N}^{i-1}\circ u^i$, for all $i\in \Z$. Therefore, \vspace{.5pt} $F(f^i)^j=F(u^{i+1})^j\circ F(d_{\hspace{-1pt}_M}^i)^j+F(d_{\hspace{-1pt}_N}^{i-1})^j\circ F(u^i)^j,$ for all $j\in \Z$. Since $F(u^i)^\cdt: F(M^i)^\cdt\to F(N^{i-1})^\cdt$ is a complex morphism, we see that $(-1)^i F(u^i)^{j+1}\circ d_{F(M^i)}^j  + (-1)^i d_{F(N^{i-1})}^j \circ F(u^i)^j=0,$ for all $j\in \Z$.
That is, $F^{DC}\hspace{-1.5pt}(f^{\ydt})$ is horizontally null-homotopic.
The proof of the proposition is completed.

\vspace{2pt}

Originally formulated for module categories in \cite[(3.7)]{BaL2}, the following statement can be routinely verified using Lemmas \ref{Tot} and \ref{Cone} and Proposition \ref{cal-F}.

\begin{Prop}\label{F-extension}

Let $\cA$ and $\cB$ be concrete additive $k$-categories such that $\cB$ has countable direct sums. Then, every functor $F: \cA\to C(\cB)$ \vspace{.5pt} extends to a functor
$F^{\hspace{.5pt}C}=\T\circ F^{DC}\hspace{-1.8pt}: C(\cA) \to C(\cB)$ with the following properties.

\vspace{-1pt}

\begin{enumerate}[$(1)$]

\item If $M$ is an object in $\cA$, then $F^{\hspace{.5pt}C}\hspace{-1.5pt}(M) = F(M).$

\item If $M^\ydt$ is a complex in $C(\cA)$, then $F^{\hspace{.5pt}C}\hspace{-1pt}(M^\ydt[1]) = F^{\hspace{.5pt}C}\hspace{-1pt}(M^\ydt)[1].$

\item If $f^\ydt$ is a morphism in $C(\cA)$, \vspace{-1pt} then $F^{\hspace{.5pt}C}\hspace{-1.4pt}(C_{\hspace{-1pt}f^\pdt})=C_{F^C\hspace{-1.4pt}(f^\pdt)};$ and in case $f^\ydt$ is null-homotopic, $F^{\hspace{.5pt}C}\hspace{-1.4pt}(f^\ydt)$ is null-homotopic.

\end{enumerate} \end{Prop}

\vspace{2pt}

As shown below, this extension of functors is compatible with the composition of functors. This is essential for our later investigation.

\begin{Prop}\label{F-composition}

Let $\cA, \cB$ and $\cC$ be concrete additive $k$-categories such that $\cB$ and $\cC$ have countable direct sums. If $F: \cA\to C(\cB)$ and $G: \cB\to C(\cC)$ are functors, then $(G^{\hspace{.5pt}C} \circ F)^{\hspace{.3pt}C} =G^{\hspace{.5pt}C} \circ F^{\hspace{.5pt}C}.$

\end{Prop}

\noindent{\it Proof.} Consider functors $F: \cA\to C(\cB)$ and $G: \cB\to C(\cC)$. We obtain a composite functor
$G^{\hspace{.5pt}C}  \circ F : \cA \to C(\cC)$, which extends to a functor $(G^{\hspace{.5pt}C} \circ F)^{\hspace{.5pt}C} : C(\cA) \to C(\cC)$. Fix $M^\ydt\in C(\cA)$ and $n\in \Z$. By definition, we have \vspace{-2pt}
$$(G^{\hspace{.5pt}C} \circ F)^C\m(M^\cdt)^n = \oplus_{i\in \Z} \hspace{.5pt} G^{\hspace{.5pt}C}\m(F(M^i)^\ydt)^{n-i} = \oplus_{(i,p)\in \mathbb{Z} \m \times \m \mathbb{Z}}\hspace{.8pt} G(F(M^i)^p)^{n-i-p}. \vspace{-2pt}
$$
Writing $(G^{\hspace{.3pt}C} \hspace{-2pt} \circ\hspace{-1pt} F)^C\!(M^\cdt)^{n+1}= \oplus_{(j,q)\in \mathbb{Z}\times\mathbb{Z}} G(F(M^j)^q)^{n+1-j-q}$, \vspace{.5pt}  we can routinely verify that the differential
$d^n_{(G^{\hspace{.5pt}C} \circ F)^{\hspace{.5pt}C}\m\left(M^\ydt\right)}: (G^{\hspace{.5pt}C}  \circ F)^{\hspace{.5pt}C}\m (M^\cdt)^n \to (G^{\hspace{.5pt}C}  \circ F)^C\m(M^\cdt)^{n+1}$ \vspace{.5pt} is given by the matrix $(d^n_{(G^{\hspace{.5pt}C}\circ F)^{\hspace{.5pt}C}(M^\cdt)}(j,q;i,p))_{(j,q;i,p) \in \Z^4},$ where \vspace{-3pt}
$$d^n_{( G^{\hspace{.5pt}C}\circ F)^{\hspace{.5pt}C}(M^\ydt)}(j,q;i,p):
G(F(M^i)^p)^{n-i-p} \to G(F(M^j)^q)^{n+1-j-q}\vspace{-1pt}$$ is defined by \vspace{-8pt}
$$
d^{\hspace{.5pt}n}_{(G^{\hspace{.5pt}C}\circ  F)^{\hspace{.5pt}C}(M^\ydt)}(j,q;i,p)
=\left\{\hspace{-3pt} \begin{array}{ll}
(-1)^{i+p} d_{G(F(M^i)^p)}^{\hspace{.5pt} n-i-p}, & j=i; q=p\hspace{.4pt};  \vspace{2pt}  \\
(-1)^i G(d_{F(M^i)}^{\hspace{.5pt}p})^{n-i-p}, & j=i; q=p+1\hspace{.4pt}; \vspace{1pt} \\
 G( F(d_M^i)^p)^{n-i-p},       & j=i+1, q=p\hspace{.4pt}; \vspace{0.5pt} \\
0,                                                      & \text{otherwise}.
\end{array}\right.$$
On the other hand, we have
$$G^{\hspace{.5pt} C \hspace{-1pt}}(F^{\hspace{.5pt}C}\hspace{-1pt}(\m M^\cdt\hspace{-1pt})\m)^n =
\oplus_{s\in \Z}\hspace{.5pt} G(F^{\hspace{.5pt}C}\hspace{-1.5pt}(M^\cdt\hspace{-.5pt})^s)^{n-s}= \oplus_{(i,s) \in \mathbb{Z} \times \mathbb{Z}} \, G(F\hspace{-,8pt} M^i)^{s-i})^{n-s}.
$$
Writing $G^{\hspace{.5pt}C}(F^{\hspace{.5pt}C}\hspace{-1.4pt}(M^\cdt)\m)^{n+1}=\oplus_{(j,t)\in \mathbb{Z}\times \mathbb{Z}} \hspace{.5pt} G(F(M^j)^{t-j})^{n+1-t},$ we can routinely verify that \vspace{.5pt}
the differential
$d^n_{G^{\hspace{.5pt}C\hspace{-1pt}}(F^{\hspace{.5pt}C}\hspace{-1pt}(M^\cdt))}:
G^{\hspace{.5pt}C\hspace{-1pt}}(\m F^{\hspace{.5pt}C}\hspace{-1pt}(M^\cdt)\m)^n\to
G^{\hspace{.5pt}C\hspace{-1pt}}(F^{\hspace{.5pt}C}\hspace{-1pt}( M^\cdt))^{n+1}$
is given by the matrix $(d^{\hspace{.5pt}n}_{G^{C \hspace{-.5pt}}(F^C\hspace{-.5pt}(M^\cdt))}(j,t;i,s))_{(j,t;i,s)\in \mathbb Z^4},$
where \vspace{-3pt}
$$
d^{\hspace{.5pt}n}_{G^{C\hspace{-.5pt}}(F^C\hspace{-.5pt}(M^\cdt))}(j,t;i,s):
G(F(M^i)^{s-i})^{n-s}\to G(F(M^j)^{t-j})^{n+1-t} \vspace{-2pt}$$
is given by
\vspace{-8pt}
$$
d^{\hspace{.5pt}n}_{G^{C\hspace{-.5pt}}(F^C\hspace{-.5pt}(M^\cdt))}(j,t;i,s)=
\left\{\hspace{-5pt} \begin{array}{ll}
(-1)^{s} d_{G(F(M^i)^{\hspace{.5pt} s-i})}^{\hspace{.4pt}n-s}, & t=s, j=i;  \vspace{2pt} \\
(-1)^i G(d^{\hspace{.5pt} s-i}_{F(M^i)})^{n-s},             & t=s+1, j=i; \vspace{1pt} \\
 G( F(d_M^i)^{s-i})^{n-s}, & t=s+1, j=i+1; \vspace{0.5pt} \\
0,                          & \text{otherwise}.
\end{array}\right.$$
Setting $p=s-i$, we see that $$G^{\hspace{.5pt}C\hspace{-1pt}}(\m F^{\hspace{.5pt}C\hspace{-1pt}}(M^\cdt))^n
= \oplus_{(i,p)\in \mathbb{Z}^2} \, G(F(M^i)^p)^{n-i-p} = (G^{\hspace{.5pt}C} \ncirc F)^{\hspace{.5pt}C\hspace{-1pt}}(M^\cdt)^n.$$ And setting $q=t-j$, we see that
$$d^{\hspace{.5pt}n}_{G^{\hspace{.5pt} C\hspace{-1pt}}(F^C\hspace{-1pt}(\m M^\cdt\hspace{-.8pt})\m)}: G^{\hspace{.5pt}C\hspace{-1pt}}(F^{\hspace{.5pt}C\hspace{-1pt}}(M^\cdt\hspace{-1pt}))^n \to G^{\hspace{.5pt}C\hspace{-1pt}}(F^{\hspace{.5pt}C\hspace{-1pt}}(M^\cdt))^{n+1}=\oplus_{(j,q)\in \mathbb{Z}\times \mathbb{Z}} G(F(M^j)^q)^{n+1-j-q}$$ is given by the matrix $\left(d^{\hspace{.5pt}n}(j,q;i,p)\right)_{(j,q;i,p)\in \mathbb Z^4},$ where
\vspace{-1pt}
$$d^{\hspace{.5pt}n}(j,q;i,p): G(F(M^i)^p)^{n-i-p} \to G(F(M^j)^q)^{n+1-j-q} \vspace{-2pt}$$ is such that $
d^{\hspace{.5pt}n} (j,q;i,p) = d^{\hspace{.5pt}n}_{G^{C \hspace{-1pt}}(F^C\hspace{-1pt}(M^\cdt))}(j, q+j; i, p+i)=d^{\hspace{.5pt}n}_{( G^{C\hspace{-.6pt}}\circ F)^{C\hspace{-.6pt}}(M^\ydt)}(j,q;i,p).$ \vspace{1pt}
This shows that $(G^{\hspace{.3pt}C}\circ F)^C(M^\pdt)=(G^{\hspace{.3pt}C}\circ F^C)(M^\pdt)$. And given a morphism $f^\cdt$ in $C(\cA)$, we verify in a similar manner that $(G^{\hspace{.3pt}C} \hspace{-2pt} \circ\hspace{-1pt} F)^C(f^\ydt)=(G^{\hspace{.3pt}C} \hspace{-2pt} \circ\hspace{-1pt} F^C)(f^\ydt)$.
The proof of the proposition is completed.

%

\smallskip

Now, we show how to extend functorial morphisms.

\begin{Lemma}\label{main-lemma1}

Let $\cA, \cB$ be concrete additive $k$-categories such that $\cB$ has countable direct sums. \vspace{1pt} Given functors $F,  G: \cA\to C(\cB)$, any functorial morphism $\eta: F\to G$ extends to functorial morphisms $\eta^{DC}\!\m: F^{DC}\to G^{DC}$ and $\eta^C\!\m: F^{\hspace{.3pt}C}\to G^{\hspace{.3pt}C}.$

\end{Lemma}

\noindent{\it Proof.} Let $\eta=(\eta_{_{\hspace{-.6pt}M}}^\ydt)_{M\in \mathcal A}: F\to G$ be a functorial morphism between two functors $F,  G: \cA\to C(\cB)$. Fix a complex $M^\cdt\in C(\cA)$. Given $i, j\in \Z$, since $\eta_{_{\hspace{-.6pt}M}}^\ydt$ is natural in $M$, we obtain a commutative diagram \vspace{-0pt}
$$\xymatrixcolsep{16pt}\xymatrixrowsep{8pt}\xymatrix{
&  G(M^i)^{j+1}\\
F(M^i)^{j+1} \ar[ur]^{\eta_{M^i}^{\hspace{.5pt}j+1}}&  \\
& G(M^i)^j \ar[rr]^{ G(d_M^i)^j}\ar[uu]_{(-1)^id^j_{ G(M^i)}}&&  G(M^{i+1})^j.\\
F(M^i)^j \ar[uu]^{(-1)^id^j_{F(M^i)}}\ar[ur]^-{\eta_{M^i} ^{\hspace{.5pt}j} } \ar[rr]^-{F(d_M^i)^j} && F(M^{i+1})^j\ar[ur]_-{\eta_{M^{i+1}}^{\hspace{.5pt}j}}
}\vspace{-1pt}$$

This yields a morphism  $\eta_{_{\hspace{-.6pt}M^\ydt}}^{\ydt}=(\eta_{_{\hspace{-1pt}M^i}}^{\hspace{.5pt}j})_{i, j\in \Z}: F^{DC}(M^\cdt)\to G^{DC}(M^\cdt)$ \vspace{-1pt} in $DC(\cB)$. Applying $\mathbb{T}$, we obtain a morphism $\eta_{_{\hspace{-.6pt}M^\ydt}}=\T(\eta_{_{\hspace{-1pt}M^\ydt}}^\ydt): F^{\hspace{.5pt}C}(M^\ydt)\to G^{\hspace{.5pt}C}(M^\ydt)$ in $C(\cB)$.
Clearly, $\eta_{_{\hspace{-.6pt}M^\ydt}}^\ydt$ \vspace{1.5pt} and $\eta_{_{\hspace{-.6pt}M^\ydt}}$ are natural in $M^\ydt$.
\vspace{-2pt} Thus, we have desired functorial morphisms $\eta^{DC}\!=\! (\eta_{_{\hspace{-.6pt}M^\ydt}}^{\ydt})_{M^\ydt\in C(\cA)}$ and $\eta^{\hspace{.5pt}C}\!=\!(\eta_{_{\hspace{-.6pt}M^\ydt}})_{M^\ydt\in C(\cA)}$.
The proof of the lemma is completed.

\vspace{5pt}

Even if $F: \cA\to C(\cB)$ is exact, the extended functor $F^{\hspace{.5pt}C}$ does not necessarily send all acyclic complexes to acyclic ones. So it only descends to categories derived from some suitable derivable subcategories of $C(\cA)$.

\begin{Theo}\label{Der-functor}

Let $\cA, \cB$ be concrete additive $k$-categories such that $\cB$ has countable direct sums. Consider an exact functor $F: \cA \to C(\cB)$ such that $F^{\hspace{.5pt}C}$ sends a derivable subcategory $\mathscr{A}$ of $C(\cA)$ into a derivable subcategory $\mathscr{B}$ of $C(\cB)$.

\begin{enumerate}[$(1)$]

\item If $F^{DC}\hspace{-2.5pt}$ sends complexes in $\mathscr{A}$ to diagonally bounded-below double complexes, then $F^{\hspace{.4pt}C}$ sends acyclic complexes in $\mathscr{A}$ to acyclic ones.

\vspace{.5pt}

\item If $F^{\hspace{.5pt}C}$ sends acyclic complexes in $\mathscr{A}$ to acyclic ones, then it induces a diagram \vspace{-5pt} $$\xymatrixrowsep{20pt}\xymatrix{\mathscr{A} \ar[r] \ar[d]_{F^{\hspace{.5pt}C}} & \mathcal{K} (\mathscr{A}) \ar[r] \ar[d]^{F^K} & \mathcal{D}(\mathscr{A})\ar[d]^{F^D}\\ \mathscr{B} \ar[r] & \mathcal{K}(\mathscr{B}) \ar[r] &  \mathcal{D}(\mathscr{B}),} \vspace{-3pt}$$ which is commutative with $F^K$ and $F^D$ being triangle-exact.

\end{enumerate}

\end{Theo}

\noindent{\it Proof.} (1) Suppose that $F^{DC}\hspace{-2.5pt}$ sends complexes in $\mathscr{A}$ to diagonally bounded-below double complexes. Let $M^\ydt$ be an acyclic complex in $\mathscr{A}$. Since $F$ is exact, $F(M^\ydt)^\cdt$ has exact rows. By Proposition \ref{Homology-zero}, $\mathbb{T}(F(M^\ydt)^\cdt)$, that is $F^C\hspace{-1pt}(M^\ydt)$, is acyclic.


(2) By Proposition \ref{F-extension}, we have a triangle-exact functor $F^K: \mathcal{K}(\mathscr{A}) \to \mathcal{K}(\mathscr{B})$ making the left square of the diagram stated in Statement (2) commute.
Suppose that $F^{\hspace{.5pt}C}$ sends acyclic complexes in $\mathscr{A}$ to acyclic ones. It is well-known that there exists a triangle-excat functor $\mathcal{F}^D: \mathcal{D}(\mathscr{A}) \to \mathcal{D}(\mathscr{B})$ making right square of the diagram commute. The proof of the theorem is completed.

\vspace{-2pt}

\section{Generalized Koszul dualities}

The main objective of this section is to describe the generalized Koszul dualities, that is a 2-real-parameter family of pairs of mutually quasi-inverse equivalences between categories derived from graded modules over a Koszul algebra and over its Koszul dual, and one of the pairs is the classical Koszul duality of Beilinson, Ginzburg and Soergel; see \cite[(2.12.1)]{BGS} and \cite[Theorem 30]{MOS}.


Throughout this section, $\La=kQ/R$ stands for a quadratic algebra, where $Q$ is a locally finite quiver. We recall and introduce some notation. Given $u\in e_y\La_r e_x$ with $r\in \Z$ and $x,y\in Q_0$ and $M\in \GrLa$, the right multiplication by $u$ yields a graded $\La$-linear morphisms $P[u]: P_y\langle n \rangle \to P_x\langle n+r\rangle $, while the left multiplication by $u$ yields a $k$-linear map $M(u): M_n(x) \to M_{n+r}(y)$, for every $n\in \Z$. For each $x\in Q_0$, we write $P_x^{\hspace{.5pt}!}=\La^! e_x$ and $I_x^!=\mk D(\hspace{3pt}\hat{\hspace{-3.5pt}\La}e_x),$ where $\hspace{3pt}\hat{\hspace{-3.5pt}\La}=(\La^!)^{\rm o}=kQ/(R^!)^{\rm o}$.

Now, we define two Kozsul functors. The first is the so-called {\it right Koszul functor} $\cF: \GrLa\to C(\GrLa^!)$ defined as follows; compare \cite[(3.1)]{BaL2}. Given $M\in \GrLa$, as will be shown below, we have a complex $\cF(M)^\pdt \in C(\GrLa^!)$ such, for all $n\in \Z$, that
$\cF(M)^n \hspace{-2.5pt}=\hspace{-1.5pt} \oplus_{x\in Q_{0}}P_x^{\hspace{.5pt}!}\nla n\nra\otimes M_{n}(x)$ and \vspace{-2pt} $$d^{n}_{\cF(M)}=({\textstyle\sum}_{\alpha \in Q_1(x, y)} P[\bar\alpha^!]\otimes M(\bar\alpha)\hspace{-1pt})_{(y,x)\in Q_0\times Q_0}: \hspace{2pt} \cF(M)^n\to \cF(M)^{n+1} \vspace{-2pt}
$$ with $P[\bar\alpha^!]\otimes M(\bar\alpha)\hspace{-1pt}: P_x^{\hspace{.5pt}!} \nla n\nra\otimes M_{n}(x)\to
P_y^{\hspace{.5pt}!} \langle n\!+\!1\rangle \otimes M_{n+1}(y),$ for every $\alpha\in Q_1(x,y).$ \vspace{1pt}
And given a morphism $f: M \to N$ in $\GrLa$, we have a complex morphism $\cF(f)^\cdt: \cF(M)^\cdt\to \cF(N)^\cdt$ such, for all $n\in \Z$, that \vspace{-3pt} $$\cF(f)^{n}=\oplus_{x\in Q_{0}}({\rm id} \otimes f_{n,x}): \oplus_{x\in Q_{0}} P_x^{\hspace{.5pt}!} \tla n\tra\otimes M_n(x) \to \oplus_{x\in Q_{0}} P_x^{\hspace{.5pt}!} \nla n\nra\otimes N_n(x),\vspace{-2pt}$$ where $f_{n, x}: M_n(x)\to N_n(x)$ is the $k$-linear map obtained by restricting $f$.

The second is the so-called {\it left Koszul functor} $\cG: \GrLa\to C(\GrLa^!)$ defined as follows. Given $M\in \GrLa$, we will have a complex $\cG(M)^\pdt \in C(\GrLa^!)$ such, for all $n\in \Z$, that $\cG(M)^{n}=\oplus_{x\in Q_{0}}I_{x}^{\hspace{.5pt} !}\langle n \rangle \otimes M_{n}(x)$ and \vspace{-2pt}
$$d^{n}_{\cG(M)}=(\textstyle\sum_{\alpha \in Q_1(x,y)}I[\bar\alpha^!]\otimes M(\bar \alpha))_{(y,x)\in Q_{0}\times Q_{0}}:\hspace{2pt} \cG(M)^{n}\to \cG(M)^{n+1} \vspace{-2pt}$$ with
$I[\bar\alpha^!]\otimes M(\bar \alpha): I_x^{\hspace{.5pt}!} \langle n \rangle \otimes M_{n}(x)\to I_{y}^{\hspace{.5pt}!} \langle n+1\rangle \otimes M_{n+1}(y)$, for every $\alpha \in Q_1(x,y)$. And given a morphism $f: M\to N$ in $\GrLa$, we will have a complex morphism $\cG(f)^\cdt: \cG(M)^\cdt\to \cG(N)^\cdt$ such, for all $n\in \Z$, that \vspace{-3pt}
$$\cG(f)^{n}=\oplus_{x\in Q_{0}} ({\rm id}\otimes f_{n,x}): \oplus_{x\in Q_{0}}I_x^{\hspace{.5pt}!} \langle n \rangle \otimes M_n(x)\to \oplus_{x\in Q_0} I_x^{\hspace{.5pt}!} \langle n \rangle \otimes N_n(x).$$

\vspace{0pt}

\begin{Prop}\label{F-property}

Let $\La=kQ/R$ be a quadratic algebra with $Q$ a locally finite quiver. The above construction yields exact functors
$\cF: \GrLa \to C(\GrLa^!)$ and $\cG:\GrLa \to C(\GrLa^!)$.

\end{Prop}

\noindent{\it Proof.} Fix $M\in \GrLa$. For $n\in \Z$, write $\cF(M)^n = \oplus_{x\in Q_{0}} P_x^{\hspace{.4pt}!}\tla n \nra\otimes M_{n}(x)$ and $\cF(M)^{n+2} = \oplus_{z\in Q_{0}}P_z^{\hspace{.4pt}!}\nla n+2\tra\otimes M_{n+2}(z)$ with $d^{n+1}_{\cF(M)} \hspace{-2pt} \circ \m d^{n}_{\cF(M)} \!\m = \! (d^{n}_{z,x})_{(z,x)\in Q_0\times Q_0},$ where \vspace{.2pt}
$d^{n}_{z,x}\!:\!  P_x^! \nla n \nra \m\otimes\! M_{n}(x)\to P_z^!\nla n\!+\!2\tra \m\otimes\! M_{n+2}(z)$. Fix $(z, x)\in Q_0\times Q_0$. Write $Q_2(x, z)=\{\alpha_1\beta_1, \ldots, \alpha_s\beta_s\}$, where $\alpha_i, \beta_i\in Q_1$. \vspace{.5pt} In view of the definition of $\cF$, \vspace{.5pt} we see that $d^{n+1}_{z,x}={\textstyle\sum}_{i=1}^s P[\hspace{.5pt}\bar\beta_i^{\hspace{.5pt}!} \hspace{.4pt} \bar\alpha_i^{\hspace{.5pt}!}\hspace{.5pt}] \otimes M(\bar\alpha_i\hspace{.4pt} \bar\beta_i).$

\vspace{.5pt}

As seen in the proof of Lemma \ref{perp}, we have $k$-bases $\{\rho_1, \ldots, \rho_r, \rho_{r+1}, \ldots, \rho_s\}$ and
$\{\eta_1, \ldots, \eta_r, \eta_{r+1}, \ldots, \eta_s\}$ of $kQ_2(x,z)$ such that $\{\rho_1, \ldots, \rho_r\}$ is a basis of $R_2(x,z)$; $\{\eta_{r+1}^{\rm o}, \ldots, \eta_s^{\rm o}\}$ is a basis of $R^{\hspace{.5pt}!}_2(z,x)$ and $\{\eta_1^*, \ldots, \eta_s^*\}$ is the dual basis of $\{\rho_1, \ldots, \rho_s\}$. Then, $\bar\rho_i=0$ for $1\le i\le r$, and $\bar\eta^!_j=0$ and for $r<j\le s$. By Corollary \ref{Cor 1.2}, we have a $k$-linear isomorphism $\sigma: D(kQ_2(x,z))\otimes  kQ_2(x,z)\to \End_k(kQ_2(x,z))$. We easily see that $\textstyle \sigma(\sum_{i=1}^s (\alpha_i\beta_i)^*\otimes \alpha_i\beta_i)={\rm id}=\sigma({\sum}_{i=1}^s \eta_i^* \otimes \rho_i).$ Thus, \vspace{-2pt} $$\textstyle \sum_{i=1}^s (\alpha_i\beta_i)^*\otimes \alpha_i\beta_i= {\sum}_{i=1}^s \eta_i^* \otimes \rho_i. \vspace{-3pt}$$

Applying the $k$-linear isomorphism $D(kQ_2(x,z))\to kQ^{\rm o}_2(z,x): \xi^*\mapsto \xi^{\rm o}$, we obtain $\sum_{i=1}^s (\alpha_i\beta_i)^{\rm o}\otimes \alpha_i\beta_i= {\sum}_{i=1}^s \eta_i^{\rm o} \otimes \rho_i$. And applying the tensor product of the canonical projections $kQ_2(x,z)\to e_z\La_2 e_x$ and $kQ_2^{\rm o}(z,x)\to e_x\La_2^!e_z$, we obtain \vspace{-2pt} $$\textstyle\sum_{i=1}^s \bar\beta_i^{\hspace{.5pt}!} \hspace{.4pt} \bar\alpha_i^{\hspace{.5pt}!} \otimes \bar{\alpha}_i\bar{\beta}_i= \sum_{i=1}^s \bar\eta_i^{\hspace{.5pt}!} \otimes \bar{\rho}_i.\vspace{-4pt}$$

Finally, we have a $k$-isomorphism $e_z\La_2^! e_x \to \Hom_{\it\Lambda}(P_z^{\hspace{.4pt}!}\mla n\!-\!2\tra, P_x^{\hspace{.4pt}!}\mla n\tra): u \mapsto P[u];$ see \cite[(2.3.3)]{LLi}, and a $k$-linear map $ e_x \La_2 e_z \to \Hom_k(M_{n-2}(z), M_{n}(x)): u\mapsto M(u)$. This yields a $k$-linear map \vspace{-2pt} $$\psi: e_z\La_2^! e_x \otimes e_x \La_2 e_z \to \Hom_{\it\Lambda}(P_z^!\mla n\!-\!2\tra, P_x^!\mla n\tra) \otimes \Hom_k(M_{n-2}(z), M_{n}(x)),\vspace{-2pt} $$ sending $\bar\eta^!\otimes \bar \rho $ to $P[\bar\eta^!]\otimes M(\bar\rho)$. Applying it to the above equation, we obtain \vspace{-1pt} $${\textstyle\sum}_{i=1}^s P[\hspace{.5pt}\bar\beta_i^{\hspace{.5pt}!} \hspace{.4pt} \bar\alpha_i^{\hspace{.5pt}!}\hspace{.5pt}] \otimes M(\bar{\alpha}_i\bar{\beta}_i)=
{\textstyle\sum}_{i=1}^s P[\bar\eta_i^!] \otimes M(\bar{\rho}_i)=0.\vspace{-1pt}$$

Therefore, $d^{n-1}_{\cF(M)}\circ d^{n-2}_{\cF(M)}=0$. That is, $\cF(M)^\cdt\in C(\GrLa^!)$. Given a morphism $f: M\to N$, it is easy to verify that $\cF(f)^{n} \circ d_{\cF(M)}^{n-1}=d_{\cF(N)}^{n-1}\circ \cF(f)^{n-1}$, for $n\in \Z$. Thus, we have a morphism $\cF(f)^\cdt: \cF(M)^\cdt\to \cF(N)^\cdt$. So $\cF$ is a functor, which is exact the tensor product is over a field. The proof of the proposition is completed.

\vspace{3pt}

Let $M^\ydt$ be a complex and $f^\ydt: M^\ydt\to N^\ydt$ a morphism in $C({\rm GMod}\hspace{.5pt} \La)$. Given $i\in \Z$,
restricting $M^\ydt$ and $f^\ydt$ to the degree $i$, we obtain a complex $M^\ydt_i$ and a morphism $f^\ydt_i: M^\ydt_i\to N^\ydt_i$ in $C({\rm Mod}\hspace{.5pt}k)$. Note that
 ${\rm H}^n(M^\cdt)\cong \oplus_{i\in \Z} {\rm H}^n(M_i^\cdt)$, for $n\in \Z$. And $f^\ydt$ is a quasi-isomorphism if and only if so is $f^\ydt_i$ for every $i\in \Z.$ A similar consideration is given to a double complex $M^{\cdt\hspace{.4pt}\cdt}$ and a morphism $f^{\hspace{.4pt}\cdt\hspace{.4pt}\cdt}$ in
$DC(\GrLa)$ so that $\T(M^{\cdt\hspace{.4pt}\cdt})_i=\T(M^{\cdt\hspace{.4pt}\cdt}_i)$ and $\T(f^{\hspace{.4pt}\cdt\hspace{.4pt}\cdt})_i=\T(f^{\hspace{.4pt}\cdt\hspace{.4pt}\cdt}_i)$, for all $i\in \mathbb Z$.

\vspace{1pt}

Given $s\in \Z$, the {\it grading $s$-shift} $M^\ydt\sla s\sra$ of $M^\ydt$ is defined by $(M^\ydt\sla s\sra)^n=M^n\sla s\sra$ and $d^n_{\hspace{-1pt} M^\cdt \m \sla s\sra}=d^n_{\hspace{-.8pt} M}\sla s\sra$ for all $n\in \Z$. In view of the above consideration, we see that ${\rm H}^n(M^\cdt\sla s\sra) \!= \! {\rm H}^n(M^\cdt)$ for all $n\in \Z$. And it is clear that $\mathfrak{t}(M^\cdt\sla s\sra)=\mathfrak{t}(M^\cdt)\sla s\sra$, where $\mathfrak t$ is the twist functor. The following statement follows from a routine verification.

\begin{Lemma}\label{Ksl-Fun-Homology}

Let $\La=kQ/R$ be a quadratic algebra with $Q$ a locally finite quiver. Given any $M\in \GrLa$ and $s\in \Z$, we have $\cF(M)^\cdt [s]=\mathfrak{t}^s( \cF(M\sla s \sra)^\cdt\sla s \sra)$ and $\cG(M)^\cdt[s]=\mathfrak{t}^s(\cG(M\sla s \sra)^\cdt\sla s \sra)$.

\end{Lemma}


\vspace{2pt}

In view of Proposition \ref{F-extension}, we can extend the two Koszul functors as follows.

\begin{Lemma}\label{cplx-kosz-fun}

Let $\La=kQ/R$ be a quadratic algebra with $Q$ a locally finite quiver.

\begin{enumerate}[$(1)$]

\vspace{-1pt}

\item \hspace{-5pt} The right Koszul functor $\cF\hspace{-1pt}$ extends to $\cF^{\hspace{.5pt}C}\!: C\hspace{-.4pt}(\hspace{-.6pt}\GrLa) \to C\hspace{-.4pt}(\hspace{-.8pt}\GrLa^!),$ called the {\em right complex Koszul functor} such, for any $M^\cdt\in C\hspace{-.4pt}(\hspace{-.8pt}\GrLa)$, that \vspace{-1.5pt}
$$\cF\hspace{.6pt}^C\hspace{-1pt}(\m M^\cdt\m)^n=\oplus_{(i, x)\in \Z\times Q_0} P^{\hspace{.4pt}!}_x\tla n\!-\!i \tra\otimes M^i_{n-i}(x); \; n\in \Z.
\vspace{-2.5pt}$$

\item \hspace{-5pt} The left Koszul functor $\cG\hspace{-1pt}$ extends $\cG^{\hspace{.6pt}C}\!: C(\m\GrLa) \m\to\m C(\m\GrLa^!),$ called the {\em left complex Koszul functor} such, for any $M^\cdt\in C\hspace{-.4pt}(\hspace{-.8pt}\GrLa)$, that
\vspace{-2pt} $$\cG\hspace{.6pt}^C\hspace{-1pt}(M^\cdt)^n=\oplus_{(i, x)\in \Z\times Q_0} I^{\hspace{.4pt}!}_x\tla n\!-\!i\tra\otimes M^i_{n-i}(x); \; n\in \Z.
\vspace{-0pt}$$

\end{enumerate}\end{Lemma}


\vspace{0pt}

We compose Koszul functors and complex Koszul functors as follows.


\begin{Lemma}\label{Comp_Kfunctors}

Let $\La=kQ/R$ be a quadratic algebra with $Q$ a locally finite quiver.

\begin{enumerate}[$(1)$]

\item The functor $\cF^{\hspace{.6pt}C}\!\circ \cG: \GrLa \to C(\GrLa)$ is such, for $M\in \GrLa$, that \vspace{-3.5pt}
$$(\cF^{\hspace{.6pt}C}\!\m\circ\m \cG)(M)^n = \oplus_{i\in \mathbb{Z}; \hspace{.7pt} a, x\in Q_0} P_a\tla n\!-\!i\tra \otimes (I_x^!)_n(a) \otimes M_i(x); \;  n\in \Z.$$

\vspace{-1.5pt}

\item The functor $\cG^{\hspace{.6pt}C}\circ \cF: \GrLa \to C(\GrLa)$ is such, for $M\in \GrLa$, that \vspace{-3.5pt}
$$(\cG^{\hspace{.6pt}C} \circ \cF)(M)^n=\oplus_{i\in \Z; \hspace{.7pt} x, a\in Q_0} I_a\tla n\!-\!i\tra\otimes (P_x^!)_n(a) \otimes M_i(x); \; n\in \Z.\vspace{-3pt}$$

\end{enumerate}

\end{Lemma}

\noindent{\it Proof.} We only verify Statement (1). Consider a module $M\in \GrLa$. By definition, $(\cF^{\hspace{.4pt}C}\!\!\circ\m \cG)(M)^\cdt = \cF^{\hspace{.4pt}C}(\cG(M)^\cdt) = \T(\cF(\cG(M)^\cdt\hspace{-1pt})^\cdt\hspace{-.8pt}).$ Thus, for any integer $n$, we have $(\cF^{\hspace{.4pt}C}\!\m \circ\m \cG)(M)^n
\m =\m \oplus_{i\in \Z}\,\cF(\cG(M)^i)^{n-i}$. Since $\cG(M)^i=\oplus_{x\in Q_0} I_x^!\langle i \rangle \otimes M_i(x)$, we see from the definition of $\cF$ that \vspace{-3pt}
$$\cF(\cG(M)^i)^{n-i}  \m = \m \oplus_{a\in Q_0} P_a \m \langle n-i \rangle \otimes \hspace{1.5pt} \cG(M)^i_{n-i}(a) \m = \m \oplus_{a, x\in Q_0} \m P_a \m \langle n-i \rangle \otimes (I_x^!)_n(a) \otimes M_i(x). \vspace{-2pt}$$
The proof of the lemma is completed.

\vspace{3pt}

The above two composite functors are extended as follows.

\begin{Lemma}\label{Extend_comp_Kfunctor}

Let $\La=kQ/R$ be a quadratic algebra with $Q$ a locally finite quiver. Consider a complex $M^\cdt \!\in\! C(\GrLa)$ and an integer $n$.

\begin{enumerate}[$(1)$]

\vspace{1pt}

\item The functor $(\cF^{\hspace{.6pt}C} \ncirc \cG)^C : C(\GrLa) \!\to\! C(\GrLa)$ is such, for $n\in \Z$, that \vspace{-2pt}
$$(\cF^{\hspace{.6pt}C} \ncirc \cG)^{\hspace{.6pt}C}\hspace{-1.5pt}( M^\cdt)^n =\oplus_{i,j\in \mathbb{Z}; \hspace{.5pt} a, x \in Q_0} P_a\tla n\!-\!i\!-\!j\tra \otimes (I_x^{\hspace{.4pt}!})_{n-i}(a) \otimes M^i_j(x).$$

\vspace{-2pt}

\item The functor $(\cG^{\hspace{.6pt}C} \circ \cF)^{\hspace{.6pt}C} : C(\GrLa) \!\to\! C(\GrLa)$ is such, for $n\in \Z$,
that \vspace{-2pt}
$$(\cG^{\hspace{.6pt}C} \circ \cF)^{C}\hspace{-1pt}(M^\ydt)^n = \oplus_{i,j\in \mathbb{Z}; \hspace{.5pt} a, x \in Q_0}\, I_x\mla n\!-\!i\!-\!j \tra \otimes (P^{\hspace{.4pt}!}_x)_{n-i}(a) \otimes M^i_j(a). \vspace{-3.5pt}$$

\end{enumerate}

\end{Lemma}

\noindent{\it Proof.} We shall only verify Statement (1). \vspace{1pt} Let $M^\cdt\in C(\GrLa)$. By definition, $(\cF^{\hspace{.5pt}C} \ncirc \cG)^{C} \m (M^\ydt)=\T((\cF^{\hspace{.5pt}C} \circ \cG)(M^\ydt)^\ydt).$ So
$(\cF^{\hspace{.5pt}C}  \circ \cG)^C\hspace{-1.5pt}(\m M^\ydt)^n = \oplus_{i\in \Z}\, (\cF^{\hspace{.5pt}C} \ncirc \cG)(M^i)^{n-i}$, for any $n\in \Z$. In view of Lemma \ref{Comp_Kfunctors}, we see that \vspace{-2pt}
$$(\cF^{\hspace{.6pt}C} \ncirc \cG)(M^i)^{n-i} = \oplus_{j\in \mathbb{Z}; \hspace{.7pt} a, x\in Q_0} P_a\tla n\!-\!i\!-\!j\tra \otimes (I_x^!)_{n-i}(a) \otimes M^i_j(x). \vspace{-2pt}$$
This implies Statement (1). The proof of the lemma is completed.

\vspace{3pt}

Next, we show that the complex Koszul functors descend to categories derived from some subcategories of $C(\GrLa)$. To introduce these subcategories, we will view a complex $M^\cdt$ of graded modules $M^i\!=\!\oplus_{j\in \Z} M^i_j$ as a bigraded $k$-vector spaces $M^i_j$ with $i, j\in \Z$.

\begin{Defn}

Let $\La=kQ/R$ be a quadratic algebra with $Q$ a locally finite quiver.
Given $p, q\in \mathbb{R}$ with $p\ge 1$ and $q\ge 0$, we denote

\begin{enumerate}[$(1)$]

\vspace{-2pt}

\item by $C^{\,\downarrow}_{p,q}(\GrLa)$ \vspace{0.5pt} the full subcategory of $C(\GrLa)$ of complexes $M^\cdt$ such that $M^i_j=0$ for $i+pj \gg 0$ or $i-qj \ll 0\,;$ in other words, $M^\cdt$ concentrates in a lower triangle formed by two lines of slopes $-\frac{1}{p}$ and $\frac{1}{q},$ respectively;

\vspace{1pt}

\item by $C^{\,\uparrow}_{p,q}(\GrLa)$ \vspace{0.5pt} the full subcategory of $C(\GrLa)$ of complexes $M^\cdt$ such that $M^i_j=0$ for $i+pj\ll 0$ or $i-qj\gg 0;$ in other words, $M^\cdt$ concentrates in a upper triangle formed by two lines of slopes $-\frac{1}{p}$ and $\frac{1}{q}$, respectively;

\vspace{1pt}

\item by $C^{\,\downarrow}_{p,q}({\rm gmod}\hspace{.4pt}\La^!)$ \vspace{0.5pt} and $C^{\,\uparrow}_{p, q}({\rm gmod}\hspace{.4pt}\La)$ the full subcategories of $C^{\,\downarrow}_{p,q}(\GrLa)$ and $C^{\,\uparrow}_{p,q}(\GrLa)$ respectively of piecewise finite dimensional modules.

\end{enumerate}

\end{Defn}

\noindent{\sc Remark.} (1) Taking $p=1$ and $q=0$, we recover the categories $C^\downarrow(\La)$ and $C^\uparrow(\La)$ defined in \cite[(2.12)]{BGS}; see also, \cite[(2.4)]{MOS}.

\vspace{1pt}

\noindent (2) The categories $C^{\,\downarrow}_{p,q}(\GrLa)$ \vspace{-.5pt} are pairwise distinct derivable subcategories of $C({\rm GMod}^-\hspace{-3pt}\La)$ containing $C^b({\rm GMod}^-\hspace{-3pt}\La),$ and the $C^{\,\uparrow}_{p,q}({\rm GMod}\hspace{.5pt}\La)$ are pairwise distinct derivable subcategories of $C({\rm GMod}^+\hspace{-3pt}\La)$ containing $C^b({\rm GMod}^+\hspace{-3pt}\La).$


\vspace{3pt}

Let $\mk A \!=\!\GrLa$ or ${\rm gmod}\La$. In the sequel, we shall denote by $K^{\,\downarrow}_{p,q}(\mk A)$ and $K^{\,\uparrow}_{p,q}(\mk A)$ \vspace{0.5pt} the quotient categories modulo null-homotopic morphisms of $C^{\,\downarrow}_{p,q}(\mk A)$ and $C^{\,\uparrow}_{p,q}(\mk A)$ respectively; and by $D^{\,\downarrow}_{p,q}(\mk A)$ \vspace{0.5pt} and $D^{\,\uparrow}_{p,q}(\mk A)$ the localizations at quasi-isomor\-phisms of $K^{\,\downarrow}_{p,q}(\mk A)$ and $K^{\,\uparrow}_{p,q}(\mk A)$ respectively.

\begin{Theo}\label{F-diag}

Let $\La=kQ/R$ be a quadratic algebra with $Q$ a locally finite quiver.
Consider $p, q\in \mathbb R$ with $p\ge 1$ and $q\ge 0$.

\begin{enumerate}[$(1)$]

\vspace{-0pt}

\item The right complex Koszul functor $\cF^C$ induces a commutative diagram 
\vspace{-3pt}
$$\xymatrixrowsep{18pt}\xymatrixcolsep{18pt}\xymatrix{
\hspace{18pt}C^{\,\downarrow}_{p,q}\hspace{.4pt}(\GrLa) \ar[d]_{\cF^C_{p,q}} \ar[r]
& \hspace{2pt} K^{\,\downarrow}_{p,q}\hspace{.4pt}(\GrLa\hspace{.3pt})\ar[d]^{\cF^K_{p,q}} \hspace{2pt} \ar[r]
& \hspace{2pt} D^{\,\downarrow}_{p,q}\hspace{.4pt}(\GrLa) \ar[d]^{\cF^D_{p,q}} \hspace{20pt} \\
C^{\,\uparrow}_{q+1, p-1}(\GrLa^!)\ar[r]
& K^{\,\uparrow}_{q+1, p-1}(\GrLa^!) \ar[r]
& D^{\,\uparrow}_{q+1, p-1}(\GrLa^!).}\vspace{-3pt}$$

\item The left complex Koszul functor $\cG^{\hspace{.3pt}C}$ induces a commutative diagram
\vspace{-3pt}
$$\xymatrixrowsep{18pt}\xymatrixcolsep{18pt}\xymatrix{
\hspace{18pt} C^{\,\uparrow}_{p,q}\hspace{.4pt}(\GrLa)\ar[r]
\ar[d]_{\cG^{\hspace{.3pt}C}_{p,q}} & \hspace{2pt} K^{\,\uparrow}_{p,q}\hspace{.4pt}(\GrLa) \hspace{2pt} \ar[r]
\ar[d]^{\cG^K_{p,q}} & D^{\,\uparrow}_{p,q}\hspace{.4pt}(\GrLa)\ar[d]^{\cG^D_{p,q}} \hspace{20pt}  \\
C^{\,\downarrow}_{q+1, p-1}(\GrLa^!) \ar[r]
& K^{\,\downarrow}_{q+1, p-1}(\GrLa^!) \ar[r]
& D^{\,\downarrow}_{q+1, p-1}(\GrLa^!).}$$\vspace{-8pt}

\item In the above two statements,  $\GrLa$ and $\GrLa^!$ can be replaced simul\-taneously by ${\rm gmod}\La$ and ${\rm gmod}\hspace{.4pt}\La^!\hspace{-2pt},$ respectively.

\vspace{1pt}

\end{enumerate}

\end{Theo}

\noindent{\it Proof.} Note that $K^b({\rm GProj}\La)$ and $K^b({\rm GInj}\La)$ are full triangulated subcategories of $D^b(\GrLa)$, and $K^b({\rm gproj}\La)$ and $K^b({\rm ginj}\La)$ are full triangulated subcategories of $D^b({\rm gmod}\La)$; see \cite[(10.4.7)]{Wei}. Consider $M^\cdt \!\m\in\! C^{\,\downarrow}_{p,q}(\GrLa).$ There exist $s, t\!\in\! \Z$ such that $M^i_j=0$ when $i+pj \! > \! s$ or $i-qj \!< \! t.$ Fix $n, m\in \Z$. \vspace{.5pt} By Lemma \ref{cplx-kosz-fun}, \vspace{-2.5pt}
$$\cF^{\hspace{.4pt}C\hspace{-1pt}}(M^\ydt)^{n}_m=\oplus_{i\le n+m; \hspace{.4pt} x\in Q_0}(P_x^!)_{n+m-i}\otimes M^{i}_{n-i}(x).\vspace{-2.5pt}$$

Fix some $i\le n+m$. If $n+(q+1)m<t$, then $i-q(n-i)<t;$ and if $n-(p-1)m>s$, then $i+p(n-i)\!>\!s$. Thus, $\cF^{\hspace{.4pt}C\hspace{-1pt}}(M^\cdt)^n_m\!=\!0$ in case $n\!+\!(q+1)m\m<\m t$ or $n-\!(p-\!1)m\m>\m s$. That is, $\cF^{\hspace{.4pt}C\hspace{-.6pt}}(M^\cdt)^n\in
C^{\,\uparrow}_{\!q+1\m, p-1}\hspace{-1pt}(\m\GrLa^!\m).$ This yields a functor $\cF^{\hspace{.4pt}C\hspace{-1pt}}_{\!p,q}\hspace{-1pt}:\hspace{-1pt} C^{\,\downarrow}_{p,q}\m(\m\GrLa\m) \! \to\! C^{\,\uparrow}_{\!q+1\m, p-1}\hspace{-1pt}(\m\GrLa^!\m).$ \vspace{.5pt} Furthermore, the $n$-diagonal of $\cF(M^\cdt\hspace{-1pt})^\cdt$ consists of $\cF(M^i)^{n-i}=\oplus_{x\in Q_0}\, P_x^!\nla n\!-\!i\nra\otimes M^i_{n-i}(x)$ with $i\in \Z$. Since $M^i_{n-i}=0$ for $i< (nq +t)(1+q)^{-1}$, we see that $\cF(M^\cdt\m)^\cdt$ is diagonally bounded-below. By Theorem \ref{Der-functor}, $\cF^{\hspace{.4pt}C\hspace{-1pt}}_{\!p,q}$ induces a commutative diagram as stated in Statement (1).

\vspace{.5pt}

Similarly, $\cG^{\hspace{.3pt}C}$ restricts to a functor $\cG^{\hspace{.4pt}C\hspace{-1pt}}_{p,q}\!:\! C^{\,\uparrow}_{p,q}(\GrLa) \!\to\! C^{\,\downarrow}_{q+1, p-1}(\GrLa^!)$. Let $N^\cdt\in C(\GrLa^!)$ be acyclic. Fix $m\in \Z$. Then $\cG^{\hspace{.4pt}C\hspace{-1.3pt}}(N^\cdt)_m=\T(\cG(N^\ydt)^\ydt_m)$. Since $\cG$ is exact, $\cG(N^\ydt)^\ydt$ has acyclic rows, and so does $\cG(N^\ydt)^\ydt_m$. Given $n\in \Z$, in view of Lemma \ref{cplx-kosz-fun}, the $n$-diagonal of $\cG(N^\cdt)^\ydt_m$ consists of $$\cG(N^i)^{n-i}_m\!=\m\oplus_{x\in Q_0} (I^!_x)_{n+m-i}\otimes N^i_{n-i}(x); i\in \Z.$$  If $i<n+m$, \vspace{.5pt} then $\cG(N^i)^{n-i}_m=0$. That is, $\cG(N^\cdt)^\cdt\hspace{-1pt}_m$ is diagonally bounded-below. By Proposition \ref{Homology-zero}, $\cG^{\hspace{.4pt}C\hspace{-1.3pt}}(N^\cdt)_m$ is acyclic. Hence, $\cG^{\hspace{.4pt}C\hspace{-1.3pt}}(N^\cdt\m)$ is acyclic. Thus, $\cG^{\hspace{.4pt}C}$ induces a commutative diagram as stated in Statement (2); see (\ref{Der-functor}).

Finally, assume that $M^\ydt\in C^{\,\downarrow}_{p,q}({\rm gmod}\hspace{.4pt}\La)$. \vspace{.5pt}
Fix $n,m\in \Z$ and $y\in Q_0$. Then, $\textstyle \cF^C\m(M^\cdt)^n_m(y)=\oplus_{i\in \mathbb{Z};\, x\in Q_0} e_y\La_{n+m-i}^!e_x\otimes M^i_{n-i}(x).$ Clearly, $\cF^C\!(M^\cdt)^n_m(y)\ne 0$ only if $(qn+t)(q+1)^{-1}\le i \le n+m$. Hence, $\textstyle \cF^C\!(M^\cdt)^n_m(y)$ is a finite direct sum of finite dimensional $k$-spaces. So, $\cF^C\!(M^\cdt)^\cdt\in C({\rm gmod}\hspace{.4pt}\La^!)$. As seen above, $\cF^C\!(M^\cdt)$ lies in $C^{\,\uparrow}_{q+1, p-1}({\rm gmod}\hspace{.4pt}\La^!)$. This yields a functor $\cF^{\hspace{.4pt}C\hspace{-1pt}}_{p,q}: C^{\,\downarrow}_{p,q}({\rm gmod}\La)\to C^{\,\uparrow}_{q+1, p-1}({\rm gmod}\La^!).$ Similarly, we obtain a functor $\cG^{\hspace{.4pt}C\hspace{-1pt}}_{p,q}: C^{\,\uparrow}_{p,q}({\rm gmod}\La) \!\to\! C^{\,\downarrow}_{q+1, p-1}({\rm gmod}\La^!)$. \vspace{.5pt} The rest of the proof of Statement (3) are similar to those of the first two statements. The proof of the theorem is completed.

\vspace{3pt}

\noindent{\sc Remark.} (1) In case $p=0$ and $q=1$, Theorem \ref{F-diag} has been established for positively graded quadratic categories; see \cite[Proposition 20]{MOS}.

\vspace{1pt}

\noindent (2) In the sequel, we shall call $\cF^D_{pq}$ and $\cG^D_{pq}$ the {\it right } and the {\it left derived Koszul functors}, respectively.

\vspace{2pt}

We shall show that the complex Koszul functors always descend to bounded derived category of finitely piece-supported graded modules.

\begin{Theo}\label{F-diag-fd}

Let $\La=kQ/R$ be a qudratic algebra with $Q$ a locally finite quiver.

\vspace{-.0pt}

\begin{enumerate}[$(1)$]

\item The Koszul functor $\cF\!:\m \GrLa \!\to\m C(\GrLa^!)$ induces a commutative diagram\vspace{-2pt}
$$\xymatrixrowsep{18pt}\xymatrixcolsep{20pt}\xymatrix{ C^{\hspace{.5pt}b}\m({\rm GMod}^{b\hspace{-3pt}}\La) \ar[d]_{\cF^C} \ar[r] & K^{\hspace{.5pt}b}\m({\rm GMod}^{b\hspace{-3pt}}\La)\ar[d]_{\cF^K} \ar[r]& D^{\hspace{.5pt}b}\m({\rm GMod}^{b\hspace{-3pt}}\La) \ar[d]_{\cF^D} \\ C^{\hspace{.5pt}b}\m({\rm GProj}\La^!)\ar[r] & K^{\hspace{.5pt}b}\m({\rm GProj}\La^!) \ar[r] & D^{\hspace{.5pt}b}\m({\rm GMod}\La^!).} \vspace{-3pt}$$

\item The Koszul functor $\cG\!:\m \GrLa \!\to\m C(\GrLa^!)$ induces a commutative diagram \vspace{-6pt} $$\xymatrixrowsep{18pt}\xymatrixcolsep{20pt}\xymatrix{
C^{\hspace{.5pt}b}\m({\rm GMod}^{b\hspace{-3pt}}\La) \ar[d]_{\cG^{\hspace{.3pt}C}} \ar[r] & K^{\hspace{.5pt}b}\m({\rm GMod}^{b\hspace{-3pt}}\La)\ar[d]_{\cG^K} \ar[r] & D^{\hspace{.5pt}b}\m({\rm GMod}^{b\hspace{-3pt}}\La) \ar[d]_{\cG^D} \\ C^{\hspace{.5pt}b}\m({\rm GInj}\La^!)\ar[r] & K^{\hspace{.5pt}b}\m({\rm GInj}\La^!) \ar[r] & D^{\hspace{.5pt}b}\m({\rm GMod}\La^!).} \vspace{-1pt}$$

\item In the above two statements, ${\rm GMod}^{b\hspace{-3pt}}\La$, ${\rm GProj}\La^!\!$, ${\rm GInj}\La^!\!$ and ${\rm GMod}\La^!\hspace{-2pt}$ can be replaced simultaneously by ${\rm gmod}^{b\hspace{-3pt}}\La$, ${\rm gproj}\La^!\!$, ${\rm ginj}\La^!\hspace{-2pt}$ and ${\rm gmod}\La^!\hspace{-2pt},$ respectively.

\end{enumerate}

\end{Theo}

\noindent{\it Proof.} Let $M^\cdt \!\in\m C^b\m({\rm GMod}^{b\hspace{-3.5pt}}\La)$, say $M^i_j\ne 0$ only if $-s\le i \le s$ and $-t\le j\le t$, for some $s, t>0$. By Lemma \ref{cplx-kosz-fun}, $\cF^{\hspace{.5pt}C}\hspace{-1pt}(M^\cdt)^n=\oplus_{(i,x)\in \Z\times Q_0} P_x^!\tla n\!-\!i\nra\otimes M^i_{n-i}(x)$, for all $n\in \Z$. Since $M^i_{n-i}\ne 0$ only if $-s-t\le n\le s+t$, the complex $\cF^{\hspace{.5pt}C}\m(M^\cdt)$ is bounded. And since the $M^i$ are finitely piece-supported, $M^i_{n-i}(x)\ne 0$ only for finitely many $(i, x)$ with $-s\le i \le s$ and $x\in Q_0$. Thus, $\cF^{\hspace{.5pt}C}\!(M^\cdt)\in C^b({\rm GProj}\La^!).$ So, we have a functor $\cF^C\!\!:\m C^b\m({\rm GMod}^{b\hspace{-3pt}}\La) \m\to\m C^b({\rm GProj}\La^!)$ which, as seen in the proof of Theorem \ref{F-diag}, induces a commutative diagram as stated in Statement (1). Similarly, Statement (2) holds.

Now, suppose that $M^\cdt \!\m\in\! C^b\m({\rm gmod}^{b\hspace{-3pt}}\La)$. Since the $M^i_{n-i}(x)$ are finite dimensional, $\cF^{\hspace{.5pt}C}\m(M^\cdt) \m\in\! C^b({\rm gproj}\La)$. Hence, we have a functor $\cF^C\!\!:\! C^b({\rm gmod}^{b\hspace{-3pt}}\La) \m\to\m C^b({\rm gproj}\La^!)$ which, as seen in the proof of Theorem \ref{F-diag}, induces a commutative diagram as stated in Statement (1) with ${\rm GMod}^{b\hspace{-3pt}}\La$, ${\rm GProj}\La^!$ and ${\rm GMod}\La^!$ replaced by ${\rm gmod}^{b\hspace{-3pt}}\La$, ${\rm gproj}\La^!$ and ${\rm gmod}\La^!$, respectively. 
The proof of the theorem is completed.

\vspace{2pt}

In the locally bounded cases, the Koszul functors induce derived functors between bounded derived categories of finitely piece-supported graded modules.

\begin{Theo}\label{F-diag-fd-1}

Let $\La=kQ/R$ be a qudratic algebra with $Q$ a locally finite quiver.

\begin{enumerate}[$(1)$]

\vspace{-2pt}

\item \hspace{-2pt} If $\La^!$ is locally left bounded, the right Koszul functor $\cF\!:\hspace{-1pt} \GrLa \!\to\m C(\GrLa^!)$ induces a commutative diagram of functors \vspace{-4pt}
$$\xymatrixrowsep{18pt}\xymatrixcolsep{20pt}\xymatrix{
C^{\hspace{.5pt}b}\m({\rm GMod}^{\hspace{.3pt}b\hspace{-3pt}}\La) \ar[d]_{\cF^C_b} \ar[r]
& K^{\hspace{.5pt}b}\m({\rm GMod}^{\hspace{.3pt}b\hspace{-3pt}}\La)\ar[d]_{\cF^K_b} \ar[r]
& D^{\hspace{.5pt}b}\m({\rm GMod}^{\hspace{.3pt}b\hspace{-3pt}}\La) \ar[d]_{\cF^D_b} \\
C^{\hspace{.5pt}b}\m({\rm GProj}\La^!)\ar[r]
& K^{\hspace{.5pt}b}\m({\rm GProj}\La^!) \ar[r]
& D^{\hspace{.5pt}b}\m({\rm GMod}^{\hspace{.3pt}b\hspace{-3pt}}\La^!).}$$

\vspace{-3pt}

\item \hspace{-2pt} If $\La^!$ is locally right bounded, the left Koszul functor $\cG\!:\hspace{-1pt} \GrLa \!\to\m C(\GrLa^!)$ induces a commutative diagram  of functors \vspace{-4pt}
$$\xymatrixrowsep{18pt}\xymatrixcolsep{20pt}\xymatrix{
C^{\hspace{.5pt}b}\m({\rm GMod}^{\hspace{.3pt}b\hspace{-3pt}}\La) \ar[d]_{\cG^{\hspace{.3pt}C}_b} \ar[r]
& K^{\hspace{.5pt}b}\m({\rm GMod}^{\hspace{.3pt}b\hspace{-3pt}}\La)\ar[d]_{\cG^K_b} \ar[r]
& D^{\hspace{.5pt}b}\m({\rm GMod}^{\hspace{.3pt}b\hspace{-3pt}}\La) \ar[d]_{\cG^D_b} \\
 C^{\hspace{.5pt}b}\m({\rm GInj}\La^!)\ar[r]
& K^{\hspace{.5pt}b}\m({\rm GInj}\La^!) \ar[r]
& D^{\hspace{.5pt}b}\m({\rm GMod}^{\hspace{.3pt}b\hspace{-3pt}}\La^!).}$$

\vspace{-6pt}

\item In the above two statements, we can simultaneously replace ${\rm GMod}^{\hspace{.3pt}b\hspace{-3pt}}\La$, ${\rm GMod}^{\hspace{.3pt}b\hspace{-3pt}}\La^!\!,$ ${\rm GProj}\La^!$ and ${\rm GInj}\La^!$ by ${\rm gmod}^{\hspace{.3pt}b\hspace{-3pt}}\La$, ${\rm gmod}^{\hspace{.3pt}b\hspace{-3pt}}\La^!$, ${\rm gproj}\La^!$ and ${\rm ginj}\La^!\!,$ respectively.

\end{enumerate} \end{Theo}

\noindent{\it Proof.} Let $\La^!$ be locally left bounded. Then, the $P_x^!$ with $x\in Q_0$ are all finite dimensional. Hence, ${\rm GProj}\La^!\subseteq {\rm GMod}^{\hspace{.3pt}b\hspace{-3pt}}\La^!$ and ${\rm gproj}\La^!\subseteq {\rm gmod}^{\hspace{.3pt}b\hspace{-3pt}}\La^!$. Therefore, $K^b({\rm GProj}\La^!)$ and $K^b({\rm gproj}\La^!)$ are full triangulated subcategories of $D^b({\rm GMod}^{\hspace{.3pt}b\hspace{-3pt}}\La^!)$ and $D^b({\rm gmod}^{\hspace{.3pt}b\hspace{-3pt}}\La^!)$, respectively. Now, by the argument used in the proof of Theorem \ref{F-diag-fd}, we can prove the three statements stated in the theorem. The proof of the theorem is completed.



\vspace{3pt}

\noindent{\sc Remark.} In the sequel, we shall call $\cF^D_b$ and $\cG^D_b$ the {\it right} and the {\it left bounded derived Koszul functors}, respectively.

\vspace{2pt}

Next, we shall show that all derived Koszul functors are triangle equivalences in the Koszul case. We start with the following important property of Koszul functors; see \cite[(1.2.6)]{BGS} and \cite[Theorem 30]{MOS}.

\begin{Lemma}\label{inj-im}

Let $\La=kQ/R$ be a Koszul algebra with $Q$ a locally finite quiver. Given
$a\in Q_0$, the graded simple $\La^!$-module $S^{\hspace{.3pt}!}_a$ has $\cF(I_a\m)^\cdt$ as a linear projective resolution and $\cG(P_a\m)^\cdt$ as a colinear injective coresolution.

\end{Lemma}

\noindent{\it Proof.} Fix $a\in Q_0$. Since $\La^!$ is Koszul; see (\ref{Opp-Koszul}), $S_a^!$ has a colinear injective coresolution $\mathcal{I}_{a^!}^\pdt$; see (\ref{Koz-cplx-dual}). Since $(\La^!)^!=\La$; see (\ref{q-dual}), by definition, $\cG(P_a)^\pdt =\mathcal{I}_{a^!}^\pdt$. Next, by Lemma \ref{k-cplx-iso}, $S_a^!$ has a linear projective resolution $\mathcal{P}_{a^!}^\pdt$ as follows$\,:$\vspace{-6pt}
$$\xymatrixrowsep{18pt}\xymatrix{\cdots \ar[r] & \mathcal{P}_{a^!}^{-n} \ar[r]^-{\ell^{-n}} & \mathcal{P}_{a^!}^{1-n} \ar[r] &
\cdots \ar[r] & \mathcal{P}_{a^!}^{-1} \ar[r]^{\ell^{-1}} & \mathcal{P}_{a^!}^0 \ar[r] & 0\ar[r] & \cdots,} \vspace{-4pt}
$$
where $\mathcal{P}_{a^!}^{-n} \!\m =\! \oplus_{x\in Q_0} \hspace{-.8pt} P_x^!\nla -n\nra \otimes D(e_a\La_{n}e_x\hspace{-.5pt})$ and $\mathcal{P}_{a^!}^{1-n} \!=\!\oplus_{y\in Q_0} \hspace{-1pt} P_y^!\nla 1-n\nra \otimes D(e_a\La_{n-1}e_y \hspace{-.5pt})$ with
$\ell^{-n}=({\textstyle\sum}_{\alpha\in Q_1(x,y)} P[\bar \alpha^!] \otimes DP[\bar\alpha])_{(y,x)\in Q_0\times Q_0}\vspace{1pt}$.

On the other hand, $\cF(I_a)^{-n}=\oplus_{x\in Q_0} P_x^!\nla -n\nra \otimes D(e_x \La_n^{\rm o} e_a)$, for all $n\in \Z$. In particular,
$\cF(I_a)^{-n}=\mathcal{P}_{a^!}^{-n}$, for all $n<0$. Fix some integer $n\ge 1$. Write $\cF(I_a)^{1-n}=\oplus_{y\in Q_0} P_y^!\nla 1 - n \nra \otimes D(e_y \La_{n-1}^{\rm o} e_a)$. Since $I_a(\bar \alpha)=DP_a^{\rm o}(\bar\alpha^{\rm o})$ for any $\alpha\in Q_1(x,y)$, we see that $d_{\cF(I_a)}^{-n}=(\textstyle\sum_{\alpha\in Q_1(x,y)} P[\bar\alpha^!]\otimes DP_a^{\rm o}(\bar\alpha^{\rm o}))_{(y,x)\in Q_0\times Q_0}$.

\vspace{1.5pt}

Consider the $k$-linear isomorphism
$\theta^n_{x}: e_a \La_n e_x \to e_x \La_n^{\rm o} e_a: \xi\mapsto \xi^{\rm o}$, which induces a $k$-linear isomorphism
$D\theta^n_x: D(e_x \La_n^{\rm o} e_a)\to D(e_a \La_n e_x)$. \vspace{.5pt} Given $\alpha\in Q_1(x, y)$, since
$\theta^{n}_x \circ P[\bar\alpha] = P_a^{\hspace{.4pt} \rm o}(\bar\alpha^{\rm o}) \circ \theta^{n-1}_y$, we have
$DP[\bar\alpha] \circ D \theta^{n}_x = D\theta^{n-1}_y \circ  DP_a^{\rm o}(\bar\alpha^{\rm o}).$
Thus, the graded $\La$-linear isomorphisms
$\oplus_{x\in Q_0} ({\rm id}\otimes D\theta^n_{x})$ with $n\in \Z$ give rise to a complex isomorphism $\cF(I_a\hspace{-1pt})^\cdt\cong \mathcal{P}_{a^!}^\pdt$.
The proof of the lemma is completed.

\vspace{2pt}

More generally, applying the left Koszul functor and the right complex Koszul functor yields graded projective resolutions for any graded modules.

\begin{Lemma}\label{natural-tran}

Let \vspace{.5pt} $\La=kQ/R$ be a Koszul algebra with $Q$ a locally finite quiver. Given $M \m\in \m \GrLa$, we have a natural quasi-isomorphism $\eta_{\hspace{-.8pt}_M}^\cdt\hspace{-3.5pt}:
(\cF^C \hspace{-2pt} \circ \cG)(M)^\ydt \m\to\m M.$

\vspace{-2pt}

\end{Lemma}

\noindent{\it Proof.} Consider $M \m\in \m \GrLa$. By definition, $(\cF^{\hspace{.4pt}C}\!\m\circ\m \cG)(M)^\cdt \!=\! \T(\cF(\cG(M)^\cdt\hspace{-1pt})^\cdt\hspace{-.8pt})$. Given $n\in \Z$, by Lemma \ref{Comp_Kfunctors}(1), the $n$-diagonal of $\cF(\cG(M)^\cdt\hspace{-.8pt})^\cdt$ consists of \vspace{-2pt}
$$\cF(\cG(M)^i)^{n-i} = \oplus_{a, x\in Q_0} P_a \m \langle n-i \rangle \otimes (I_x^!)_n(a)
\otimes M_i(x); \; i\in \Z, \vspace{-3pt} $$
and so, $(F^{\hspace{.4pt}C}  \circ G)(M)^n = \oplus_{i\in \mathbb{Z}; \hspace{.7pt} a, x\in Q_0} P_a\tla n\!-\!i\tra \otimes (I_x^!)_n(a) \otimes M_i(x).$ \vspace{.5pt} In particular,
$(\cF^{\hspace{.4pt}C} \!\m \circ\cG)(M)^n=0$ for $n>0$. We divide the rest of the proof into two statements.

\vspace{2pt}

{\sc Statement 1.} {\it If $n<0$, then ${\rm H}^n(\m(\cF^{\hspace{.4pt}C} \!\m\circ\m \cG)(M)^\cdt\m)=0$.}

\vspace{1pt}

Fix an integer $n<0.$ For any $i\in \Z$, recall that the $i$-th column of $\cF(\cG(M)^\cdt\m)^\cdt$ is $\mathfrak{t}^i(\cF(\cG(M)^i\m)^\cdt)=\oplus_{x\in Q_0} \mathfrak{t}^i\m(\m \cF(I_x^!\tla i \tra)^\cdt\m)\otimes M_i(x),$ where $\mathfrak{t}$ is the twist functor. In view of Lemmas \ref{Ksl-Fun-Homology} and \ref{inj-im}, we see that
\vspace{-3pt}$${\rm H}^{n-i}(\mathfrak{t}^i\m(\m\cF(I_x^!\tla i \tra)^\cdt)) \!\cong \!
{\rm H}^{n-i}(\mathfrak{t}^i\m(\m\cF(I_x^!\tla i \tra)^\cdt \tla i \tra)) = {\rm H}^{n-i}(\cF(I_x^!)[i])=
{\rm H}^n(\cF(I_x^!)^\cdt)= 0. \vspace{-3pt} $$ So,
${\rm H}^{n-i}(\mathfrak{t}^i(\cF(\cG(M)^i)^\cdt))\cong \oplus_{x\in Q_0}{\rm H}^{n-i}(\mathfrak{t}^i\m(\m\cF(I_x^!\tla i \tra)^\cdt))\otimes M_i(x)=0.$ Fix $p\in \Z$. Consider the double complex $\cF(\cG(M)^\ydt)^\ydt_p\hspace{.4pt}$, whose $i$-th column is $\mathfrak{t}^i(\cF(\cG(M)^i)^\ydt)_p$ with ${\rm H}^{n-i}(\mathfrak{t}^i(\cF(\cG(M)^i)^\ydt)_p)=0$, and whose $n$-diagonal consists of \vspace{-3pt}
$$\cF(\cG(M)^i)^{n-i}_{p}=\oplus_{a, x\in Q_0} \La_{p+n-i}e_a \otimes D(e_{a}\hspace{2pt}\hat{\hspace{-3pt}\La}_{-n}e_x) \otimes M_{i}(x); \; i\in \Z.$$ So,
$\cF(\cG(M)^\ydt)^\ydt_p$ is $n$-diagonally bounded-above. In view of Lemma \ref{Homology-zero}, we see that ${\rm H}^n((\cF^C \hspace{-2.5pt} \circ \hspace{-1.5pt} \cG)(M)^\ydt_p)={\rm H}^n(\T(\cF(\cG(M)^\ydt)^\ydt_p))=0.$ Therefore, ${\rm H}^n\m((\cF^C \hspace{-2.5pt} \circ \hspace{-1.5pt} \cG)(M)^\cdt\m)=0$. This establishes Statement (1).

\vspace{1pt}

It remains to show that ${\rm H}^0((\cF^C \!\!\circ\m \cG)(M)^\cdt)$ is naturally isomorphic to $M$. Recall that $\La^!=
\{\bar\gamma^{\hspace{.4pt}!} \nid \gamma \in kQ\}$, where $\bar\gamma^{\hspace{.4pt}!}=\gamma^{\rm o}+R^!$; and
$\hspace{2.5pt}\hat{\hspace{-3.5pt}\La}
=kQ/(R^!)^{\rm o}=
\{\hspace{2pt}\hat{\hspace{-2pt}\gamma} \nid \gamma \in kQ\}$, where $\hspace{2pt}\hat{\hspace{-2pt}\gamma}=\gamma+(R^!)^{\rm o}$. Observe that the $1$-diagonal of $\cF(\cG(M)^\cdt)^\cdt$ is null. Since $(I_x^!\otimes M)(\bar\alpha^!)=I_x^!(\bar\alpha^!)\otimes {\rm id}_M$ for $\alpha\in Q_1(a, x)$, the $0$-diagonal and the $(-1)$-diagonal of $\cF(\cG(M)^\ydt)^\ydt$ are illustrated as \vspace{-3pt}
$$\xymatrixcolsep{22pt}\xymatrixrowsep{18pt}\xymatrix{
\oplus_{b\,\in Q_0} P_b \sla\!-i\tra\otimes (I_b^!)_0(b)\otimes M_i(b)\\
\oplus_{a, x\in Q_0} P_a\mla -i\!-\!1\mra\otimes (I_x^!)_{-1}(a) \hspace{-2pt} \otimes \hspace{-2pt} M_i(x) \ar[u]^-{v^{i, -i-\hspace{-2pt}1}} \ar[r]^-{\hspace{1pt}h^{i,-i-\hspace{-1.5pt}1}} &
\oplus_{c \in Q_0} P_{c}\mla -i\!-\!1\mra \hspace{-2pt} \otimes \hspace{-2pt} (I_c^!)_0(c)\otimes M_{i+1}(c),
}\vspace{-3pt}$$
where $v^{i,-i-1}=(v^{i,-i-1}(b, a,x))_{(b, a, x)\in Q_0 \times  Q_0 \times Q_0}\vspace{1.5pt}$ with \vspace{-2pt}
$$v^{i,-i-1}(b,a, x)=\left\{\begin{array}{ll}
{\textstyle\sum}_{\alpha\in Q_1(x, a)}(-1)^i P[\bar{\alpha}] \otimes I_x^!(\bar\alpha^{!})\otimes {\rm id}, & \mbox{ if } b=x; \vspace{-1pt} \\
0, & \mbox{ if } b\ne x, \end{array}\right. \vspace{-2pt} $$
and $h^{i,-i-1}=(h^{i,-i-1}(c,a, x))_{(c,a, x)\in Q_0 \times Q_0\times Q_0}$ \vspace{1pt} with \vspace{-3pt}
$$h^{i,-i-1}(c,a, x)=\left\{\begin{array}{ll}
{\textstyle\sum}_{\alpha\in Q_1(x, a)} {\rm id}\otimes I[\bar\alpha^!]\otimes M(\bar \alpha),  & \mbox{ if } c=a; \vspace{-1pt} \\
0, & \mbox{ if } c\ne a.
\end{array}\right.
\vspace{-3pt}$$

In particular, \vspace{.5pt} $(\cF^C \!\m\circ\m \cG)(M)^{-1} \hspace{-2pt} = \oplus_{i\in \mathbb{Z}; \hspace{.7pt} a, x\in Q_0} \hspace{.5pt} P_a\mla -i\!-\!1\mra\otimes (I^!_x)_{-1}(a) \otimes M_{i}(x),$ where $(I_x^!)_{-1}(a)= D(e_{a}\hspace{2pt}\hat{\hspace{-3pt}\La}_1 e_x)$ has a $k$-basis $\{\hat\beta^\star \nid \beta\in Q_1(x, a)\}$, that is the dual basis of $\{\hat \beta \nid \beta\in Q_1(x,a)\}$. Moreover, $(\cF^C \!\m\circ\m \cG)(M)^0 = \oplus_{(i, b)\in \mathbb{Z}\times Q_0} \hspace{.5pt} P_b\tla\!-i\tra\otimes (I^!_b)_0(b) \otimes M_i(b),$ where \vspace{1pt} $(I_b^!)_0(b)= D(e_b\hspace{2.5pt}\hat{\hspace{-3pt}\La}_0e_b)$ has a $k$-basis $\{\hat{e}_b^\star\}$.

\vspace{2pt}

{\sc Statement 2.} {\it There exists a natural graded $\La$-linear epimorphism \vspace{-3.5pt}
$$\eta_{_M}\hspace{-2.5pt}: \! (\cF^C \hspace{-4pt} \circ \hspace{-1pt} \cG)(M)^0 \!\to \! M\hspace{-3pt}: \! {\textstyle \sum}_{(i, b)\in \mathbb{Z}\times Q_0}\m u_{i,b} \m \otimes \m \hat{e}^\star_b \m \otimes \m m_{i,b} \mapsto \!\m {\textstyle\sum}_{i\in \mathbb{Z}; \hspace{.7pt} b\in Q_0}\m (\m-\hspace{-1.5pt} 1\m)^{\frac{i(i+1)}{2}} \! u_{i,b} m_{i,b},\vspace{-3pt}$$ 
such that $\eta_{_M}^0 \!\! \circ \! d^{-1}=0$, where $d^{-1}$ is the differential of degree $-1$ of $(\cF^C \!\m\circ\m \cG)(M)^\cdt$.}

The existence and the naturalness of of $\eta_{_M}$ are evident. Let $\omega\in (\cF^C \hspace{-3pt} \circ \hspace{0.3pt} \cG)(M)^{-1}$. We may assume that
$\omega \in P_a\mla -i\!-\!1\mra\otimes I_x^!\tla i \tra_{-i-1}(a)\otimes M_i(x)$ for some $i\in \Z$ and $a, x\in Q_0$. Further, we may assume that
$\omega=u_0\otimes \hat\beta^\star_0\otimes m_0,$ for some $u_0\in P_a\mla -i\!-\!1\mra$, $\beta_0\in Q_1(x, a)$ and $m_0\in M_i(x)$.
Write $\hat{P}_x=\hspace{2pt}\hat{\hspace{-3pt}\La}e_x$. For $\alpha\in Q_1(x, a)$, since $(\bar\alpha^!)^{\rm o}=\hat\alpha$, we obtain
$I_x^!(\bar\alpha^!)=D\hat{P}_x(\hat\alpha)$ and $I_x^![\bar\alpha^!]=DP[\hat\alpha]$. Thus, $I_x^!(\bar\alpha^!)(\hat{\beta}^\star_0)(e_x) 
=\hat{\beta}^\star_0(\hat{\alpha})$ and $I[\bar\alpha^!](\hat{\beta}^\star_0)(e_a) 
=\hat{\beta}^\star_0(\hat\alpha).$
Hence, $I_x^!(\bar\alpha^!)(\hat{\beta}^\star_0)=I[\bar\alpha^!](\hat{\beta}^\star_0)=\hat{e}_x^\star$ in case $\alpha=\beta_0$; and otherwise, $I_x^!(\bar\alpha^!)(\hat{\beta}^\star_0)=I[\bar\alpha^!](\hat{\beta}^\star_0)=0$. This yields \vspace{-4pt}
$$\hspace{-2pt} \begin{array}{rclr}
d^{-1} (\omega)  
& = & (-1)^i{\textstyle\sum}_{\alpha\in Q_1(x, a)} (P[\bar\alpha]\otimes I_x^!(\bar\alpha^!)\otimes {\rm id}) (u_0\otimes \hat{\beta}^\star_0\otimes m_0) \vspace{2.5pt}\\
& & + \hspace{2pt} {\textstyle\sum}_{\alpha\in Q_1(x, a)} ({\rm id} \otimes I[\bar\alpha^!] \otimes M(\bar\alpha)) (u_0\otimes \hat{\beta}^\star_0\otimes m_0) & \hspace{30pt} (*) \vspace{2.5pt}\\
& = & (-1)^i (u_0\bar\beta_0) \otimes \hat{e}_x^\star\otimes m_0 + u_0 \otimes \hat{e}_a^\star \otimes (\bar\beta_0 m_0).
\end{array}\vspace{-3pt}$$
Since $u_0\in P_a\mla -i\!-\!1\mra$ and $u_0\bar\beta_0\in P_x\tla\!-i\tra$, we obtain \vspace{-3pt}
$$(\eta_{_M} \!\!\circ\m d^{-1}) (\omega) = (-1)^{\frac{i(i+1)}{2}+i} (u_0 \bar\beta_0 m_0) +  (-1)^{\frac{(i+1)(i+2)}{2}} (u_0\bar\beta_0 m_0) = 0.
\vspace{-3pt}$$
This establishes Statement 2.

To conclude, we need to verify that 
${\rm Ker}(d^{-1})\subseteq {\rm Im}(\eta_{_M})$. Fix $\omega\in {\rm Ker}(\eta_{_M})$. Since $\eta_{_M}$ is graded, we may assume that there exists $(p,a)\in \Z\times Q_0$ such that \vspace{-3pt}
$$\omega\in (\cF^C \!\circ \cG)(M)^0_p(a)=\oplus_{i\le p; \hspace{.4pt} x \hspace{.4pt} \in Q_0} e_aP_x\sla -i \tra_p \otimes (I_x^!)_0(x)\otimes M_i(x),
\vspace{-4pt}$$ where $e_aP_x\sla -i \tra_p=e_a\La_{p-i}e_x.$ Then, we may find some $i_s\le \cdots\le i_2\le i_1= p$ and $x_1, \ldots, x_s\in Q_0$ such that $\omega=\textstyle\sum_{j=1}^s \bar{\gamma\hspace{1.6pt}}\hspace{-2pt}_j \otimes \hat{e}_{x_j}^\star\otimes m_j,$ where $\gamma_j\in Q_{p-i\hspace{-.5pt}_j}(x_j, a)$ such that
the $\bar{\gamma\hspace{1.6pt}}\hspace{-2pt}_j$ are pairwise distinct and $m_j\in M_{i_j}(x_j)$. In particular, $\gamma_1=\varepsilon_a$. We shall proceed by induction on the minimal integer $n_\omega$ for which $\omega$ can be written in this form and $n_\omega=\sum_{j=1}^s(p-i\hspace{-.8pt}_j)$.

If $n_\omega=0$, then $s=1$ and $m_1=\pm \hspace{.4pt} \eta_{_M}(\omega)=0$, and hence, $\omega=0$. Suppose that $n_\omega>0$. Since $\gamma_j\ne \gamma_1=\varepsilon_a$, we may write $\gamma_j=\sigma_j\beta_j$ with $\beta_j\in Q_1(x_j, y_j)$ and $\sigma_j\in Q_{p-i_j-1}(y_j, a)$, for $2\le j\le s.$ Set $\sigma=\sum_{j=2}^s (-1)^{i\hspace{-1pt}_j} \bar\sigma_j \otimes \hat \beta^\star _j \otimes m_j,$ where $\bar\sigma_j\in P_{y_j}\mla -i_j\!-\!1\mra_p$.
In view of the equations $(*)$, we obtain \vspace{-2pt}
$$\begin{array}{rcl}
d^{-1}(\sigma) &=& {\textstyle\sum}_{j=2}^s ((-1)^{2i\hspace{-1pt}_j} (\bar\sigma_j \, \bar\beta_j) \otimes \hat{e}_{x_j}^\star \otimes m_j +  (-1)^{i\hspace{-1pt}_j} \bar\sigma_j \otimes \hat{e}_{y_j}^\star \otimes (\bar\beta_j \, m_j)) \vspace{2pt}\\
 &=& \omega + e_a\otimes \hat{e}_a^\star \otimes (-m_1) +  \sum_{j=2}^s \bar\sigma_j \otimes \hat{e}_{y_j}^\star \otimes (-1)^{i\hspace{-1pt}_j} (\bar\beta_j \, m_j),
\end{array} \vspace{-4pt}$$
Put $\omega'=d^{-1}(\sigma)-\omega=e_a\otimes \hat{e}_a^\star \otimes (-m_1) +  \sum_{j=2}^s \bar\sigma_j \otimes \hat{e}_{y_j}^\star \otimes (-1)^{i\hspace{-1pt}_j} (\bar\beta_j \, m_j)$. Then, $\omega'\in {\rm Ker}(\eta_{_M}\!)$ with $n_{\omega'}<n_\omega$. Thus $\omega'\in {\rm Im}(d^{-1})$, and hence, $\omega\in {\rm Im}(d^{-1})$. The proof of the lemma is completed.

\vspace{2pt}

Similarly, applying the right Koszul functor and the left complex Koszul functor yields graded injective coresolutions for bounded-above graded modules.

\begin{Lemma}\label{natural-tran-2}

Let \vspace{.5pt} $\La=kQ/R$ be a Koszul algebra with $Q$ a locally finite quiver. Given $M\in {\rm GMod}^-\hspace{-3.5pt}\La$, we have a natural quasi-isomorphism $\zeta_{\hspace{-.5pt}_M}^\cdt\hspace{-4pt}:\hspace{-1.5pt} M \!\to\! (\cG^{\hspace{.3pt}C}\!\circ\m \cF)(M)^\cdt\!$.

\end{Lemma}

\vspace{-2pt}

\noindent{\it Proof.} Let $M \!\in {\rm GMod}^-\hspace{-3.5pt}\La$, say $M_{i}=0$ for all $i\ge r$, where $r\in \Z$. By definition, $(\cG^{\hspace{.3pt}C}\!\circ \cF)(M)^\cdt=\T(\cG(\cF(M)^\pdt\m)^\pdt)$. Given any $n\in \Z$, by Lemma \ref{Comp_Kfunctors}(2), we see that the $n$-diagonal of $\cG(\cF(M)^\ydt)^\ydt$ consists of  \vspace{-2pt}
$$\cG(\cF(M)^i)^{n-i} =\oplus_{\hspace{.7pt} x, a\in Q_0} I_x \langle n\m -\m i \rangle \otimes e_x\La_n^!e_a\otimes M_i(a); \; i\in \Z, \vspace{-2pt}
$$ and so,
$(\cG^{\hspace{.3pt}C}\!\circ \cF)(M^\cdt)^n=\oplus_{i\in \Z; \, x, a\in Q_0} I_x \langle n\m -\m i \rangle \otimes e_x\La_n^!e_a\otimes M_i(a).$
In particular, $\cG(\cF(M)^\ydt)^\ydt$ is diagonally bounded-above and $(\cG^{\hspace{.3pt}C} \hspace{-2.5pt}\circ\hspace{-1pt} \cF)(M)^n=0$ for $n<0.$

\vspace{2pt}

{\sc Statement 1.} {\it Consider $\mathfrak{t}^i(\cG(\cF(M)^i)^\cdt\hspace{.4pt})$, the $i$-th column of $\cG(\cF(M)^\ydt)^\ydt$. Given any $n\in \Z$, we have
${\rm H}^{n-i}(\mathfrak{t}^i(\cG(\cF(M)^i)^\cdt\hspace{.4pt})) \cong \oplus_{a\in Q_0}\, {\rm H}^n(\cG(P_a^{\hspace{.4pt}!})^\cdt\hspace{-.5pt}) \otimes M_i(a)$.
}

\vspace{1pt}

Indeed fix $n, i\in \Z$. Then, $\cF(M)^i=\oplus_{a\in Q_0} P_a^{\hspace{.4pt}!}\tla i \tra \otimes M_i(a)$. It follows easily that $\mathfrak{t}^i(\cG(\cF(M)^i)^\cdt) = \oplus_{a\in Q_0} \mathfrak{t}^i(\cG(P_a^{\hspace{.4pt}!}\tla i \tra)^\cdt) \otimes M_i(a).$ In view of Lemma \ref{Ksl-Fun-Homology}, we see that \vspace{-2pt}
$$
{\rm H}^{n-i}(\mathfrak{t}^i(\cG(P_a^{\hspace{.4pt}!}\tla i \tra)^\cdt))
\cong {\rm H}^{n-i}(\mathfrak{t}^i(\cG(P_a^{\hspace{.4pt}!}\tla i \tra)^\cdt \tla i \tra)) \cong {\rm H}^{n-i}(\cG(P_a^{\hspace{.4pt}!})^\cdt[i]) \cong  {\rm H}^n(\cG(P_a^{\hspace{.4pt}!})^\cdt).\vspace{-2pt}$$
This establishes Statement 1.

\vspace{2.5pt}

{\sc Statement 2.} {\it If $n>0$, then ${\rm H}^n((\cG^{\hspace{.3pt}C} \hspace{-2.5pt}\circ\hspace{-1pt} \cF)(M)^\cdt)=0.$}

Fix $n>0$. By Statement (1) and Lemma \ref{inj-im}, ${\rm H}^{n-i}(\mathfrak{t}^i(\cG(\cF(M)^i)^\cdt)) =0$ for all $i\in \Z$. Since $\cG(\cF(M)^\ydt)^\ydt$ is diagonally bounded-above, we deduce from  Lemma \ref{Homology-zero} that
${\rm H}^n((\cG^{\hspace{.3pt}C}\hspace{-2.5pt} \circ \hspace{-1.5pt}\cF)(M)^\cdt)={\rm H}^n(\T(\cG(\cF(M)^\cdt)^\cdt))=0$.
This establishes Statement 2.

\vspace{1pt}

It remains to construct a natural graded isomorphism $M\to {\rm H}^0((\cG^{\hspace{.3pt}C}\hspace{-2.5pt} \circ \hspace{-1pt} \cF)(M)^\cdt)$.
Note that \vspace{.5pt} $(\cG^{\hspace{.3pt}C} \hspace{-2.5pt} \circ \hspace{-1pt} \cF)(M)^0=\oplus_{i\in \mathbb{Z};\, a\in Q_0}\, I_a\tla\!-i\tra\otimes e_a\La_0^!e_a\otimes M_i(a).$ Given $(i, a)\in \Z\times Q_0$, we shall construct a morphism $f^i_a: M\to I_a\tla\!-i\tra\otimes e_a\La_0^!e_a\otimes M_i(a)$ in $\GrLa$. For this, we shall first define a $k$-linear map $f_{a,j}^i: M_j\to  I_a\tla\!-i\tra_j\otimes e_a\La_0^!e_a\otimes M_i(a)$ for every $j\in \Z$, where $I_a\tla\!-i\tra_j=D(\La_{i-j}^{\rm o} e_a).$ Indeed, we set $f_{a,j}^i=0$ for $j>i$. Fix $j$ with $j\le i$. We have a $k$-linear map \vspace{-2.5pt}
$$\psi_{a,j}^i: M_j\to \Hom_k(\La_{i-j}^{\rm o} e_a, \, e_a\La_0^!e_a\otimes M_i(a)): w\mapsto \psi_{a,j}^{i}(w),\vspace{-2pt}$$
where $\psi_{a,j}^{i}(w)$ sends $\bar{\hspace{-1pt}\gamma}^{\rm o}$ to $e_a\otimes \bar{\hspace{-1pt}\gamma} w,$ for all $\gamma\in kQ_{i-j}(-,a).$ Since $\La_{i-j}^{\rm o} e_a$ is finite dimensional, in view of Corollary \ref{Cor 1.2}(1), we have a $k$-linear isomorphism \vspace{-2pt}
$$\theta_{a,j}^{i}: D(\La_{i-j}^{\rm o} e_a)\otimes e_a\La_0^!e_a\otimes M_i(a) \to \Hom_k(\La_{i-j}^{\rm o} e_a, \, e_a\La_0^!e_a\otimes M_i(a)).\vspace{-4pt}$$

Now, put $f_{a,j}^i=(\theta_{a,j}^{\hspace{.6pt}i})^{-1}\circ \psi_{a,j}^{i}: M_j\to I_a\tla\!-i\tra_j \otimes e_a\La_0^!e_a\otimes M_i(a),$
which can be computed in the following way.

\vspace{1pt}

{\sc Statement 3.} \vspace{0pt} {\it Let $\{\bar{\hspace{-1pt}\gamma}^{\rm o}_1,\ldots, \bar{\hspace{-1pt}\gamma}^{\rm o}_s\}$ with $\gamma_p\in kQ_{i-j}(-, a)$ be a $k$-basis of $\La_{i-j}^{\rm o}e_a$ with dual basis $\{\bar{\hspace{-1pt}\gamma}_1^{{\rm o}, \star}, \ldots, \bar{\hspace{-1pt}\gamma}_s^{{\rm o}, \star}\}$. Then $f_{a,j}^{i}(w)=\sum_{p=1}^s \bar{\hspace{-1pt}\gamma}_p^{{\rm o}, \star} \otimes e_a\otimes \bar{\hspace{-1pt}\gamma}_p w,$ for $w\in M_j$.}

\vspace{1pt}

Indeed, every $\bar{\hspace{-1pt}\gamma}^{\rm o}\in \La_{i-j}^{\rm o}e_a$ is written as $\bar{\hspace{-1pt}\gamma}^{\rm o}=\sum_{t=1}^s\, \lambda_t \, \bar{\hspace{-1pt}\gamma}^{\rm o}_t$ with $\lambda_t\in k$. Given $w\in M_j$, by the definition given in Corollary \ref{Cor 1.2}(1), we obtain \vspace{-2pt}
$$\textstyle\theta_{a,j}^i(\sum_{p=1}^s \bar{\hspace{-1pt}\gamma}_p^{{\rm o},\star} \otimes e_a\otimes \bar{\hspace{-1pt}\gamma}_p w)(\bar{\hspace{-1pt}\gamma}^{\rm o}) = e_a\otimes (\sum_{t=1}^s \lambda_t \hspace{.4pt} \bar{\hspace{-1pt}\gamma}_t^{\rm o}) \hspace{.4pt} w
= \psi_{a,j}^{i}(w)(\bar{\hspace{-1pt}\gamma}^{\rm o}).
\vspace{-2pt}$$
Thus, $\theta_{a,j}^{i}(\sum_{p=1}^s \bar{\hspace{-1pt}\gamma}_p^{{\rm o}, \star}  \otimes e_a\otimes \bar{\hspace{-1pt}\gamma}_p w)= \psi_{a,j}^{i}(w)$. 
This establishes Statement 3.

\vspace{2pt}

{\sc Statement 4.} \vspace{1pt} {\it Given $(i, a)\in \Z\times Q_0$, there exists a natural graded $\La$-linear morphism $f^i_a: M\to I_a\tla\!-i\tra\otimes e_a\La_0^!e_a\otimes M_i(a)$ such that $(f^i_a)_j=f_{a,j}^i$, for all $j\in \Z$.}

\vspace{.5pt}

Given $\alpha\in Q_1$ and $j\le i$, it is easy to see that we have a commutative diagram \vspace{-5pt}
$$\xymatrixcolsep{20pt}\xymatrixrowsep{22pt}\xymatrix{
\hspace{5pt} M_j \ar[r]^-{\psi_{a,j}^{i}} \ar[d]^{M(\bar\alpha)} & \hspace{1pt} \Hom((P_a^{\rm o})_{i-j}, \, e_a\La_0^!e_a \otimes M_i(a)) \ar[d]^{\Hom(P_a^{\rm o}(\bar\alpha^{\rm o})\hspace{-.4pt}, \hspace{.6pt} e_a\mathit\Lambda_0^!e_a \otimes M_i(a))} & I_a\tla\!-i\tra_j\otimes e_a\La_0^!e_a\otimes M_i(a) \ar[l]_-{\theta_{a,j}^{i}} \ar[d]^{I_a\tla\!-i\tra(\bar\alpha)\otimes {\rm id} \otimes {\rm id} } \hspace{8pt} \\
M_{j+1} \ar[r]^-{\psi_{a,j+1}^{i}} & \Hom(\hspace{-1pt}(P_a^ {\rm o})_{i-j-1}\hspace{-1pt}, \hspace{-.4pt}e_a\La_0^!e_a \otimes \hspace{-2pt} M_i(a)).
& I_a\tla\!-i\tra_{j+1}\otimes e_a\La_0^!e_a \otimes \hspace{-2pt} M_i(a)\ar[l]_-{\hspace{5pt}\theta_{a,j+1}^i}
}\vspace{-4.5pt}
$$
Therefore, $f_a^i$ is a graded $\La$-linear morphism.
%
Similarly, one can verify that $f^i_a$ is natural in $M$. This establishes Statement 4.

\vspace{1pt}

Fix $i\in \Z$. Given $a\in Q_0$, by Statement 4, we obtain a natural graded $\La$-linear morphism $g_a^i: M\to I_a\tla\!-i\tra\otimes e_a\La_0^!e_a\otimes M_i(a)$ with $(g^i_a)_j=(-1)^{\frac{(i-1)i}{2}}f_{a,j}^i$, which will be written as $g^i_{a,j}$, for all $j\in \Z$. Let $w=\sum_{x\in Q_0; \hspace{.4pt} j\in \mathbb{Z}} w_{x,j}\in M$ with $w_{x,j}\in e_xM_j$. If $g_a^i(w_{j, x})=g_{a, j}^i(w_{j,x})\ne 0$ for some $a\in Q_0$, then $f_{a, j}^i(w_{j,x})\ne 0$. Hence, $j\le i$, and by Statement 3, $kQ_{i-j}(x, a)\ne 0$. Since $Q$ is locally finite, $g_a^i(w)=0$ for all but finitely many $a\in Q_0$. Hence, we have a graded $\La$-linear morphism \vspace{-2.2pt}
$$g^i=(g^i_a)_{a\in Q_0}: M\to \cG(\cF(M)^i)^{-i}=\oplus_{a\in Q_0} I_a\tla\!-i\tra\otimes e_a\La_0^!e_a\otimes M_i(a)
\vspace{-2.5pt}$$ such that $g^i_j=(g^i_{a,j})_{a\in Q_0}: M_j\to \oplus_{a\in Q_0} (I_a\tla\!-i\tra)_j\otimes e_a\La_0^!e_a\otimes M_i(a)$, for all $j\in \Z$.

\vspace{2pt}

{\sc Statement 5.} {\it We have a natural graded $\La$-linear monomorphism \vspace{-1pt}
$$\zeta_{\hspace{-.5pt}_M}\m=(g^i)_{i\in \Z}: M\to (\cG^{\hspace{.3pt}C} \!\!\circ\m \cF)(M)^0=\oplus_{i\in \Z\,} \cG(\cF(M)^i)^{-i}.
\vspace{-1pt}$$}

Indeed, observe that $\cG(\cF(M)^i)^{-i}=0,$ for all $i\ge r$. Let $w=\sum_{j\in \Z} w_j\in M$ with $w_j\in M_j$. If $g^i(w_j)=\sum_{a\in Q_0} g^i_a(w_j)=g^i_{a,j}(w_j) \ne 0$ for some $i$, then $j\le i$, and hence, $j\le i<r$. As a consequence, $g^i(w)=0$ for all but finitely many $i\in \Z$. Thus, we obtain a graded $\La$-linear morphism $\zeta_{\hspace{-.5pt}_M}\! =(g^{i})_{i\in \Z}: M\to (\cG^{\hspace{.3pt}C}\!\!\circ\m \cF)(M)^0$, which is clearly natural in $M$. Assume that $\zeta_{\hspace{-.5pt}_M\m}(w)=0$, for some $w\in M_j$ with $j\in \Z$. In particular, $g^j(w)=0$, that is,
$g^j_j(w_j)=\sum_{a\in Q_0} g^j_{a,j}(w)=0$. Thus, $g_{a,j}^j(w)= 0$, and hence, $f_{a,j}^{j}(w)=0$, for all $a\in Q_0$. Since $\{e_a\}$ is a basis of $\La^{\rm o}_{j-j}e_a$, by Statement 3, ${e}_a^{{\rm o}, \star} \otimes e_a \otimes e_aw=0$, and hence, $e_aw=0$, for all $a\in Q_0$. That is, $w=0$. This implies that $\zeta_{\hspace{-.5pt}_M\m}$ is a monomorphism. Statement 5 is established.

\vspace{0pt}

Observe that the $(-1)$-diagonal of $\cG(\cF(M)^\ydt\hspace{.4pt})^\ydt$ is null, while the $0$-diagonal and the $1$-diagonal can be illustrated as follows$\,:$ \vspace{-5pt}
$$\xymatrixrowsep{20pt}\xymatrixcolsep{18pt}\xymatrix{
\hspace{-5pt}\oplus_{b\,\in Q_0}\, I_b\tla\!-i\tra \hspace{-1pt} \otimes \hspace{-2pt} e_b\La_0^! e_b \hspace{-1pt} \otimes \hspace{-2pt} M_i(b) \ar[r]^-{\hspace{2pt}h^{i,-i}} & \oplus_{a, x\in Q_0} I_x\tla\!-i\tra  \!\otimes\! e_x\La_1^!e_a \!\otimes\! M_{i+1}(a) \hspace{20pt}  \\
& \hspace{-4pt} \oplus_{c\, \in Q_0}  I_{c}\mla -i\!-\!1\mra  \!\otimes\! e_c\La_0^!e_c \!\otimes\! M_{i+1}(c) \ar[u]_-{v^{i+1, -i-1}},
}\vspace{-5pt} $$ where $h^{i,-i}=(h^{i,-i}(a,x,b))_{a,x,b\,\in Q_0}$ with \vspace{-1pt}
$$h^{i,-i}(a,x,b)=\left\{\begin{array}{ll}
\! {\textstyle\sum}_{\alpha\in Q_1(x,a)} \, {\rm id} \otimes P[\bar\alpha^!]\otimes M(\bar\alpha), & \mbox{ if } b=x; \vspace{0pt} \\
\! 0, & \mbox{ if } b\ne x, \end{array} \right.$$ and
$v^{i+1,-i-1} \!=\! (v^{i+1,-i-1}(a,\hspace{-1pt}x,\hspace{-1pt}c))_{a,x,c\,\in Q_0}$ with \vspace{-1pt}
$$v^{i+1,-i-1}\hspace{-1pt}(a, \hspace{-1pt} x, \hspace{-1pt} c) \!= \! \left\{\begin{array}{ll}
\hspace{-4pt} {\textstyle\sum}_{\alpha\in Q_1(x,a)}(-1)^{i+1} I[\bar\alpha] \!\otimes\! P_a^{\,!}(\bar \alpha^!) \!\otimes\! {\rm id}, &  \mbox{ if } c=a; \vspace{-0pt}  \\
\hspace{-4pt} 0, &  \mbox{ if } c\ne a.
\end{array} \right.$$

\vspace{0pt}

{\sc Statement 6}. {\it We have $d^{\hspace{.8pt}0} \!\circ \zeta_{\hspace{-.5pt}_M}= 0$, where $d^{\hspace{.8pt}0}\m$ denotes the differential of degree $0$ of the complex $(\cG^{\hspace{.3pt}C} \!\circ \cF)(M)^\cdt.$}

\vspace{1pt}

Indeed, it amounts to show, for any $i\in \Z$, that the diagram\vspace{-4pt}
$$\xymatrixcolsep{55pt}\xymatrixrowsep{20pt}\xymatrix{
\oplus_{x\,\in Q_0} I^!_x\tla\!-i\tra \hspace{-2pt}\otimes \hspace{-2pt} e_x\La_0e_x \hspace{-2pt}\otimes \hspace{-2pt} M_i(\hspace{-.5pt}x\hspace{-.5pt}) \ar[r]^-{\oplus_{x\in Q_0} \hspace{-.5pt}h^{i,-\hspace{-1pt}i}\hspace{-1pt}(\hspace{-.5pt}a,x,x\hspace{-.5pt})} & \oplus_{a, x\,\in Q_0} I^!_x\tla\!-i\tra \hspace{-2pt}\otimes \hspace{-2pt}  e_x\La_1e_a \hspace{-2pt}\otimes \hspace{-2pt} M_{i+1}(\hspace{-.5pt}a\hspace{-.5pt}) \\
M \ar[r]^-{(g_a^{i+1})_{a\in Q_0}} \ar[u]^-{(g_x^i)_{x\in Q_0}} & \oplus_{a\, \in Q_0} I^!_a\mla -i\!-\!1\mra \hspace{-2pt}\otimes \hspace{-2pt} e_a\La_0e_a \hspace{-2pt}\otimes \hspace{-2pt} M_{i+1}(a), \ar[u]^-{\oplus_{a\in Q_0}v^{i+1,-i-1}(a,x,a)}}
\vspace{-2.5pt}$$
is anti-commutative, or equivalently, the diagram \vspace{-2.5pt}
$$\xymatrixcolsep{45pt}\xymatrixrowsep{20pt}\xymatrix{
\oplus_{x\,\in Q_0} I_x\tla\!-i\tra\m_j \m\otimes\! e_x\La_0^!e_x \m\otimes\! M_i(x) \ar[r]^-{\oplus h^{i,-i}(a,x,x)\m_j} & \oplus_{a, x\,\in Q_0}I_x\tla\!-i\tra\m_j \m\otimes\! e_x\La_1^!e_a \m\otimes\! M_{i+1}(a) \\
M_j\ar[r]^-{(g_{a,j}^{i+1})_{a\in Q_0}} \ar[u]^-{(g_{x,j}^i)_{x\in Q_0}} & \oplus_{a\, \in Q_0} I_a\mla -i\!-\!1\mra\m_j \otimes e_a\La_0^!e_a\otimes M_{i+1}(a) \ar[u]_-{\oplus v^{i+1,-i-1}(a,x,a)_j}}
\vspace{-2.5pt}
$$
is anti-commutative for all $i,j\in \Z$, where $I_x\tla\!-i\tra_j=D(\La^{\rm o}_{i-j}e_x)$. This is evident in case $j>i$. Fix some $i\ge j$ and $a, x\in Q_0$. Then, we have a $k$-isomorphism $\theta^i_j: D(\La^{\rm o}_{i-j}e_x)\otimes e_x\La_1^!e_a \otimes M_{i+1}(a)\to \Hom_k(\La_{i-j}^{\rm o} e_x, \, e_x\La_1^!e_a\otimes M_{i+1}(a))\vspace{.5pt}$
as stated in  Corollary \ref{Cor 1.2}(1). Consider $\alpha\in Q_1(x, a)$ and $w\in M_j$. We choose a $k$-basis $\{\bar\delta^{\rm o}_1, \ldots, \bar\delta^{\rm o}_s\}$ for $\La_{i-j}^{\rm o} e_x$. By Statement 3, we obtain
\vspace{-2pt} $$\textstyle ({\rm id}\otimes P[\bar \alpha^!]\otimes M(\bar\alpha))(g^i_{x,j}(w))=(-1)^{\frac{(i-1)i}{2}} \sum_{p=1}^s \bar{\delta}_p^{{\rm o}, \star} \otimes \bar\alpha^! \otimes \bar\alpha\bar\delta_pw.\vspace{-2pt}$$
As a consequence, we see that \vspace{-1.5pt}
$$\textstyle \theta_j^i \hspace{-2pt}\left[({\rm id}\otimes P[\bar \alpha^!]\otimes M(\bar\alpha)) (g^i_{x,j}(w))\right]\hspace{-2pt} (\bar\delta_p^{\hspace{.4pt}\rm o})
=(-1)^{\frac{(i-1)i}{2}}(\bar \alpha^! \otimes \bar\alpha \bar\delta_pw), \; p=1, \ldots, s. \vspace{-1.5pt}
$$
On the other hand, for any $k$-basis $\{\bar\gamma^{\rm o}_1, \ldots, \bar\gamma^{\rm o}_t\}$ of $\La_{i+1-j}^{\rm o} e_a$, by Statement 3, \vspace{-2pt} $$\textstyle (I[\bar\alpha]\otimes P_a(\bar\alpha^!)\otimes {\rm id}) (g^{i+1}_{a,j}(w))=(-1)^{\frac{i(i+1)}{2}}\sum_{q=1}^t (\bar{\hspace{-1pt}\gamma}^{{\rm o}, \star}_q \circ P[\bar\alpha^{\rm o}])\otimes \bar \alpha^! \otimes \bar\gamma_q w. \vspace{-2pt}$$
And hence, for any $1\le p\le s,$ we obtain \vspace{-2pt}
$$\textstyle\theta^i_j\hspace{-2pt}\left[(I[\bar\alpha] \hspace{-1pt} \otimes \hspace{-1pt} P_a(\bar\alpha^!) \hspace{-1pt} \otimes {\rm id}) (g^{i+1}_{a,j}\hspace{-1pt}(w)) \right]\hspace{-2pt} (\bar\delta_p^{\hspace{.4pt}\rm o}) =
(\hspace{-1pt}-1\hspace{-.5pt})^{\frac{i(i+1)}{2}} \hspace{-2pt} \sum_{q=1}^t \hspace{-1.5pt} \bar{\hspace{-1pt}\gamma}^{{\rm o}, \star}_q(\bar\delta_p^{\hspace{.4pt}\rm o} \bar\alpha^{\rm o}) \cdot (\bar\alpha^! \otimes \bar{\hspace{-1pt}\gamma}_q w).$$
Fix $1\le p\le s$. If $\bar\delta_p^{\hspace{.4pt}\rm o}\, \bar\alpha^{\rm o}=0$, then $\bar\alpha\bar\delta_p=0$. In this case, we see trivially that \vspace{-2pt}
$$\textstyle\theta_j^i \hspace{-2pt} \left[({\rm id} \hspace{-2pt} \otimes \hspace{-2pt} P[\bar \alpha^!] \hspace{-2pt} \otimes \hspace{-2pt} M\hspace{-1pt}(\bar\alpha))(g^i_{j,x}(w))\right]\hspace{-2pt} (\bar\delta_p^{\hspace{.4pt}\rm o})= (\hspace{-1pt}-1\hspace{-.5pt})^i\theta^i_j\hspace{-2pt}\left[(I[\bar\alpha] \hspace{-2pt} \otimes \hspace{-2pt} P_a(\bar \alpha^!)  \hspace{-2pt} \otimes \hspace{-2pt} {\rm id}) (g^{i+1}_{j,a}\hspace{-1pt}(w)) \right]\hspace{-2pt} (\bar\delta_p^{\hspace{.4pt}\rm o}).\vspace{-1pt}$$

If $\bar\delta_p^{\hspace{.4pt}\rm o}\, \bar\alpha^{\rm o}\ne 0$, then $\La_{i+1-j}^{\rm o} e_a$ has a basis $\{\bar\gamma^{\rm o}_1, \ldots, \bar\gamma^{\rm o}_t\}$, where $\bar\gamma^{\rm o}_1=\bar\delta_p^{\hspace{.4pt}\rm o}\, \bar\alpha^{\rm o}$. Noting that $\bar\gamma_1=\bar\alpha\bar\delta_p$, we obtain \vspace{-1.5pt}
$$\hspace{-5pt}\begin{array}{rcl}
\theta^i_j\hspace{-2pt}\left(I[\bar\alpha] \hspace{-2pt} \otimes \hspace{-2pt} P_a(\bar\alpha^!) \hspace{-2pt} \otimes \hspace{-1.5pt}  {\rm id}) (g^{i+1}_{j,a}\hspace{-1pt}(w)) \right)\hspace{-2pt} (\bar\delta_p^{\hspace{.4pt}\rm o})
\hspace{-7pt} & = & \hspace{-7pt} (-1)^{\frac{i(i+1)}{2}}(\bar\alpha^! \hspace{-2pt} \otimes \hspace{-1pt} \bar\alpha \bar\delta_p\, w)\vspace{1.5pt} \\
\hspace{-7pt} & = & \hspace{-7pt} (-1)^i\theta_j^i \hspace{-2pt}\left(({\rm id} \hspace{-1pt} \otimes \hspace{-2pt} P[\bar \alpha^!] \hspace{-2pt} \otimes \hspace{-2pt}  M(\bar\alpha))(g^i_{j,x}(w))\right)\hspace{-2pt} (\bar\delta_p^{\hspace{.4pt}\rm o}).
\end{array}\vspace{-2pt}$$
Thus, $(I[\bar\alpha]\otimes P_a(\bar\alpha^!)\otimes {\rm id})(g^{i+1}_{j,a}(w))= (-1)^i({\rm id}\otimes P[\bar \alpha^!]\otimes M(\bar\alpha)) (g^i_{j,x}(w)).$ Now, it is easy to see that \vspace{-2pt} $$(h_j^{i,-i}(a,x,x)\circ g^i_{j,x})(w)+(v_j^{i+1,-i-1}(a,x,a)\circ g^{i+1}_{j,a})(w)=0.$$ This establishes our claim. Thus, Statement 6 holds. 

\vspace{1pt}

To conclude, we shall show that ${\rm Ker}(d^{\hspace{.7pt}0})\subseteq {\rm Im}(\zeta_{\hspace{-.5pt}_M\m})$.
Let $\omega=\sum_{i\in \mathbb{Z}}\, \omega^i\in {\rm Ker}(d^{\hspace{.6pt}0})$, where $\omega^i\in \cG(\cF(M)^i)^{-i}=\oplus_{c\in Q_0}\, I_c\tla\!-i\tra\otimes e_c\La_0^!e_c\otimes M_i(c).$ Note that $\omega^i=0$ for all $i\ge r$. We proceed by induction on the maximal integer
$n_\omega\le r$ such that $\omega^i=0$ for all $i<n_\omega$. If $n_\omega=r$, then $\omega=0\in {\rm Im}(\zeta_{\hspace{-.5pt}_M\m})$. Assume that $n_\omega<r$ and write $n=n_\omega$. Since $d^0(\omega)=0$, we have $v^{n,-n}(\omega^n)=-h^{n-1,1-n}(\omega^{n-1})=0$. By Statement 1, \vspace{-3pt} $${\rm H}^{-n}(\mathfrak{t}^{n}(\cG(\cF(M)^n)^\cdt\hspace{-.5pt})))\cong \oplus_{c\in Q_0}\, {\rm H}^0(\cG(P_a)^\cdt\hspace{-.5pt}) \otimes M_n(c).\vspace{-1.5pt}$$
By Lemma \ref{inj-im}, $\omega^n\in {\rm Soc}(\cG(\cF(N)^n)^{-n})$. Thus, $\omega^n=\sum_{c\in Q_0}\, e^{{\rm o}, \star}_c\otimes e_c\otimes u_c$, where  $u_c\in M_n(c)=e_cM_n$. \vspace{.5pt} Put $u=(-1)^{\frac{(n-1)n}{2}} \sum_{c\in Q_0} u_c\in M_n$. Then, $\zeta_{\hspace{-.5pt}_M\m}(u)\m=\m\sum_{i\in \mathbb{Z}} u^i$, where $u^i\m=\m g^i(u)\m=\m g^i_n(u)\in  G(F(M)^i)^{-i}$. And by Statement 3, we get \vspace{-4pt}
$$u^n = \textstyle \sum_{c, a\in Q_0} (-1)^{\frac{(n-1)n}{2}} g^n_{a,n}(w_c)= \sum_{c, a\in Q_0} e^{{\rm o}, \star}_c \otimes e_a\otimes e_aw_c = \sum_{c\in Q_0} e^{{\rm o}, \star}_c \otimes e_c\otimes w_c.\vspace{-3pt}$$

Put $\nu=\omega-\zeta_{\hspace{-.5pt}_M\m}(u)=\sum_{{i\in \mathbb{Z}}}(\omega^i-u^i) \in {\rm Ker}(d^{\hspace{.7pt}0})$. If $i<n$, then $g^i_n(u)=0$, and hence, $\omega^i-u^i=0$. Since $u^n=\omega^n$, we have $n_\omega< n_\nu$. Thus, $\nu\in {\rm Im}(\zeta_{\hspace{-.5pt}_M\m})$, and hence, $\omega\in {\rm Im}(\zeta_{\hspace{-.5pt}_M\m})$.
The proof of the lemma is completed.

\vspace{2pt}

As a consequence, we obtain the following interesting statement.

\begin{Cor}\label{proj_rel_fd_mod}

Let $\La=kQ/R$ be a Koszul algebra with $Q$ a locally finite quiver. Then every finite dimensional graded $\La$-module admits a graded projective resolution over ${\rm gproj}\La$ and a graded injective co-resolution over ${\rm ginj}\La$.

\end{Cor}

\noindent{\it Proof.} We shall only prove the first part of the statement. Let $M\in {\rm gmod}^{\hspace{.5pt}b\hspace{-3pt}}\La$. By Lemmas \ref{natural-tran} and \ref{Extend_comp_Kfunctor}, we obtain a quasi-isomorphism
$\eta_{\m_M}^\pdt\m\!:\m (\cF^{\hspace{.8pt}C} \ncirc \cG)(M)^\cdt\to M$ such that
$(\cF^{\hspace{.8pt}C} \circ \cG) (M)^n\!=\!\oplus_{(i, x)\in \mathbb{Z}\times Q_0} P_a \langle n - i \rangle \otimes (I_x^!)_{n-i}(a)
\otimes \m M_i(x),$ for all $n\in \Z$.
Since $M$ is finite dimensional, $(\cF^{\hspace{.8pt}C} \ncirc \m \cG) (M)^n\in {\rm gproj}\La$ and
$(\cF^{\hspace{.8pt}C} \ncirc \m \cG) (M)^n=0$ for $n\gg 0.$ That is, $(\cF^C \circ \cG)(M)^\cdt$ is a graded projective resolution of $M$ over ${\rm gproj}\La$. The proof of the corollary is completed.

\vspace{2pt}

More generally, applying the left complex Koszul functor and the right complex Koszul functor yields graded projective resolutions for complexes in $C^{\,\uparrow}_{p,q}(\GrLa)$.

\begin{Lemma}\label{Pres-icres-2}

Let $\La=kQ/R$ be a Koszul algebra with $Q$ a locally finite quiver. Consider
$M^\cdt \in C^{\,\uparrow}_{p,q}(\GrLa)$, where $p,q\in \mathbb{R}$ with $p\ge 1$ and $q\ge 0$. Then, we have a natural quasi-isomorphism $\eta_{M^\cdt}^{\hspace{.4pt}C}\m: (\cF^{\hspace{.5pt}C} \!\circ \cG)^C\m(M^\cdt) \to M^\cdt\hspace{-2pt}$.

\end{Lemma}

\noindent{\it Proof.} Consider the functor $\cF^{\hspace{.5pt}C} \! \circ\, \cG: {\rm Mod}\hspace{.5pt}\La \rightarrow C({\rm Mod}\hspace{.5pt}\La)$ and the embedding functor $\kappa: {\rm Mod}\hspace{.5pt}\La \rightarrow C({\rm Mod}\hspace{.5pt}\La)$. In view of Lemma \ref{natural-tran}, we obtain a functorial morphism $\eta=(\eta_{\hspace{-.4pt} M}^\cdt\hspace{-2.5pt})_{M\in {\rm Mod}\hspace{.5pt}\mathit\Lambda}: \cF^{\hspace{.5pt}C} \m \circ \hspace{.4pt} \cG \to \kappa$. \vspace{1.5pt} By Lemma \ref{main-lemma1}, it extends to a functorial morphism $\eta^{\hspace{.5pt}C}\hspace{-2pt}: (\cF^{\hspace{.5pt}C} \!\!\circ \m \cG)^C\to \kappa^C={\rm id}_{C({\rm Mod}\it\Lambda)}$ \vspace{-.4pt}
such that $\eta^{\hspace{.5pt}C}_{M^\cdt}\!=\!\T(\eta_{\m M^\cdt}^\cdt\hspace{-3.5pt}):
(\cF^{\hspace{.5pt}C} \hspace{-3pt} \circ \hspace{-1.6pt} \cG)^C(M^\cdt)\to M^\cdt,$ where $\eta_{\m M^\cdt}^\ydt=(\eta^j_{\hspace{-.8pt} M^i}\m)_{i,j\in \mathbb Z}: (\cF^{\hspace{.5pt}C} \!\circ \cG)(M^\cdt\m)^\cdt \to \kappa(M^\cdt\m)^\cdt$. We claim that $\eta^{\hspace{.5pt}C}_{M^\cdt}\!=\!\T(\eta_{\m M^\cdt}^\cdt\hspace{-3.5pt})$ is a quasi-isomorphism, or equivalently, $\T(\eta_{\m M^\cdt}^\cdt\hspace{-3.5pt})_s$ is a quasi-isomorphism for any $s\in \Z$.

Now $\T(\eta_{\m M^\cdt}^\ydt\hspace{-1.5pt})_s  \!=\! \T((\eta_{\hspace{-1pt} M^\cdt}^\cdt\hspace{-4pt})_s)$, where \vspace{1pt} $(\eta_{\m M^\cdt}^\cdt\hspace{-4pt})_s \!=\! ((\eta_{\m M^i}^j\hspace{-1.5pt})_s)_{i,j\in \mathbb Z}: (\cF^{\hspace{.5pt}C} \! \circ \cG)(M^\cdt\hspace{-1pt})^\ydt_s \!\to\! \kappa(M^\cdt \hspace{-1pt})^\ydt_s$.
For any $i\in \Z$, by Lemma \ref{natural-tran}, $(\eta_{_{\hspace{-.8pt}M^i}}^\cdt\hspace{-3.5pt})_s\hspace{-1pt}: \mk{t}^i((\cF^{\hspace{.5pt}C} \hspace{-2pt} \circ \cG)(M^i)^\ydt)_s \!\to \mk{t}^i(\kappa(M^i)^\cdt\hspace{-.5pt})_s$ is a quasi-isomorphism. On the other hand, $\kappa(M^\cdt\m)_s^\cdt$ is clearly diagonally bounded-above. And given $n\in \Z$, by Lemma \ref{Comp_Kfunctors}(1), the $n$-diagonal of $(\cF^{\hspace{.5pt}C} \!\m\circ\m \cG)(M^\cdt\hspace{-1.2pt})^\ydt$ consists of \vspace{-3.5pt}
$$(\cF^{\hspace{.5pt}C} \!\m\circ\m \cG)(M^i)^{n-i} = \oplus_{j\in \mathbb{Z}; \hspace{.7pt} a, x\in Q_0} P_a\tla n\!-\!i\!-\!j\tra \otimes (I_x^!)_{n-i}(a)
\otimes M^i_j(x); \; i\in\Z
\vspace{-3.5pt}$$
As a consequence, the $n$-diagonal of $(\cF^{\hspace{.5pt}C} \!\!\circ\m \cG)(M^\cdt\m)^\cdt\!\!_s$ consists of \vspace{-3pt}
$$(\cF^{\hspace{.5pt}C} \!\!\circ\m \cG)(M^i)^{n-i}_s=\oplus_{j\in \mathbb{Z}; \, a, x \in Q_0}\, \La_{n+s-i-j}e_a \otimes (I_x^!)_{n-i}(a)
\otimes M_j^i(x); \; i\in \Z.
\vspace{-3pt}
$$

Let $t\in \Z$ be such that $M^i_j=0$ for $i-qj>t$. Fix any $i> (q(n+s)+t)(q+1)^{-1}$. If $j\hspace{-1pt} >\hspace{-1pt}  n+s-i$, then $\La_{n+s-i-j}=0$; otherwise, $M^i_j=0$ since $i-qj\hspace{-1pt} \ge\hspace{-1pt}  i-q(n+s-i)\hspace{-1pt}>\hspace{-1pt} t.$ So $(\cF^{\hspace{.5pt}C} \! \circ \cG)(M^\cdt)^\ydt_s$ is $n$-diagonally bounded-above. Hence, $(\cF^{\hspace{.5pt}C} \!\circ  \cG)(M^\cdt)^\ydt_s$ is diagonally bounded-above.
By Lemma \ref{V-cone}, $\T((\eta_{\m M^\cdt}^\cdt\hspace{-4pt})_s\m)$ is a quasi-isomorphism. This proves our claim. The proof of the lemma is completed.

\vspace{3pt}

Applying the right complex Koszul functor and the left complex Koszul functor yields injective coresolutions for complexes of bounded-above graded modules.

\begin{Lemma}\label{Pres-icres}

Let $\La=kQ/R$ be a Koszul algebra with $Q$ a locally finite quiver. For
$M^\cdt \! \in C({\rm GMod}^{-\hspace{-3.5pt}}\La),$ we have a natural quasi-isomorphism $\zeta^C_{M^\cdt}\m: M^\ydt \hspace{-1pt}\to \hspace{-1pt} (\cG^{\hspace{.3pt}C} \m \circ \cF)^C\m(M^\cdt)$.
\end{Lemma}

\noindent{\it Proof.} Consider the embedding functor $\kappa: {\rm GMod}^-\hspace{-3.5pt}\La \rightarrow C({\rm GMod}\La)$ and the functor
$\cG^{\hspace{.3pt}C} \!\m \circ \m \cF: {\rm GMod}^-\hspace{-3pt}\La \!\to\! C({\rm GMod}\La)$. In view of Lemma \ref{natural-tran-2}, we have a functorial morphism \vspace{.5pt} $\zeta \! = \! (\zeta_{\m M}^\ydt\hspace{-1pt})_{M\in {\rm Mod}^-\hspace{-2pt}\it\Lambda}\hspace{-1pt}: \hspace{-1pt} \kappa \to \cG^{\hspace{.3pt}C} \!\circ\! \cF$ which, by Lemma \ref{main-lemma1}, \vspace{1pt} extends to a functorial morphism
$\zeta^C\hspace{-1pt}: {\rm id}_{{\rm GMod}^-\hspace{-2.5pt}\mathit\Lambda} = \kappa^C \!\to\! (\cG^{\hspace{.3pt}C}\!\!\circ\m \cF)^C$.

\vspace{1pt}

Let $M^\cdt \!\in C({\rm GMod}^-\hspace{-2.5pt}\La)$. Then, $(\cG^{\hspace{.3pt}C} \!\m \circ\m \cF)^C\!(M^\cdt\hspace{-.8pt})=\T((\cG^{\hspace{.3pt}C} \!\m\circ\m \cF)(M^\cdt\hspace{-1pt})^\cdt)$. By Lemma \ref{main-lemma1}, $\zeta_{M^\cdt}^C \! =\! \T(\zeta_{\hspace{-.4pt} M^\cdt}^\cdt\hspace{-4pt})\!:\! M^\cdt \!\to\! (\cG^{\hspace{.3pt}C} \!\m\circ\m \cF)^C(M^\cdt\m)$ \vspace{.5pt} with $\zeta_{\m M^\cdt}^\cdt \!=\!(\zeta_{\m M^i}^j\m)_{i,j\in \mathbb Z}\m:\m \kappa(M^\cdt)^\cdt \!\to\m (\cG^{\hspace{.3pt}C} \!\m\circ \m \cF)(M^\cdt)^\cdt$. It is evident that $\kappa(M^\cdt)^\cdt$ is diagonally bounded-above. Given $n\in \Z$, by Lemma \ref{Comp_Kfunctors}(2), we deduce that the $n$-diagonal of $(\cG^{\hspace{.3pt}C} \!\m\circ\m \cF)(M^\cdt)^\cdt$ consists of \vspace{-2pt} $$(\cG^{\hspace{.3pt}C} \!\m\circ\m \cF)(M^i)^{n-i}=\oplus_{j\in \Z; \hspace{.7pt} x, a\in Q_0} I_x\tla n\!-\!i\!-\!j\tra\otimes e_x \La_{n-i}^! e_a \otimes M^i_j(a); \;\; i\in\Z.\vspace{-4pt}$$
Since $\La^!_{n-i}=0$ for $i>n$, we see that $(\cG^{\hspace{.3pt}C} \!\!\circ\m \cF)(M^\cdt\m)^\cdt$ is $n$-diagonally bounded-above. So $(\cG^{\hspace{.3pt}C} \!\!\circ\m \cF)(M^\cdt\m)^\cdt$ is diagonally bounded above. Further for any $i\in \Z$, by Lemma \ref{natural-tran-2}, $\zeta_{M^i}^\cdt\m:\m \mk{t}^i(\kappa(M^i)^\cdt) \to \mk{t}^i((\cG^{\hspace{.3pt}C} \!\m \circ \m \cF)(M^i)^\cdt)$ is a quasi-isomorphism. \vspace{.5pt}
Thus, by Lemma \ref{V-cone}, $\zeta_{M^\ydt}^C$ is a quasi-isomorphism. The proof of the lemma is completed.

\vspace{2pt}

We are ready to prove the main result of this section, which includes the classical Koszul duality of Belinson, Ginzburg and Soergel; see \cite{BGS}.

\begin{Theo}\label{Main}

Let $\La=kQ/R$ be a Koszul algebra, where $Q$ is a locally finite quiver. Consider $p, q\in \mathbb R$ with $p\ge 1$ and $q\ge 0$.

\begin{enumerate}[$(1)$]

\vspace{-1pt}

\item The right derived Koszul functor \vspace{.0pt} $\mathcal{F}^D_{p,q}\!\m:\! D^{\,\downarrow}_{p,q}(\GrLa) \!\!\to\!\m  D^{\,\uparrow}_{q+1, p-1}(\GrLa^!)$ and the left derived Koszul functor $\mathcal{G}^D_{q+1\m,\hspace{.4pt}p-1}\!\m:\m D^{\,\uparrow}_{q+1,p-1}(\GrLa^!) \!\to\! D^{\,\downarrow}_{p,q}(\GrLa)$ are mutually quasi-inverse.


\item The left derived Koszul functor  \vspace{.5pt} $\cG^D_{p,q}\!:\! D^{\,\uparrow}_{p,q}(\GrLa) \!\m\to\! D^{\,\downarrow}_{q+1, p-1}(\GrLa^!)$ and the right derived Koszul functor $\cF^D_{q+1, p-1}\!:\! D^{\,\downarrow}_{q+1,p-1}(\GrLa^!) \!\m\to\! D^{\,\uparrow}_{p,q}(\GrLa)$ are mutually quasi-inverse.


\item In the above two statements, $\GrLa$ and $\GrLa^!$ can be simulta\-neously replaced by ${\rm gmod}\La$ and ${\rm gmod}\hspace{.4pt}\La^!\hspace{-2pt},$ respectively.

\end{enumerate}

\end{Theo}

\noindent{\it Proof.} Recall that $C^{\,\downarrow}_{p,q}(\GrLa)\subseteq C({\rm GMod}^-\!\!\La)$. By Lemmas \ref{Pres-icres} and \ref{Pres-icres-2}, we have natu\-ral quasi-isomorphisms
$\zeta_{N^\ydt}^C\!:\! N^\cdt \!\to\! (\cG^{\hspace{.3pt}C} \!\!\circ\m \cF)^C\m(N^\cdt)$ \vspace{-1pt} for $N^\cdt \!\in\! C^{\,\downarrow}_{p,q}(\GrLa)$, and $\eta^C_{M^\cdt}\!:\! (\mathcal{F}^{\hspace{.5pt}C} \!\!\circ\m \mathcal{G})^C(M^\cdt) \!\to\! M^\cdt$ for $M^\cdt \!\in\m C^{\,\uparrow}_{q+1, p-1}(\GrLa^!)$. \vspace{.5pt} And by Proposition \ref{F-composition},
$(\mathcal{G}^{\hspace{.5pt}C} \circ \mathcal{F})^{\hspace{.5pt}C} \m(N^\cdt\hspace{-1pt}) \m =\m (\cG^{\hspace{.3pt}C}_{q+1,p-1} \circ \cF^{\hspace{.5pt}C}_{p,q})(N^\cdt\m)$ \vspace{.5pt} and
$(\cF^{\hspace{.5pt}C} \circ \cG)^C\!(M^\cdt\m) \m =\m (\cF^{\hspace{.5pt}C}_{p,q} \circ \cG^{\hspace{.3pt}C}_{q+1,p-1})(M^\cdt\m)$. This gives rise to natural isomorphisms \vspace{-.5pt} $\zeta^D_{N^\cdt}\!:\! N^\cdt \!\to\! (\cG^D_{q+1,p-\m1} \! \circ\m \cF^D_{p,q})(N^\cdt\hspace{-.8pt})$ for $N^\cdt$ in $D^{\,\downarrow}_{p,q}(\GrLa)$ and
$\eta^D_{M^\cdt}\!:\! (\cF^D_{p,q} \m \circ \cG^D_{q+1,p-1})(M^\cdt\m) \!\to\! M^\cdt$ for $M^\cdt$ in $\m D^{\,\uparrow}_{q+1,p-1}(\GrLa^!)$. So, $\cF^D_{p,q}$ and $\hspace{.3pt}\cG^D_{q+1,p-1}$ are mutually quasi-inverse. This establishes Statement (1). Similarly, we can prove Statement (2). Finally using the same argument, we deduce Statement (3) from Theorem \ref{F-diag}(3). The proof of the theorem is completed.

\vspace{3pt}

\noindent{\sc Remark.} Taking $p=1$ and $q=0$, Theorem \ref{Main}(1) and (2) have been established in \cite[(2.12.1)]{BGS} with a sophisticated proof, while Theorem \ref{Main}(3) has been proved for a positively graded Koszul category in \cite[Theorem 30]{MOS}.

\vspace{3pt}

We shall show that the bounded derived Koszul functors are also triangle equivalences in the locally bounded case,

\vspace{-.5pt}

\begin{Theo}\label{Main-2}

Let $\La=kQ/R$ be a Koszul algebra with $Q$ a locally finite quiver.

\begin{enumerate}[$(1)$]

\vspace{-1.5pt}

\item If $\hspace{-.5pt}\La\hspace{-.5pt}$ is locally right bounded and $\hspace{-.5pt}\La^!\hspace{-.5pt}$ is locally left bounded, then the derived Koszul functors $\cF^D_b\hspace{-2.5pt}:\hspace{-2pt} D({\rm GMod}^b\hspace{-3pt}\La) \hspace{-2pt}\to\hspace{-2pt} D({\rm GMod}^b\hspace{-3pt}\La^!)$ and $\cG^{\hspace{.3pt}D}_b\hspace{-2.5pt}:\hspace{-2pt} D({\rm GMod}^b\hspace{-3pt}\La^!) \hspace{-2pt}\to\hspace{-2pt} D({\rm GMod}^b\hspace{-3pt}\La)$ are mutually quasi-inverse.

\item If $\hspace{-.5pt}\La\hspace{-.5pt}$ is locally left bounded and $\hspace{-.5pt}\La^!\hspace{-.5pt}$ is locally right bounded, then the derived Koszul functors $\cG^{\hspace{.3pt}D}_b\hspace{-2.5pt}:\hspace{-2pt} D({\rm GMod}^b\hspace{-3pt}\La) \hspace{-2pt}\to\hspace{-2pt} D({\rm GMod}^b\hspace{-3pt}\La^!)$ and
$\cF^D_b\hspace{-2.5pt}:\hspace{-2pt} D({\rm GMod}^b\hspace{-3pt}\La^!) \hspace{-2pt}\to\hspace{-2pt} D({\rm GMod}^b\hspace{-3pt}\La)$ are mutually quasi-inverse.

\item In the above two statements, ${\rm GMod}^b\hspace{-3pt}\La$ and ${\rm GMod}^b\hspace{-3pt}\La^!$ can be simulta\-neously replaced by ${\rm gmod}^b\hspace{-3pt}\La$ and ${\rm gmod}^b\hspace{-3pt}\La^!\hspace{-2pt},$ respectively.

\end{enumerate}\end{Theo}

\noindent{\it Proof.} Suppose that $\La$ is right locally bounded and $\La^!$ is left locally bounded. Since $\La^!$ is Koszul with $(\La^!)^!=\La$,
by Theorem \ref{F-diag-fd-1}, we obtain triangle exact functors $\cF^D_b\!\!:\m D^b({\rm GMod}^{\hspace{.5pt}b}\!\m\La) \!\to\! D^b({\rm GMod}^{\hspace{.5pt}b}\!\m\La^!)$ and $\cG^D_b\!\!:\m D^b({\rm GMod}^{\hspace{.5pt}b}\!\m\La^!) \!\to\!\m D^b({\rm GMod}^{\hspace{.8pt}b}\!\m\La)$. As argued in the proof of Theorem \ref{Main}, we see that they are mutual quasi-inverse. This establishes Statement (1). Similarly, one can verify Statements (2) and (3). The proof of the theorem is completed.

\vspace{2pt}

\noindent{\sc Remark.} In case $\La$ is of finite dimensional and $\La^!$ is left noetherian, Beilinson, Ginzburg and Soergel established Theorem \ref{Main-2}(3); see \cite[(2.12.6)]{BGS}.


\vspace{-4pt}

\section{Graded almost split triangles}

The objective of this section is to study almost split triangles in derived categories of graded modules over a Koszul algebra. We shall show that the Auslander-Reiten translations and the Serre functors are related to derived Koszul functors.

\vspace{1pt}

Throughout this section let $\La=kQ/R$ be a Koszul algebra, where $Q$ is a locally finite quiver. It is known that ${\rm gproj}\La$ and ${\rm ginj}\La$ are Hom-finite and Krull-Schmidt, and so are $K^{\hspace{.5pt}b\hspace{-.5pt}}({\rm gproj}\La)$ and $K^{\hspace{.5pt}b\hspace{-.5pt}}({\rm ginj}\La)$; \vspace{.5pt} see \cite[(2.12.2), (4.1.1)]{LLi}.
Moreover, there exists a Nakayama functor $\nu: {\rm gproj}\La\to \GrLa^!$, which restricts to an equivalence $\nu: {\rm gproj}\La\to {\rm ginj}\La$ such that $\nu(P_x\sla s\sra \otimes V) \cong  I_x\sla s \sra \otimes V,$ for $(s,x)\in \Z\times Q_0$ and $V \!\m \in \! {\rm mod}k\hspace{.2pt}$; see \cite[(3.2.1)]{LLi}. Applying $\nu$ componentwise, we obtain a triagnle-equivalence $\nu: K^{\hspace{.5pt}b\hspace{-.5pt}}({\rm gproj}\La) \to K^{\hspace{.5pt}b\hspace{-.5pt}}({\rm ginj}\La);$ see \cite[(4.1.1)]{LLi}. 

\vspace{2pt}

First, by making use of the derived Koszul functors, we may describe some almost split triangles in $D^{\hspace{.5pt}b\hspace{-.5pt}}({\rm gmod}\La)$ in terms of bounded complexes of finite dimensional graded $\La^!$-modules ; compare \cite[(5.2)]{BaL2}.


\begin{Prop}\label{Main-ast-1}

Let $\La=kQ/R$ be a Koszul algebra with $Q$ a locally finite quiver. If
$M^\cdt \hspace{-1pt} \in \hspace{-1pt} D^{\hspace{.6pt}b\hspace{-.5pt}}({\rm gmod}^{\hspace{.5pt}b\hspace{-2.8pt}}\La^!)$ with $\cF^D\hspace{-1pt}(M^\cdt)$ or $\cG^D\hspace{-1pt}(M^\cdt)$ indecomposable, \vspace{-1.5pt} then we have an almost split triangle $\xymatrixcolsep{16pt}\xymatrixrowsep{14pt}\xymatrix{\cG^D\hspace{-1pt}(M^\cdt)[-1] \ar[r] & N^\cdt \ar[r] & \cF^D\hspace{-1pt}(M^\cdt)\ar[r]
& \cG^D\hspace{-1pt}(M^\cdt)\hspace{.5pt}}\hspace{-2pt}\vspace{-3pt}$ in $D^{\hspace{.6pt}b\hspace{-.5pt}}({\rm gmod}\La). \vspace{-0.5pt}$

\end{Prop}

\noindent{\it Proof.} Consider a complex $M^\cdt\in C^{\hspace{.5pt}b\hspace{-.5pt}}({\rm gmod}^{\hspace{.5pt}b\hspace{-3pt}} \La^!)$. By Theorem \ref{F-diag-fd} and Lemma \ref{cplx-kosz-fun}, $\mathcal{F}^{\hspace{.6pt}C\hspace{-.5pt}}(M^\cdt)\in C^{\hspace{.6pt}b\hspace{-.5pt}}({\rm gproj}\La^!)$
such that $\cF^{\hspace{.5pt}C\hspace{-.5pt}}(M^\cdt)^n=\oplus_{(i, x)\in \Z\times Q_0} P^{\hspace{.4pt}!}_x\tla n\!-\!i \tra\otimes M^i_{n-i}(x),$ for all $n\in \Z$. Then, \vspace{.5pt} $\nu \cF^{\hspace{.6pt}C\hspace{-1pt}}(M^\cdt)^n\cong \oplus_{(i, x)\in \Z\times Q_0} I^{\hspace{.4pt}!}_x\nla n\!-\!i \tra\otimes M^i_{n-i}(x) = \cG^{\hspace{.6pt}C\hspace{-1pt}}(M^\cdt)^n$, for all $n\in \Z$. That is, $\cG^{\hspace{.6pt}C\hspace{-1pt}}(M^\cdt)\cong \nu \cF^{\hspace{.6pt}C\hspace{-1pt}}(M^\cdt)$. In particular, $\cF^{\hspace{.6pt}K\hspace{-1pt}}(M^\cdt)$ is indecomposable in $K^b({\rm gproj}\La)$ if and only if $\cG^{\hspace{.6pt}K\hspace{-1pt}}(M^\cdt)$ is indecomposable in $K^b({\rm ginj}\La)$. If $\cF^{\hspace{.6pt}D\hspace{-1pt}}(M^\cdt)$ or $\cG^{\hspace{.6pt}D\hspace{-1pt}}(M^\cdt)$ is indecomposable in $D^{\hspace{.6pt}b\hspace{-.5pt}}({\rm gmod}\La)$, then $\cF^{\hspace{.6pt}K\hspace{-1.5pt}}(M^\cdt)$ is indecomposable in $K^b({\rm gproj}\La)$. And in this case, by Theorem 4.1.2 in \cite{LLi}, we obtain a desired almost split triangle in $D^{\hspace{.6pt}b\hspace{-.5pt}}({\rm gmod}\La)$. The proof of the proposition is completed.

\vspace{3pt}

\noindent{\sc Remark.} Our terminology of left and right Koszul functors is explained by the almost split triangle stated in Proposition \ref{Main-ast-1}.

\vspace{3pt}

\noindent{\sc Example.} Let $\La=kQ/R$ be a Koszul algebra with $a\in Q_0$.  It is well known that $S_a$ is indecomposable in $D^{\hspace{.6pt}b\hspace{-.5pt}}({\rm gmod}\La)$; see \cite[(III.3.4.7)]{Mil}. If $I^!_a$ or $P_a^!$ is finite dimensional, by Lemma \ref{inj-im}, $\mathcal{F}^{\hspace{.6pt}D\hspace{-1pt}}(I^!_a)\cong  S_a\vspace{.3pt}$ or $\mathcal{G}^{\hspace{.6pt}D\hspace{-1pt}}(P_a^!)\cong S_a\vspace{.3pt}$ in $D^{\hspace{.6pt}b\hspace{-.5pt}}({\rm gmod}\La)$ respectively, and by Proposition \ref{Main-ast-1}, there exists an almost split triangle in $D^{\hspace{.6pt}b\hspace{-.5pt}}({\rm gmod}\La)$ ending or starting with $S_a$ respectively.

\vspace{2pt}

In case $\La^!$ is locally bounded, we shall establish the existence of almost split triangles in $D^{\hspace{.6pt}b\hspace{-.5pt}}({\rm gmod}\La)$ for bounded complexes of finite dimensional graded $\La$-modules
and describe the Auslander-Reiten translates in terms of derived Koszul functors. For this purpose, we call a bounded complex over ${\rm gmod}\La$ {\it derived-indecomposable} if it is indecomposable in $D^{\hspace{.6pt}b\hspace{-.5pt}}({\rm gmod}\La)$.

\begin{Theo}\label{Main-ast-2}

Let $\La=kQ/R$ be a Koszul algebra with $Q$ a locally finite quiver.

\begin{enumerate}[$(1)$]

\vspace{-2pt}

\item Every derived-indecomposable complex $M^\cdt$ in $C^{\hspace{.6pt}b\hspace{-.5pt}}({\rm gmod}^{\hspace{.3pt}b\hspace{-3pt}}\La)$ is the ending term of an almost split triangle in $D^{\hspace{.6pt}b\hspace{-.5pt}}({\rm gmod}\La)$ if and only if \vspace{.5pt} $\La^!$ is locally right bounded$\,;$ and in this case, $\tau M^\cdt\cong \cG^{\hspace{-.3pt}D\hspace{-.6pt}}(\cG^{\hspace{-.3pt}D\hspace{-.8pt}}(M^\cdt)\m)[-1]$.

\item Every derived-indecomposable complex $M^\cdt$ in $C^{\hspace{.6pt}b\hspace{-.5pt}}({\rm gmod}^{\hspace{.3pt}b\hspace{-3pt}}\La)$ is the starting term of an almost split triangle in $D^{\hspace{.6pt}b\hspace{-.5pt}}({\rm gmod}\La)$ if and only if \vspace{.5pt}  $\La^!$ is locally left bounded$\,;$ and in this case, $\tau^-\m M^\cdt\cong \cF^{\hspace{.0pt}D\hspace{-.8pt}}(\cF^{\hspace{0pt}D\hspace{-.8pt}}(M^\cdt\hspace{-.8pt}))[1]$.


\end{enumerate}\end{Theo}

\noindent{\it Proof.} We shall only prove Statement (1). Given $a\in Q_0$, by Lemma \ref{k-cplx-iso}, $S_a$ has a minimal graded projective resolution $\mathcal{P}_a^\pdt$ with $\mathcal{P}_a^{-n}=\oplus_{x\in Q_0} P_x\mla \m-n \nra \otimes D(e_a \La^!_n e_x),$ for all $n\in \Z$. If $S_a$ is the ending term of
an almost split triangle in $D^{\hspace{.6pt}b\hspace{-.5pt}}({\rm gmod}\La)$, then $\mathcal{P}_a^\pdt$ is finite; see \cite[(5.2)]{LiN}. Then, $e_a \La^!_n=0$ for all but finitely many $n \ge 0$. So, $e_a\La^!$ is finite dimensional. This establishes the necessity of Statement (1).

\vspace{.5pt}

Suppose that $\La^!$ is locally right bounded. Then ${\rm ginj}\La^! \!\subseteq\! {\rm gmod}^{\hspace{.6pt}b\hspace{-2.8pt}}\La^!$. Consider a derived-indecomposable complex $M^\cdt\in C^{\hspace{.6pt}b\hspace{-.5pt}}({\rm gmod}^{\hspace{.3pt}b\hspace{-3pt}}\La).$ In view of Theorem \ref{F-diag-fd-1}(3), \vspace{.5pt} $\cG^{\hspace{.3pt} C\hspace{-1pt}}(M^\cdt)\!\in\! C^{\hspace{.6pt}b\hspace{-.8pt}}({\rm ginj}\La^!)$, and by Theorem \ref{F-diag-fd}(3),
$\cF^{\hspace{.5pt}C\hspace{-1pt}}(\cG^{\hspace{.5pt}C\hspace{-1pt}}(M^\cdt\m)\m) \in C^{\hspace{.6pt}b\hspace{-.8pt}}({\rm gproj}\La)$. Thus, by Lemma \ref{Pres-icres-2}(2) and Proposition \ref{F-composition}, $\cF^{\hspace{.3pt}D\hspace{-1pt}}(\cG^{\hspace{.3pt}D\hspace{-1pt}}(M^\cdt\hspace{-.3pt})\hspace{-.3pt})\!\cong \!\m M^\cdt$ in $D^{\hspace{.6pt}b\hspace{-.5pt}}({\rm gmod}\La)$. Observing that $\cG^{\hspace{.3pt} D\hspace{-1pt}}(M^\cdt)\!\in\! C^{\hspace{.6pt}b\hspace{-.8pt}}({\rm gmod}^{\hspace{.6pt}b\hspace{-2.8pt}}\La^!)$, we deduce from Proposition \ref{Main-ast-1}  \vspace{-2pt} an almost split triangle $\xymatrixrowsep{14pt}\xymatrixcolsep{18pt}\xymatrix{
\hspace{-2pt} \cG^{\hspace{-.3pt}D\hspace{-.8pt}}(\cG^{\hspace{-.3pt}D\m}(M^\cdt\m)\hspace{-.2pt})[-1] \ar[r] & N^\cdt \ar[r] & M^\cdt \ar[r] & \cG^{\hspace{-.3pt}D\hspace{-.8pt}}(\m \cG^{\hspace{-.3pt}D\!}(M^\cdt))} \vspace{-1pt} $
in $D^b\m({\rm gmod}\La).$ The proof of the theorem is completed.

\vspace{2pt}

\noindent{\sc Example.} Let $\La=kQ$, where $Q$ is a locally finite quiver. Then $\La^!=kQ^{\rm o}/R^!$, where $R^{\hspace{.4pt}!}$ is the ideal in $kQ^{\rm o}$ generated by the paths of length two. Clearly, $\La^!$ is locally bounded. By Theorem \ref{Main-ast-2}, every derived-indecomposable bounded complex over ${\rm gmod}^{\hspace{.3pt}b\hspace{-3pt}}\La$ is the starting term, as well as the ending term, of an almost split triangle in $D^{\hspace{.3pt}b\hspace{-.5pt}}({\rm gmod}\La)$.

\vspace{3pt}

To conclude, we concentrate on the bounded derived category $D^{\hspace{.5pt}b\hspace{-.5pt}}({\rm gmod}^b\hspace{-2.8pt}\La)$.

\begin{Prop}\label{KS}

Let $\La=kQ/R$ be a Koszul algebra with $Q$ a locally finite quiver.
Then $D^{\hspace{.5pt}b\hspace{-.5pt}}({\rm gmod}^b\hspace{-2.8pt}\La)$ is Hom-finite and Krull-Schmidt. 

\end{Prop}

\noindent{\it Proof.} Let $M\in {\rm gmod}^{\hspace{.5pt}b\hspace{-3pt}}\La$. By Corollary \ref{proj_rel_fd_mod},
$M$ admits a graded projective resolution over ${\rm gproj}\La$. Given any $ N \in {\rm gmod}^{\hspace{.5pt}b\hspace{-3pt}}\La$, since $\GHom_{\mathit\Lambda}(P, N)$ is finite dimensional for $P\in {\rm gproj}\La$; see \cite[(2.12.1)]{LLi},
${\rm GExt}_{\mathit\Lambda}^n(M, N)$ is finite dimensional for all $n\in \Z$. Thus, $D^{\hspace{.5pt}b\hspace{-.5pt}}({\rm gmod}^{\hspace{.5pt} b\hspace{-2.8pt}}\La)$ is Hom-finite and Krull-Schmidt; see \cite[Corollary B]{LeC}.
The proof of the lemma is completed.

\vspace{2pt}

Finally, we find conditions for $D^{\hspace{.5pt}b\hspace{-.5pt}}({\rm gmod}^{\hspace{.5pt} b\hspace{-2.8pt}}\La)$ to have almost split triangles and describe the Serre functors in terms of the derived Koszul functors.

\vspace{-1pt}


\begin{Theo}\label{Main-ast-3}

Let $\La=kQ/R$ be a locally bounded Koszul algebra, where $Q$ is a locally finite quiver.

\vspace{-2pt}

\begin{enumerate}[$(1)$]

\item There exist almost split triangles in $D^{\hspace{.5pt}b\hspace{-.5pt}}({\rm gmod}^{\hspace{.5pt} b\hspace{-2.8pt}}\La)$ on the right if and only if $\La^!$ is right locally bounded$\,;$ and in this case, $\cG^D \circ \cG^D: D^{\hspace{.5pt}b\hspace{-.5pt}}({\rm gmod}^{\hspace{.5pt} b\hspace{-2.8pt}}\La) \to D^{\hspace{.5pt}b\hspace{-.5pt}}({\rm gmod}^{\hspace{.5pt} b\hspace{-2.8pt}}\La)$ is a right Serre functor.


\item There exist almost split triangles in $D^{\hspace{.5pt}b\hspace{-.5pt}}({\rm gmod}^{\hspace{.5pt} b\hspace{-2.8pt}}\La)$ on the left if and only if $\La^!$ is left locally bounded$\,;$ and in this case, $\cF^D \circ \cF^D: D^{\hspace{.5pt}b\hspace{-.5pt}}({\rm gmod}^{\hspace{.5pt} b\hspace{-2.8pt}}\La) \to D^{\hspace{.5pt}b\hspace{-.5pt}}({\rm gmod}^{\hspace{.5pt} b\hspace{-2.8pt}}\La)$ is a left Serre functor.


\item There exist almost split triangles in $D^{\hspace{.5pt}b\hspace{-.5pt}}({\rm gmod}^{\hspace{.5pt} b\hspace{-2.8pt}}\La)$ if and only if $\La^!$ is locally bounded. In this case, $\cG^D \ncirc \cG^D\hspace{-3pt}:\hspace{-1pt} D^{\hspace{.5pt}b\hspace{-.5pt}}({\rm gmod}^{\hspace{.5pt} b\hspace{-2.8pt}}\La) \!\to\! D^{\hspace{.5pt}b\hspace{-.5pt}}({\rm gmod}^{\hspace{.5pt} b\hspace{-2.8pt}}\La)$ is a right Serre equivalence and
$\cF^D \hspace{-1.5pt} \circ  \cF^D \hspace{-3pt}:\hspace{-1pt} D^{\hspace{.5pt}b\hspace{-.5pt}}({\rm gmod}^{\hspace{.5pt} b\hspace{-2.8pt}}\La) \!\to\! D^{\hspace{.5pt}b\hspace{-.5pt}}({\rm gmod}^{\hspace{.5pt} b\hspace{-2.8pt}}\La)$ is a left Serre equivalence.

\end{enumerate}\end{Theo}

\noindent{\it Proof.} Since $\La$ is locally bounded, ${\rm gproj}\La$ and ${\rm ginj}\La$ are subcategories of ${\rm gmod}^{\hspace{.4pt}b}\hspace{-3pt}\La$. Now, using the same argument for proving Theorem \ref{Main-ast-2}(1), we can establish Statement (1). And Statement (2) can be shown in a similar fashion. Finally, Statement (3) is an immediate consequence of the first two statements.
The proof of the theorem is completed.

\vspace{2pt}

\noindent{\sc Remark.} Let $\La=kQ/(kQ^+)^2$, where $Q$ is a locally finite quiver without infinite paths. Then, $\La$ and $\La^!=kQ$ are locally bounded. Hence, $D^{\hspace{.5pt}b\hspace{-.5pt}}({\rm gmod}^{\hspace{.4pt}b\hspace{-2.8pt}}\La)$ and $D^{\hspace{.5pt}b\hspace{-.5pt}}({\rm gmod}^{\hspace{.4pt}b\hspace{-1pt}} kQ)$ are equivalent and have almost split triangles; see (\ref{Main-2}) and (\ref{Main-ast-3}). One could describe the Aulsander-Reiten components of ${\rm gmod}^{\hspace{.4pt}b}\m kQ$, and this would yield a description of the Auslander-Reiten components for $D^b\m({\rm gmod}^{\hspace{.2pt}b}\m kQ)$ and $D^{\hspace{.5pt} b \hspace{-.5pt}}({\rm gmod}^{\hspace{.5pt}b\hspace{-2.8pt}}\La)$, as is done in the ungraded case; see \cite{BaL2, BLP}.

\vspace{3pt}

{\sc Acknowledgment.} The last named author is supported in part by the Natural Sciences and Engineering Research Council of Canada.

%


\vspace{-1pt}

\begin{thebibliography}{99}

\vspace{1pt}


%

%
%
\bibitem{As2} {\sc H. Asashiba}, ``A generalization of Gabriel's Galois covering functors and derived equivalences," J. Algebra 334 (2011) 109-149.

%
%

%



%

%


\bibitem{BaL} {\sc R. Bautista and S. Liu}, ``Covering theory for linear categories with application to derived categories," J. Algebra 406 (2014) 173-225.


\bibitem{BaL2} {\sc R. Bautista and S. Liu}, ``The bounded derived categories of an algebra with radical squared zero," J. Algebra 482 (2017) 303-345.


\bibitem{BLP} {\sc R. Bautista, S. Liu and C. Paquette}, ``Representation theory of strongly locally finite quivers," Proc. London Math. Soc. 106 (2013) 97-162.

%




\bibitem{BGS} {\sc A. Beilinson, V. Ginzburg and W. Soergel}, ``Koszul duality patterns in representation theory," J. Amer. Math. Soc. 9 (1996) 473-527.

%

%


\bibitem{BoG} {\sc K. Bongartz and P. Gabriel}, ``Covering spaces in representation theory," Invent. Math.  65 (1982) 331-378.

%




\bibitem{CEi} {\sc H. Cartan and S. Eilenberg}, ``Homological Algebra," Princeton University Press, Princeton, New Jersey, 1956.

%
%

%

%

%

%

%

\bibitem{Gre} {\sc E. L. Green}, ``Graphs with relations, coverings and group-graded algebras," Trans. Amer. Math. Soc. 279 (1983) 297-310.

%


\bibitem{GMV} {\sc E.L. Green and R. Martinez-Villa}, ``Koszul and Yoneda algebras," Canad. Math. Soc. Conf. Proc. 18 (1996)  247-298.


\bibitem{GMV2} {\sc E.L. Green and R. Martinez-Villa}, ``Koszul and Yoneda algebras II," Canad. Math. Soc. Conf. Proc. 24 (1998)  227-244.




\bibitem{GrZ} {\sc E. L. Green and D. Zhacharia}, ``The cohomology ring of a monomial algebra,"  Manuscripta Math. 85 (1994) 11-23.


\bibitem{GKM} {\sc M. Goresky, R. Kottwitz and R. MacPherson}, ``Equivariant cohomology, Koszul duality, and the localization theorem," Invent. Math. 131 (1998) 25-83.


\bibitem{Gri} {\sc P.-P. Grivel},``Cat\'egories d\'eriv\'ees et foncteurs d\'eriv\'es," in {\it Algebraic D-modules}, Perspective in Mathematics 2 (Academic Press Inc., Boston, 1987) 1-108.

%
\bibitem{Ha1} {\sc D. Happel}, ``On the derived category of a finite-dimensional algebra," Comment. Math. Helv. 62 (1987) 339-389.
%



%
%


\bibitem{VHW} {\sc H.-J. Von H\"{o}hne and J. Waschb\"{u}sch}, ``Die struktur
n-reihiger Algebren," Comm. Algebra 12 (1984) 1187-1206.

\bibitem{ILP} {\sc K. Igusa, S. Liu and C. Paquette,} ``A proof of the strong no loop conjecture," Adv. Math. 228 (2011) 2731-2742.

%




\bibitem{LeC} {\sc J. Le and X. Chen}, ``Karoubianness of a triangulated category," J. Algebra 310 (2007) 452-457.

%
%



%


\bibitem{LLi} {\sc Z. Lin and S. Liu}, ``Representation theory of graded algebras given by locally finite quivers," arXiv: 2409.20392.





\bibitem{LiN} {\sc S. Liu and H. Niu}, ``Almost split sequences in tri-exact categories," J. Pure Appl. Algebra 226 (2022) 1-31.

%

%


\bibitem{Mac} {\sc S. MacLane}, ``Homology," Springer-Verlag, Berline, Heidelberg, New York, 1971.




\bibitem{Mar} {\sc R. Martinez-Villa,} ``Introduction to Koszul algebras," Rev. Un. Mat. Argentina 48 (2007) 67-95.


\bibitem{RMV1} {\sc R. Martinez-Villa,} ``Applications of Koszul algebras: The preprojective algebra," Canad. Math. Soc. Conf. Proc. 18 (1996)  487-504.


\bibitem{RMV2} {\sc R. Martinez-Villa,} ``Graded, self-injective and Koszul algebras," J. Algebra 215 (1999) 34-72.


\bibitem{MvS} {\sc R. Martinez-Villa and M. Saorin}, ``Koszul equivalences and dualities," Pacific J. Math. 214 (2004) 359-378.




\bibitem{May} {\sc J. P. May}, ``The cohomology of restricted Lie algebras and of Hopf algebras," J. Algebra 3 (1966) 123-146.


\bibitem{MOS} {\sc V. Mazorchuk, S. Ovsienko, and C. Stroppel}, ``Quadratic duals, Koszul dual functors, and applications," Trans. Amer. Math. Soc. 361 (2009) 1129-1172.


\bibitem{Mil} {\sc D. \hspace{-3pt} Milicic}, \hspace{-3pt} ``Lectures on derived categories," \hspace{-3pt} www.math.utah.edu/\~{}milicic/Eprints/dercat.pdf.

\bibitem{NCF2} {\sc C. Nastasescu and F. van Oystaeyen,} ``Methods of Graded Rings,'' Lecture Notes in Mathematics 1836 (Springer-Verlag, Berlin, Heidelberg, 2004).

%


\bibitem{Pr2} {\sc S. B. Priddy}, ``Koszul resolutions," Trans. Amer. Math. Soc. 152 (1970) 39-60.


\bibitem{RVDB} {\sc I. Reiten and M. Van den Bergh}, ``Notherian hereditary categories satisfying Serre duality," J.
Amer. Math. Soc. 15 (2002) 295-366.


\bibitem{Ric} {\sc J. Rickard,} ``Morita theory for derived categories," J. London Math. Soc. 39 (1989) 436-456.

%

%

%

%


\bibitem{RH} {\sc S. Ryom-Hansen}, ``Koszul duality of translation and Zuckermann functor," J. Lie Theory 14 (2004) 151-163.

%



%


\bibitem{Wei} {\sc C. A. Weibel}, ``An introduction to homological algebra," Cambridge Studies in Advanced Mathematics 38 (Cambridge University Press, Cambridge, 1994).

%

\vspace{-8pt}

\end{thebibliography}
\end{document}